\documentclass[a4,10pt]{article}
\usepackage{amsmath,amscd,amssymb}
\usepackage[margin=.8in]{geometry}
\usepackage{xcolor}

\newcommand{\nbigf}{\mathcal{F}}
\newcommand{\nbigg}{\mathcal{G}}
\newcommand{\nbigh}{\mathcal{H}}

\newcommand{\nbigl}{\mathcal{L}}

\newcommand{\nbigo}{\mathcal{O}}
\newcommand{\nbigp}{\mathcal{P}}

\newcommand{\nbigu}{\mathcal{U}}
\newcommand{\nbigv}{\mathcal{V}}

\newcommand{\seisuu}{{\mathbb Z}}
\newcommand{\rnum}{{\mathbb Q}}

\newcommand{\cnum}{{\mathbb C}}
\newcommand{\real}{{\mathbb R}}
\newcommand{\hyperh}{\mathbb{H}}
\newcommand{\hyperk}{\mathbb{K}}

\newcommand{\DD}{\mathbb{D}}

\newcommand{\gbigk}{\mathfrak K}

\newcommand{\gbigs}{\mathfrak S}

\newcommand{\gminia}{\mathfrak a}

\newcommand{\gminim}{\mathfrak m}

\newcommand{\vece}{{\boldsymbol e}}

\newcommand{\vecv}{{\boldsymbol v}}

\newcommand{\vecalpha}{{\boldsymbol \alpha}}
\newcommand{\veca}{{\boldsymbol a}}
\newcommand{\vecb}{{\boldsymbol b}}
\newcommand{\vecbeta}{{\boldsymbol \beta}}
\newcommand{\vecdelta}{{\boldsymbol \delta}}

\newcommand{\vecm}{{\boldsymbol m}}

\newcommand{\vecx}{{\boldsymbol x}}

\newcommand{\vecn}{{\boldsymbol n}}
\newcommand{\vecp}{{\boldsymbol p}}

\newcommand{\lrarr}{\longrightarrow}

\newcommand{\pf}{{\bf Proof}\hspace{.1in}}
\newcommand{\qed}{\mbox{\rule{1.2mm}{3mm}}}

\def\Image{\mathop{\rm Im}\nolimits}

\def\Gr{\mathop{\rm Gr}\nolimits}

\def\GL{\mathop{\rm GL}\nolimits}
\def\SL{\mathop{\rm SL}\nolimits}

\def\rank{\mathop{\rm rank}\nolimits}

\def\Ker{\mathop{\rm Ker}\nolimits}

\def\Gr{\mathop{\rm Gr}\nolimits}
\def\Sym{\mathop{\rm Sym}\nolimits}

\def\Res{\mathop{\rm Res}\nolimits}

\def\ch{\mathop{ch}\nolimits}

\def\tr{\mathop{\rm tr}\nolimits}
\def\Tr{\mathop{\rm Tr}\nolimits}

\def\can{\mathop{\rm can}\nolimits}

\def\id{\mathop{\rm id}\nolimits}

\def\ch{\mathop{\rm ch}\nolimits}

\newcommand{\del}{\partial}
\newcommand{\delbar}{\overline{\del}}

\newcommand{\nhom}{{\mathcal Hom}}

\newcommand{\baralpha}{\overline{\alpha}}
\newcommand{\alphabar}{\baralpha}

\newcommand{\Hbar}{\overline{H}}

\newcommand{\Sp}{{\mathcal Sp}}

\newcommand{\lefttop}[1]{{}^{#1}\!}

\def\Harm{\mathop{\rm Harm}\nolimits}

\newcommand{\Etilde}{\widetilde{E}}

\newcommand{\thetatilde}{\widetilde{\theta}}

\newcommand{\Vtilde}{\widetilde{V}}

\newcommand{\ftilde}{\widetilde{f}}
\newcommand{\utilde}{\widetilde{u}}
\newcommand{\Ftilde}{\widetilde{F}}

\newcommand{\nbigvtilde}{\widetilde{\nbigv}}

\newcommand{\Dtilde}{\widetilde{D}}
\newcommand{\Xtilde}{\widetilde{X}}
\newcommand{\vecy}{\boldsymbol y}

\def\Herm{\mathop{\rm Herm}\nolimits}

\def\Gal{\mathop{\rm Gal}\nolimits}

\def\exchange{\mathop{\rm ex}\nolimits}

\def\Disc{\mathop{\rm Disc}\nolimits}

\newcommand{\Pbar}{\overline{P}}

\newcommand{\varphitilde}{\widetilde{\varphi}}
\newcommand{\Ctilde}{\widetilde{C}}

\newcommand{\nbigftilde}{\widetilde{\nbigf}}

\newcommand{\Ytilde}{\widetilde{Y}}

\newcommand{\Qtilde}{\widetilde{Q}}

\newcommand{\betabar}{\overline{\beta}}

\newcommand{\ctilde}{\widetilde{c}}

\newcommand{\Kbar}{\overline{K}}

\newcommand{\Xbar}{\overline{X}}

\newcommand{\qtilde}{\widetilde{q}}

\newcommand{\vecq}{\boldsymbol q}

\newcommand{\Sigmatilde}{\widetilde{\Sigma}}

\newtheorem{thm}{Theorem}[section]
\newtheorem{cor}[thm]{Corollary}
\newtheorem{rem}[thm]{Remark}
\newtheorem{lem}[thm]{Lemma}
\newtheorem{prop}[thm]{Proposition}
\newtheorem{df}[thm]{Definition}

\begin{document}

\title{Harmonic metrics of generically regular semisimple Higgs bundles
on non-compact Riemann surfaces}

\author{Qiongling Li\thanks{Chern Institute of Mathematics and LPMC, Nankai University, Tianjin 300071, China, qiongling.li@nankai.edu.cn}
\and Takuro Mochizuki\thanks{Research Institute for Mathematical Sciences, Kyoto University, Kyoto 606-8512, Japan, takuro@kurims.kyoto-u.ac.jp}}
\date{}
\maketitle

\begin{abstract}
We prove that a generically regular semisimple Higgs bundle equipped
with a non-degenerate symmetric pairing on any Riemann surface
always has a harmonic metric compatible with the pairing.
We also study the classification of
such compatible harmonic metrics
in the case where the Riemann surface is
the complement of a finite set $D$ in a compact Riemann surface.
In particular,
we prove the uniqueness of a compatible harmonic metric
if the Higgs bundle is
wild and regular semisimple at each point of $D$.

\vspace{.1in}
\end{abstract}

\section{Introduction}

\subsection{Harmonic bundles}

Let $X$ be a Riemann surface.
Let $(E,\delbar_E,\theta)$ be a Higgs bundle on $X$.
Let $h$ be a Hermitian metric of $E$.
We obtain the Chern connection $\nabla_h=\delbar_E+\del_{E,h}$
and the adjoint $\theta^{\dagger}_h$ of $\theta$.
The metric $h$ is called a harmonic metric of
the Higgs bundle $(E,\delbar_E,\theta)$
if $\nabla_h+\theta+\theta^{\dagger}_h$ is flat,
i.e.,
$\nabla_h\circ\nabla_h+[\theta,\theta^{\dagger}_h]=0$,
and $(E,\delbar_E,\theta,h)$ is called a harmonic bundle.
It was introduced by Hitchin \cite{Hitchin-self-duality},
and it has been one of the most important
and interesting mathematical objects.

A starting point is the study 
of the existence and the classification of harmonic metrics.
If $X$ is compact,
the following definitive theorem is most fundamental,
which is due to Hitchin \cite{Hitchin-self-duality}
and Simpson \cite{s1}.
\begin{thm}[\mbox{\cite{Hitchin-self-duality, s1}}]
\label{thm;22.9.6.30}
If $X$ is compact,
$(E,\delbar_E,\theta)$ has a harmonic metric
if and only if
$(E,\delbar_E,\theta)$ is polystable of degree $0$.
If $h_1$ and $h_2$ are two harmonic metrics of
$(E,\delbar_E,\theta)$,
there exists a decomposition
$(E,\delbar_E,\theta)
=\bigoplus_{i=1}^m (E_i,\delbar_{E_i},\theta_i)$
such that (i)  the decomposition is orthogonal
with respect to both $h_1$ and $h_2$,
(ii) $h_{2|E_i}=a_i h_{1|E_i}$ for some $a_i>0$. 
\hfill\qed
\end{thm}

The study in the non-compact case was pioneered by Simpson \cite{s1,s2},
and pursued by Biquard-Boalch \cite{Biquard-Boalch}
and the second author \cite{Mochizuki-KH-Higgs}.
Let $X$ be the complement of a finite subset $D$
in a compact Riemann surface $\Xbar$.
A Higgs bundle
$(E,\delbar_E,\theta)$ on $X$
induces a coherent sheaf on the cotangent bundle $T^{\ast}X$
whose support $\Sigma_{E,\theta}$ is called the spectral curve of
the Higgs bundle.
The natural morphism $\Sigma_{E,\theta}\to X$ is finite and flat.
The Higgs bundle is called tame
if the closure of $\Sigma_{E,\theta}$
in $T^{\ast}\Xbar(\log D)$ is proper over $\Xbar$.
The Higgs bundle is called wild
if the closure of $\Sigma_{E,\theta}$
in the projective completion of
$T^{\ast}\Xbar$ is complex analytic.
If a tame (resp. wild) Higgs bundle $(E,\delbar_E,\theta)$
is equipped with a harmonic metric $h$,
$(E,\delbar_E,\theta,h)$
is called a tame (resp. wild) harmonic bundle.
From a tame (resp. wild) harmonic bundle
$(E,\delbar_E,\theta,h)$,
we obtain a regular (resp. good) filtered Higgs bundle
$(\nbigp^h_{\ast}(E),\theta)$ on $(\Xbar,D)$.
(See \S\ref{subsection;22.9.23.10}
for the notions of regular and good filtered Higgs bundles.
See \cite[\S2.5]{Mochizuki-KH-Higgs}
for the notation $\nbigp^h_{\ast}(E)$.)
The following theorem was proved by Simpson \cite{s2}
in the tame case
and generalized to the wild case
in \cite{Biquard-Boalch, Mochizuki-KH-Higgs}.
\begin{thm}
For a wild harmonic bundle $(E,\delbar_E,\theta,h)$ on $(\Xbar,D)$,
$(\nbigp^h_{\ast}(E),\theta)$ is polystable of degree $0$.
Conversely, for any polystable good filtered Higgs bundle
$(\nbigp_{\ast}\nbigv,\theta)$ of degree $0$ on $(X,D)$
such that $(\nbigv,\theta)_{|X\setminus D}=(E,\delbar_E,\theta)$,
there exists a harmonic metric $h$ of $(E,\delbar_E,\theta)$
such that $\nbigp^h_{\ast}(E)=\nbigp_{\ast}\nbigv$.
Moreover, if $h_i$ $(i=1,2)$ be two harmonic metrics of
$(E,\delbar_E,\theta)$ such that
$\nbigp^{h_1}_{\ast}(E)=\nbigp^{h_2}_{\ast}(E)$,
then there exists a decomposition
$(E,\delbar_E,\theta)
=\bigoplus_{i=1}^m (E_i,\delbar_{E_i},\theta_i)$
as in Theorem {\rm\ref{thm;22.9.6.30}}.
\hfill\qed
\end{thm}

More recently,
in \cite{Li-Mochizuki1, Li-Mochizuki2},
we studied a different type of existence
and classification results.
Let $X$ be any Riemann surface,
which is not necessarily the complement of
a finite subset in a compact Riemann surface.
Let $K_X$ denote the canonical bundle of $X$.
Let $r\geq 2$ be an integer.
We set
$\hyperk_{X,r}:=\bigoplus_{i=1}^rK_X^{(r+1-2i)/2}$.
For any $r$-differential $q_r$ on $X$,
we obtain the cyclic Higgs field
$\theta(q_r)$ of $\hyperk_{X,r}$
from $q_r:K_X^{-(r-1)/2}\to K_X^{(r-1)/2}\otimes K_X$
and the identity morphisms
$K_X^{(r+1-2i)/2}
\simeq
K_X^{(r+1-2(i+1)/2)}\otimes K_X$.

\begin{thm}[\mbox{\cite{Li-Mochizuki1}}]
If $q_r\neq 0$,
$(\hyperk_{X,r},\theta(q_r))$ has a harmonic metric. 
More precisely,
there exists a unique harmonic metric $h^c$
 such that
(i) $\det(h^c)=1$,
(ii) the decomposition
$\hyperk_{X,r}:=\bigoplus_{i=1}^rK_X^{(r+1-2i)/2}$
is orthogonal,
(iii) the metrics 
 $h^c_{|K_X^{(r+1-2i)/2}}\otimes
 (h^c_{|K_X^{(r+1-2(i-1))/2}})^{-1}$
 of $TX=K_X^{-1}$ are complete.
\hfill\qed
\end{thm}

Note that in general there are many other harmonic metrics
$h$ of $(\hyperk_{X,r},\theta(q_r))$
satisfying the conditions (i) and (ii).
In \cite{Li-Mochizuki2},
we studied a classification of such harmonic metrics
under the additional assumption that
$X$ is the complement of a finite subset $D$
in a compact Riemann surface $\Xbar$,
and that $q_r$ has at worst multi-growth orders
at each point of $D$.
See the introductions of \cite{Li-Mochizuki1, Li-Mochizuki2}
for a brief review of previous studies on this type of harmonic bundles
and related subjects.

\vspace{.1in}

After \cite{Li-Mochizuki1,Li-Mochizuki2},
it seems reasonable to study
whether a Higgs bundle with an appropriate additional symmetry
has a harmonic metric compatible with the symmetry.
Note that such a harmonic metric is interesting
in relation with the theory of minimal surfaces.
If the base space is compact, it should be
a consequence of the existence and uniqueness theorem
due to Hitchin and Simpson
(see Theorem \ref{thm;22.9.6.30}).
Namely, we know the existence of a harmonic metric
of a polystable Higgs bundle of degree $0$,
and the uniqueness should imply the compatibility
of the harmonic metric with the symmetry.
(For example, see \S\ref{subsection;22.9.22.1}.)
If the base space is the complement of
a finite subset in a compact Riemann surface,
and if the Higgs bundle is wild,
we may still apply a similar argument
to obtain the classification of harmonic metrics.
However, in a more general situation,
we do not know a useful general result
for neither existence nor uniqueness.
We are motivated to find many examples.
Even in the wild case,
it would be helpful
if we could simplify the assumption on Higgs bundles.
We also note that the analysis of the non-compact case
would be useful for the compact case.
In this paper, we shall study Higgs bundles
equipped with a non-degenerate symmetric pairing
under the assumption that
the Higgs field is generically regular semisimple.

\subsection{Higgs bundles with non-degenerate symmetric product}

Let $(E,\delbar_E,\theta)$ be a Higgs bundle
on any Riemann surface $X$.
A non-degenerate symmetric pairing of
the Higgs bundle
is a holomorphic symmetric bilinear form
$C:E\otimes E\to\nbigo_X$
such that
$C(\theta\otimes\id)=C(\id\otimes\theta)$.
When $(E,\delbar_E,\theta)$ is equipped with
a non-degenerate symmetric pairing $C$,
we say that a harmonic metric $h$ is compatible with $C$
if the induced morphism
$\Psi_C:E\to E^{\lor}$
is isometric with respect to $h$ and $h^{\lor}$.
In that case,
the harmonic map $(E,\delbar_E,\theta)$ is induced
by a $\GL(r,\real)$-harmonic bundle,
where $r=\rank E$.

We shall also impose the following generic regular semisimplicity
condition.
\begin{itemize}
\item  There exists $Q\in X$
such that the fiber of $\Sigma_{E,\theta}\to X$
over $Q$ consists of exactly $r$ points.
\end{itemize}

\begin{rem}
From the definition,
it is clear that a Higgs bundle $(E,\theta)$ is
generically regular semisimple
if one of the following holds:
(i) $\tr(\theta^r)\neq 0$
and $\tr(\theta^k)=0$ $(k=1,\ldots,r-1)$;
(ii) $\tr(\theta^{r-1})\neq 0$
and $\tr(\theta^k)=0$ $(k=1,\ldots,r-2, r)$,
where $r=\rank E$.
We also remark that
most Higgs bundles
are generically regular semisimple. 
For example, a Higgs bundle $(E,\theta)$ is
generically regular semisimple 
if one of the following holds:
(i) the characteristic polynomial of $(E,\theta)$ is irreducible;
(ii) $E$ has no proper subbundle preserved by $\theta$.
See {\rm\S\ref{subsection;22.10.10.1}}
(Proposition {\rm\ref{prop;22.10.11.1}} and
Proposition {\rm\ref{prop;22.10.11.2}})
for more details.
\hfill\qed
\end{rem}

We shall prove the following general theorem.
\begin{thm}[Theorem
\ref{thm;22.8.31.10}]
\label{thm;22.9.6.100}
Suppose that $(E,\delbar_E,\theta)$ is
generically regular semisimple
and equipped
with a non-degenerate symmetric pairing $C$.
Then, $(E,\delbar_E,\theta)$
has a harmonic metric compatible with $C$.
\end{thm}

In the compact case,
it follows from Theorem \ref{thm;22.9.6.30},
but note that we do not have to assume
that the Higgs bundle is polystable of degree $0$.
For the proof of Theorem \ref{thm;22.9.6.100}
in the non-compact case,
we use a method developed in \cite{Li-Mochizuki2}
together with a key estimate for harmonic metrics
compatible with a non-degenerate symmetric product
(Proposition \ref{prop;22.9.1.2}).
It follows from Simpson's main estimate (\cite{s2,s5})
on the norm of the Higgs field
with the aid of
a linear algebraic argument (Proposition \ref{prop;22.8.31.2}).

\begin{rem}
We can obtain another a priori estimate 
for harmonic metrics compatible with a non-degenerate symmetric pairing
from a variant of Simpson's main estimate 
{\rm\cite{s2, Decouple}}
and a linear algebraic argument.
It is useful even in the study of Higgs bundles
on compact Riemann surfaces
which are not necessarily equipped with
a global non-degenerate symmetric pairing.
In {\rm\cite{Mochizuki-Szabo}}, 
the estimate is applied to study 
large-scale solutions of Hitchin equations.
Non-degenerate symmetric pairings are also 
useful in the study of Hitchin metric
of the moduli space of Higgs bundles 
{\rm \cite{Mochizuki-Hitchin-metric}}.
\hfill\qed
\end{rem}

In particular,
suppose $(E,\delbar_E,\theta)$ is cyclic,
that is, $(E,\delbar_E)$ is 
a direct sum of holomorphic line bundles $L_k$ $(k=1,\ldots, r)$
and $\theta$ takes $L_k$ to
$L_{k+1}\otimes K_X$ $(k=1,\ldots,r)$ where $L_{r+1}:=L_1$.
Then, if $\det(\theta)\neq 0$,
$(E,\delbar_E,\theta)$ is
obviously generically regular semisimple.
We obtain the following corollary.
\begin{cor}
Suppose $(E,\delbar_E,\theta)$ is
a cyclic Higgs bundle on $X$ satisfying $\det\theta\neq 0$ 
and equipped with a non-degenerate symmetric pairing $C$.
Then, 
$(E,\delbar_E,\theta)$ has a harmonic metric compatible with $C$.
\end{cor}

If $X$ is compact,
the assumptions imply that
$(E,\delbar_E,\theta)$ is polystable of degree $0$,
and we can easily deduce Theorem \ref{thm;22.9.6.100}
from Theorem \ref{thm;22.9.6.30}.
(See \S\ref{subsection;22.9.22.1}.)
We can also obtain the uniqueness of
such compatible harmonic metrics in the compact case.
In the non-compact case,
the uniqueness does not hold, in general.

We also study more detailed classification
in the case where $X$ is the complement of
a finite subset $D$ in a compact Riemann surface $\Xbar$.
Let $(E,\delbar_E,\theta)$ be 
a generically regular semisimple Higgs bundle on $X$
which is wild at each point of $D$.
Suppose that $(E,\delbar_E,\theta)$ is equipped with
a non-degenerate symmetric pairing $C$.

\begin{thm}[Theorem
\ref{thm;22.9.5.40}]
\label{thm;22.9.6.101}
If $h$ is a harmonic metric of $(E,\delbar_E,\theta)$
compatible with $C$,
then the induced symmetric pairing
$C$ of $\nbigp^h_{\ast}(E)$ is perfect.
Conversely, 
for any good filtered Higgs bundle
$(\nbigp_{\ast}\nbigv,\theta)$ on $(\Xbar,D)$
equipped with a perfect pairing $C_{\nbigv}$
such that $(\nbigv,\theta,C_{\nbigv})_{|X}=(E,\delbar_E,\theta,C)$,
there exists a unique harmonic metric $h$
of $(E,\delbar_E,\theta)$ compatible with $C$
satisfying
the boundary condition $\nbigp^h_{\ast}(E)=\nbigp_{\ast}\nbigv$.
\end{thm}
We remark that in Theorem \ref{thm;22.9.6.101},
we do not have to assume that
$(\nbigp_{\ast}\nbigv,\theta)$ is polystable of degree $0$
because it follows from the existence of a perfect pairing.

Let us mention a condition under which
we can prove the uniqueness of a compatible harmonic metric
without the boundary condition.
For each $P\in D$,
let $(U_P,z)$ be a holomorphic coordinate neighbourhood around $P$
such that $z(P)=0$.
We can regard $U_P$ as an open subset of $\cnum$.
For $e\in\seisuu_{>0}$,
let $\varphi_e:\cnum\to\cnum$ be the map
defined by $\varphi_e(\zeta)=\zeta^e$.
We set $U^{(e)}_P:=\varphi_e^{-1}(U_P)$.
We obtain the endomorphism
$f^{(e)}_P$ of
$\varphi_e^{\ast}(E)_{|U^{(e)}_P\setminus\{0\}}$
defined by
$f^{(e)}_P\,d\zeta/\zeta=\varphi_e^{\ast}(\theta)$.
If $U$ is sufficiently small,
and if $e=r!$,
there exist holomorphic functions
$\alpha_1,\ldots,\alpha_{r}$ on $U^{(e)}\setminus\{0\}$
and the eigen decomposition
\[
 \varphi_e^{\ast}(E)
 =\bigoplus_{i=1}^r E_{i}
\]
such that
(i) $\alpha_i-\alpha_j$ $(i\neq j)$ are nowhere vanishing
on $U^{(e)}\setminus\{0\}$,
(ii) $f^{(e)}-\alpha_i\id_{E_{i}}$
is $0$ on $E_{\alpha_i}$,
(iii) $\rank E_{i}=1$.
We say that $(E,\delbar_E,\theta)$ is regular semisimple at $P$
if $(\alpha_i-\alpha_j)^{-1}$ $(i\neq j)$
are holomorphic at $0$.

\begin{thm}[Theorem \ref{thm;22.9.5.100}]
\label{thm;22.9.6.120}
Suppose $(E,\delbar_E, \theta)$ is generically regular semisimple on $X$ which is wild at each point of $D$ and equipped with a non-degenerate symmetric pairing $C$. If $(E,\delbar_E,\theta)$ is
regular semisimple at each point of $D$, there exists a unique 
 harmonic metric on
$(E,\delbar_E,\theta)$ compatible with $C$.
\hfill\qed
\end{thm}

\subsection{Hitchin section}

Our main examples are Higgs bundles
contained in the Hitchin section \cite{Hitchin-section},
which we recall by following \cite{Li-survey}.
Let $q_j$ $(j=2,\ldots,r)$ be
holomorphic $j$-differentials on $X$.
The multiplication of $q_j$ induces the following morphisms:
\[
 K_X^{(r-2i+1)/2}\to
K_X^{(r-2i+2(j-1)+1)/2}\otimes K_X
 \quad
 (j\leq i\leq r).
\]
We also have the multiplications of $i(r-i)/2$
for $i=1,\ldots,r-1$:
\[
 K_X^{(r-2i+1)/2}
 \to
 K_X^{(r-2(i+1)+1)/2}\otimes K_X.
\]
They define a Higgs field
$\theta(\vecq)$ of $\hyperk_{X,r}$.
The natural pairings
$K_X^{(r-2i+1)/2}
\otimes
K_X^{-(r-2i+1)/2}
\to \nbigo_X$
induce a non-degenerate symmetric bilinear form
$C_{\hyperk,X,r}$ of
$\hyperk_{X,r}$.
It is a non-degenerate symmetric pairing of
$(\hyperk_{X,r},\theta(\vecq))$.
The following result is a partial affirmative answer
to \cite[Question A]{Tamburelli-Wolf}.

\begin{cor}[Corollary
\ref{cor;22.9.6.102}]
If the Higgs bundle $(\hyperk_{X,r},\theta(\vecq))$
is generically regular semisimple,
then there exists a harmonic metric
of $(\hyperk_{X,r},\theta(\vecq))$
which is compatible with $C_{\hyperk,X,r}$.
It is induced by an $\SL(r,\real)$-harmonic bundle.
\hfill\qed
\end{cor}

In the other extreme case $q_i=0$ $(i=2,\ldots,r)$, 
$(\hyperk_{X,r},\theta(\vecq))$
has a harmonic metric
if and only if $X$ is hyperbolic.
For example, 
$(\hyperk_{\cnum,r},\theta(0,\ldots,0))$
does not have a harmonic metric.
(See \cite[Lemma 3.13]{Li-Mochizuki1}.)

See \S\ref{subsection;22.9.6.110}--\S\ref{subsection;22.9.6.111}
for more concrete examples.
In \S\ref{subsection;22.9.6.110},
we explain an example for which
Theorem \ref{thm;22.9.6.120} ensures the uniqueness of harmonic metrics.
In \S\ref{subsection;22.9.6.112},
we explain an example for which
the uniqueness does not hold.

\begin{rem}
In {\rm\cite{Li-Mochizuki4}},
by using different techniques,
we shall establish that the Higgs bundle $(\hyperk_{X,r},\theta(\vecq))$
has a harmonic metric for any $\vecq$ if $X$ is hyperbolic.
We shall also study the existence of harmonic metrics
of Higgs bundles in
the Collier section and Gothen section
by using the techniques in both {\rm\cite{Li-Mochizuki4}} and this paper.
\hfill\qed
\end{rem}

\subsection{Higher dimensional case}

Though we do not study the higher dimensional case in this paper,
some part of the theory can be easily generalized.
We give only a brief sketch of the statements without proof.
More details will be explained elsewhere.

Let $(E,\delbar_E,\theta)$ be a Higgs bundle of rank $r$ on $X$.
It naturally induces a coherent $\nbigo_{T^{\ast}X}$-module
whose support $\Sigma_{E,\theta}$ is called the spectral variety.
The Higgs bundle $(E,\delbar_E,\theta)$ is called
generically regular semisimple
if there exists $Q\in X$ such that
the fiber of $\Sigma_{E,\theta}\to X$ over $Q$
consists of exactly $r$ points.
A non-degenerate symmetric pairing $C$ of $E$
is called a non-degenerate symmetric pairing of $(E,\delbar_E,\theta)$
if $\theta$ is self-adjoint with respect to $C$.
Note that the existence of $C$
implies that $c_1(E)=0$ in $H^2(X,\rnum)$.

Let us consider
a generically regular semisimple Higgs bundle
$(E,\delbar_E,\theta)$ on $X$
equipped with a non-degenerate symmetric pairing $C$.
If $X$ is compact K\"ahler,
we can prove that
$(E,\delbar_E,\theta)$ is polystable.
By the theory of Simpson {\rm\cite{s1}},
we can prove that there exists
a unique Hermitian-Einstein metric $h$ of $(E,\delbar_E,\theta)$
which is compatible with $C$. 
If moreover $\ch_2(E)=0$, $h$ is pluri-harmonic.
In this sense,
Theorem {\rm\ref{thm;22.9.6.100}}
can be generalized if $X$ is compact K\"ahler.

Suppose that there exists a complex projective manifold $\Xbar$
with a simple normal crossing hypersurface $D$
such that $X=\Xbar\setminus D$.
If $(E,\delbar_E,\theta)$ is good wild along $D$
and equipped with a pluri-harmonic metric compatible with $C$,
then $C$ induces a perfect pairing of
the associated good filtered Higgs bundle
$(\nbigp^h_{\ast}(E),\theta)$.
Conversely, suppose that $(E,\delbar_E,\theta)$ and $C$
extends to 
a good filtered Higgs bundle
$(\nbigp_{\ast}\nbigv,\theta)$ with a perfect pairing $C_{\nbigv}$
such that $\ch_2(\nbigp_{\ast}\nbigv)=0$,
then there exists a unique pluri-harmonic metric
$h$ of $(E,\delbar_E,\theta)$ compatible with $C$
such that $\nbigp^h_{\ast}E=\nbigp_{\ast}\nbigv$.
In this sense, Theorem {\rm\ref{thm;22.9.6.101}}
can be generalized.

We can generalize the regular semisimplicity condition
for a good wild Higgs bundle $(E,\delbar_E,\theta)$ on $(X,D)$
at $P\in D$.
If it is satisfied at any point of $D$,
there exists a unique good filtered extension
$(\nbigp_{\ast}\nbigv,\theta)$ of $(E,\delbar_E,\theta)$
such that $C$ induces a perfect pairing of
$(\nbigp_{\ast}\nbigv,\theta)$.
If $\ch_2(\nbigp_{\ast}\nbigv)=0$,
then $(E,\delbar_E,\theta)$ has a unique pluri-harmonic
metric compatible with $C$.
If $\ch_2(\nbigp_{\ast}\nbigv)\neq 0$,
then $(E,\delbar_E,\theta)$ does not have
a pluri-harmonic metric compatible with $C$.
In this sense, Theorem \ref{thm;22.9.6.120} can be generalized.

\section{Harmonic bundles with real structure}

\subsection{Real structures and symmetric pairings of vector spaces}

\subsubsection{Real structures of vector spaces}

For complex vector spaces $V_i$  $(i=1,2)$,
an $\real$-linear morphism $f:V_1\to V_2$
is called sesqui-linear
if $f(\alpha v)=\alphabar f(v)$
for any $\alpha\in\cnum$ and $v\in V_1$.

Let $V$ be a finite dimensional complex vector space.
A real structure of $V$ is
a sesqui-linear morphism $\kappa:V\to V$
such that $\kappa\circ\kappa=\id_V$.
When a real structure $\kappa$ is provided with $V$,
we set
$V_{\real,\kappa}:=\{v\in V\,|\,\kappa(v)=v\}$.
It is naturally an $\real$-vector space.
The naturally induced morphism
$\cnum\otimes V_{\real,\kappa}\to V$ is an isomorphism
of complex vector spaces.
We have
$\kappa(\alpha\otimes v_{\real})=\alphabar\otimes v_{\real}$
under the identification.

Let $C_{\real,\kappa}$ be a positive definite
symmetric bilinear form on $V_{\real,\kappa}$.
It extends to a symmetric bilinear form  $C$
and a Hermitian metric $h$
on $V$
in the standard ways:
\begin{equation}
\label{eq;22.8.29.11}
 C(\alpha\otimes u_{\real},\beta\otimes v_{\real})
 =\alpha\beta C_{\real,\kappa}(u_{\real},v_{\real}),
\quad\quad
 h(\alpha\otimes u_{\real},\beta\otimes v_{\real})
=\alpha\betabar C_{\real,\kappa}(u_{\real},v_{\real}).
\end{equation}
Here, $\alpha,\beta\in\cnum$
and $u_{\real},v_{\real}\in V_{\real,\kappa}$.
We have
$C(\kappa(u),\kappa(v))=\overline{C(u,v)}$
and
$h(\kappa(u),\kappa(v))
=\overline{h(u,v)}=h(v,u)$.

\subsubsection{Compatibility of symmetric pairings and Hermitian metrics}

Let us recall the notion of compatibility of
a non-degenerate symmetric pairing
and a Hermitian metric
on a complex vector space $V$.
Let $V^{\lor}$ denote the dual space of $V$.
Let $\langle\cdot,\cdot\rangle:V^{\lor}\times V\to \cnum$
denote the canonical pairing.

Let $C:V\times V\to\cnum$
be a non-degenerate symmetric bilinear form.
We obtain the linear isomorphism
$\Psi_C:V\simeq V^{\lor}$
by
$\langle \Psi_C(u),v\rangle=C(u,v)$.
We obtain the symmetric bilinear form
$C^{\lor}:V^{\lor}\times V^{\lor}\to\cnum$
by
\[
 C^{\lor}(u^{\lor},v^{\lor})
 =C(\Psi_C^{-1}(u^{\lor}),\Psi_C^{-1}(v^{\lor})).
\]
We have $\Psi_{C^{\lor}}\circ\Psi_C=\id_V$.

Let $h$ be a Hermitian metric of $V$.
We obtain the sesqui-linear isomorphism
$\Psi_h:V\simeq V^{\lor}$
by
$\langle\Psi_h(u),v\rangle=h(v,u)$.
We obtain the Hermitian metric $h^{\lor}$ of $V^{\lor}$
by
\[
h^{\lor}(u^{\lor},v^{\lor})=
h\bigl(
 \Psi_h^{-1}(v^{\lor}),\Psi_h^{-1}(u^{\lor})
\bigr).
\]
It is easy to see that
$\Psi_{h^{\lor}}\circ\Psi_h=\id_V$.
\begin{df}
We say that $h$ is compatible with $C$
if $\Psi_C$ is isometric with respect to $h$ and $h^{\lor}$.
Let $\Herm(V,C)$ denote the space of Hermitian metrics
of $V$ which are compatible with $C$.
 \hfill\qed
\end{df}

\begin{lem}
\label{lem;22.8.29.10}
The following conditions are equivalent.
\begin{itemize}
 \item  $h$ is compatible with $C$.
 \item $C(u,v)=\overline{C^{\lor}(\Psi_h(u),\Psi_h(v))}$ holds
       for any $u,v\in V$.
 \item $\Psi_{C^{\lor}}\circ\Psi_h=\Psi_{h^{\lor}}\circ\Psi_C$ holds.
It is also equivalent to
$\Psi_{C}\circ\Psi_{h^{\lor}}
=\Psi_h\circ\Psi_{C^{\lor}}$.
\end{itemize}
\end{lem}
\pf
Let $f\in V^{\lor}$ and $u\in V$.
We have
\[
 \langle f,\Psi_{C^{\lor}}\circ\Psi_h(u)\rangle
 =C^{\lor}(f,\Psi_h(u))
 =\langle\Psi_h(u),\Psi_{C}^{\lor}(f)\rangle
 =h(\Psi_C^{-1}(f),u).
\]
We also have
\[
\langle f,\Psi_{h^{\lor}}\circ\Psi_C(u)\rangle
=h^{\lor}(f,\Psi_C(u))
=\overline{h^{\lor}(\Psi_C(u),f)}
=\overline{\langle
\Psi_{h^{\lor}}(f),\Psi_C(u)
 \rangle}
=\overline{C(u,\Psi_h^{-1}(f))}.
\]
Hence, we obtain the claim of Lemma \ref{lem;22.8.29.10}.
\hfill\qed

\vspace{.1in}
From a non-degenerate symmetric bilinear form $C$
and a Hermitian metric $h$,
we obtain a sesqui-linear isomorphism
$\kappa=\Psi_{C^{\lor}}\circ\Psi_h$.
\begin{lem}
$h$ is compatible with $C$
if and only if $\kappa$ is a real structure of $V$.
\end{lem}
\pf
We have
$\kappa^{-1}=\Psi_h^{-1}\circ\Psi_{C^{\lor}}^{-1}
=\Psi_{h^{\lor}}\circ \Psi_C$.
Hence, the claim follows from Lemma \ref{lem;22.8.29.10}.
\hfill\qed

\begin{lem}
\label{lem;22.9.3.2}
Suppose that $h$ is compatible with $C$.
\begin{itemize}
 \item 
 We have $h(\kappa(u),\kappa(v))=\overline{h(u,v)}=h(v,u)$
and $h(u,v)=C(u,\kappa(v))$
for any $u,v\in V$. 
 \item
      For any $u,v\in V_{\real,\kappa}$,
we have $C(u,v)=h(u,v)\in\real$.
 \item
      Let $C_{\real,\kappa}$ denote
      the $\real$-valued symmetric bilinear form on
      $V_{\real,\kappa}$
      obtained as the restriction of $C$.
      Then,
      it is also equal to the restriction of $h$ to $V_{\real,\kappa}$.
      Moreover,
      $h$ and $C$ are related with $C_{\real,\kappa}$
      as in {\rm(\ref{eq;22.8.29.11})}.            
\end{itemize}
\end{lem}
\pf
By the constructions, we have
\[
 h(\kappa(u),\kappa(v))
=h\bigl(\Psi_C^{-1}\Psi_h(u),\Psi_C^{-1}\Psi_h(v)\bigr)
=h^{\lor}(\Psi_{h}(u),\Psi_h(v))
=h(v,u).
\]
We also have
\[
 C(u,\kappa(v))
=C(\kappa(v),u)
=\langle
 \Psi_C\bigl(\Psi_C^{-1}\circ\Psi_h(v)\bigr),u
 \rangle
 =\langle \Psi_h(v),u\rangle
 =h(u,v).
\]
If $u,v\in V_{\real,\kappa}$,
we obtain
$h(u,v)=h(\kappa(u),\kappa(v))
=h(v,u)=\overline{h(u,v)}$,
and hence $h(u,v)\in\real$.
We also have
$h(u,v)=C(u,\kappa(v))=C(u,v)$.
Thus, we obtain the second claim.
The third claim follows easily.
\hfill\qed

\begin{cor}
\label{cor;22.9.22.2}
$h$ is compatible with $C$
if and only if
there exists a base of $V$
which is orthonormal with respect to
both $h$ and $C$.
\end{cor}
\pf
The ``only if'' part follows from 
Lemma \ref{lem;22.9.3.2}.
The ``if'' part is easy to see.
\hfill\qed 

\subsubsection{Hermitian automorphisms}

Let $h$ be a Hermitian metric of $V$.
Let $\nbigh(V,h)$ be the set of automorphisms $f$ of $V$
which are Hermitian with respect to $h$,
i.e., $h(fu,v)=h(u,fv)$ for any $u,v\in V$.
Let $\nbigh_+(V,h)\subset\nbigh(V,h)$
be the subset of $f\in\nbigh(V,h)$
such that any eigenvalues of $f$ are positive.
There exists the exponential map
$\exp:\nbigh(V,h)\to\nbigh_+(V,h)$
defined by
$\exp(H)=\sum_{m=0}^{\infty}\frac{1}{m!}H^m$.
We also have the well defined logarithm
$\log:\nbigh_+(V,h)\to\nbigh(V,h)$
such that
$\log\circ\exp=\id_{\nbigh(V,h)}$
and $\log\circ\exp=\id_{\nbigh_+(V,h)}$.
For $H\in \nbigp_+(V,h)$,
there exist a unitary matrix $P$
and a diagonal matrix $G$
such that $H=PGP^{-1}$.
The $(i,i)$-entries $G_{i,i}$ of $G$ are positive.
We have $\log(H)=P\log(G)P^{-1}$,
where $\log(G)$ is the diagonal matrix
whose $(i,i)$-entries are $\log(G_{i,i})$.
We recall the following well known lemma.
\begin{lem}
The exponential map
$\exp:\nbigh(V,h)\to\nbigh_+(V,h)$
is a diffeomorphism.
\end{lem}
\pf
It is easy to see that
the exponential function is $C^{\infty}$.
It is enough to check that the logarithm is also $C^{\infty}$.
For $H\in\nbigh_+(V,h)$,
we have
\begin{equation}
\label{eq;22.9.7.1}
 \log(H)=
\frac{1}{2\pi\sqrt{-1}}
 \int_{\Gamma}
 \log(z)\cdot (z\id_V-H)^{-1}dz.
\end{equation}
Here, $\Gamma$ is the union of small circles 
with the counter clock-wise direction
around the eigenvalues of $H$
in $\cnum$.
Hence, $\log:\nbigh_+(V,h)\to\nbigh(V,h)$
is also $C^{\infty}$.
\hfill\qed

\vspace{.1in}

Let $C$ be a non-degenerate symmetric bilinear form on $V$.
Suppose that $h$ is compatible with $C$.
Let $\nbigh(V,C,h)$ be the set of $f\in\nbigh(V,h)$
which are anti-symmetric with respect to $C$,
i.e., $C(fu,v)+C(u,fv)=0$ for any $u,v\in V$.
Let $\nbigh_+(V,C,h)$ be the set of $f\in\nbigh_+(V,h)$
which are isometric with respect to $C$,
i.e., $C(fu,fv)=C(u,v)$ for any $u,v\in V$.

\begin{lem}
\label{lem;22.9.3.30}
The exponential map induces
a diffeomorphism
$\exp:\nbigh(V,C,h)\to \nbigh_+(V,C,h)$.
\end{lem}
\pf
We set $n:=\dim_{\cnum}V$.
Let $I_n$ denote the identity $(n\times n)$-matrix.
Let $\nbigh_0$ denote the space of
Hermitian $(n\times n)$-matrices.
Let $\nbigh_{0,+}\subset\nbigh_0$
denote the space of positive definite Hermitian $(n\times n)$-matrices.
Let $\nbigh_1\subset\nbigh_0$ denote the subspace of
$H\in\nbigh_0$ such that $H+\lefttop{t}H=0$.
Let $\nbigh_{1,+}\subset\nbigh_{0,+}$ denote the subspace of
$H\in\nbigh_0$ such that $H\cdot\lefttop{t}H=I_n$.

The exponential map induces a diffeomorphism
$\exp:\nbigh_{0}\to\nbigh_{0,+}$.
It induces a $C^{\infty}$-map
$\nbigh_{1}\to\nbigh_{1,+}$.
Let $H\in \nbigh_{1,+}$.
We note that $\Hbar=\lefttop{t}H$.
There exist
a diagonal matrix $G$ and a unitary matrix $P$
such that $H=PGP^{-1}$.
The diagonal entries of $G$ are positive numbers.
We have
\[
G\Pbar^{-1}\cdot P G=\Pbar^{-1}\cdot P.
\]
It implies that
$(\Pbar^{-1}\cdot P)_{i,j}=0$
unless $G_{i,i}\cdot G_{j,j}=1$.
Hence, we obtain
\[
 \log(G)\cdot (\Pbar^{-1}\cdot P)
+ (\Pbar^{-1}\cdot P)\cdot\log G=0.
\]
We obtain
$\Pbar(\log G)\Pbar^{-1}
+P\cdot (\log G)P^{-1}=0$,
i.e.,
$
\lefttop{t}(\log H)+\log H
=\overline{(\log H)}+\log H=0$.
Hence, the exponential map induces
a diffeomorphism
$\nbigh_1\lrarr \nbigh_{1,+}$.

Let $\vecv$ be a base of $V_{\real,\kappa}$
which is orthonormal with respect to $C_{\real,\kappa}$.
Any $f\in\nbigh_+(V,C,h)$
(resp. $f\in\nbigh(V,C,h)$)
is represented by
a matrix in $\nbigh_{1,+}$
(resp. $\nbigh_{1}$)
with respect to $\vecv$.
Hence,
the exponential map induces
a diffeomorphism
$\nbigh(V,C,h)\simeq \nbigh_+(V,C,h)$.
\hfill\qed

\vspace{.1in}

We recall that any $f\in \nbigh_+(V,h)$ has
natural $s$-powers $f^s=\exp(s\log f)$ $(s\in\real)$.
We obtain the following lemma
from Lemma \ref{lem;22.9.3.30}.

\begin{lem}
\label{lem;22.9.3.1}
If $f\in \nbigh_+(V,C,h)$,
then $f^s\in \nbigh_+(V,C,h)$ for any $s\in\real$.
\hfill\qed
\end{lem}

We also obtain the following lemma.
\begin{lem}
\label{lem;22.9.3.3}
For any $f\in\nbigh_+(V,C,h)$,
we have $\det(f)=1$.
\end{lem}
\pf
Let $\vecv$ be a base of $V$
as in the proof of Lemma \ref{lem;22.9.3.30}.
Then,
$f$ is represented by a positive definite Hermitian matrix
$H$ such that $\Hbar\cdot H=I_n$.
It implies $\det(H)^2=\det(\Hbar)\det(H)=1$,
and hence $\det(H)=1$.
\hfill\qed

\vspace{.1in}
Let $f\in\nbigh_+(V,C,h)$.
Let $V(f,a)\subset V$ denote the eigen space of $f$
corresponding to the eigenvalue $a\in\real$.
We obtain the decomposition
$V=\bigoplus_{a>0}V(f,a)$.
We set $V(f,a,a^{-1})=V(f,a)\oplus V(f,a^{-1})$
for $a>1$.

\begin{lem}
\label{lem;22.9.3.110}
The decomposition
\begin{equation}
\label{eq;22.9.3.100}
 V=V(f,1)\oplus
 \bigoplus_{a>1}V(f,a,a^{-1})
\end{equation}
is orthogonal with respect to
both $h$ and $C$.
The real structure $\kappa$ preserves
the decomposition {\rm(\ref{eq;22.9.3.100})}. 
Moreover, $\kappa$ exchanges
$V(f,a)$ and $V(f,a^{-1})$.
\end{lem}
\pf
Let $\vecv$ be a base of $V_{\real,\kappa}$
which is orthonormal with respect to $C_{\real,\kappa}$.
We obtain the Hermitian matrix $H$
representing $f$ with respect to $\vecv$.
We have $\Hbar=H^{-1}$,
i.e., $\kappa\circ f\circ \kappa=f^{-1}$.
Hence, we obtain that $\kappa$ preserves
(\ref{eq;22.9.3.100}),
and that it exchanges $V(f,a)$ and $V(f,a^{-1})$.
It implies that (\ref{eq;22.9.3.100})
is induced by a decomposition
\begin{equation}
\label{eq;22.9.3.101}
 V_{\real,\kappa}=
  V(f,1)_{\real,\kappa}
  \oplus
  \bigoplus_{a>1}V(f,a,a^{-1})_{\real,\kappa}.
\end{equation}
We note that (\ref{eq;22.9.3.100}) is orthogonal with respect to $h$
because it is induced by the eigen decomposition of
the Hermitian automorphism $f$.
It implies that
the decomposition (\ref{eq;22.9.3.101})
is orthogonal with respect to $C_{\real,\kappa}$.
Hence, 
the decomposition (\ref{eq;22.9.3.100})
is orthogonal with respect to $C$.
\hfill\qed

\subsubsection{Difference between two compatible Hermitian metrics}
\label{subsection;22.9.4.4}

Let $h$ and $h'$ be Hermitian metrics of $V$.
There exists the unique automorphism $s(h,h')$
such that
$h'(u,v)=h(s(h,h')u,v)$ for any $u,v\in V$.
Note that $s(h,h')$ is self-adjoint with respect to
both $h$ and $h'$.
Let $s(h,h')^{\lor}$ denote the automorphism of $V^{\lor}$
obtained as the dual of $s(h,h')$.
We have
$\Psi_{h'}=\Psi_{h}\circ s(h,h')=s(h,h')^{\lor}\circ\Psi_h$.
Suppose that $h$ is compatible with
a non-degenerate symmetric pairing $C$.

\begin{lem}
\label{lem;22.8.29.20}
The following conditions are equivalent.
\begin{itemize}
 \item  $h'$ is compatible with $C$.
 \item $s(h,h')$ is an isometry with respect to $C$,
       i.e., $C(s(h,h')u,s(h,h')v)=C(u,v)$ for any $u,v\in V$.
 \item $\kappa\circ s(h,h')\circ\kappa=s(h,h')^{-1}$.
 \item Let $\vecv$ be an orthonormal frame of $V_{\real,\kappa}$
       with respect to $C_{\real,\kappa}$.
       Let $A$ be the Hermitian matrix representing
       $s(h,h')$
       with respect to $\vecv$.
       Then,
       $A\cdot \lefttop{t}A$
       equals the identity matrix,
       which is equivalent to
       $A^{-1}=\overline{A}$.
\end{itemize}
\end{lem}
\pf
To simplify the description,
we set $f=s(h,h')$.
Because 
$\Psi_{C}^{-1}\circ \Psi_{h'}=
\Psi_{C}^{-1}\circ \Psi_h\circ f
=\Psi_{h}^{-1}\circ\Psi_C\circ f$
and
$\Psi_{h'}^{-1}\circ\Psi_C
=\Psi_h^{-1}\circ(f^{\lor})^{-1}\circ\Psi_C$,
$h'$ is compatible with $C$
if and only if
$\Psi_C\circ f=(f^{\lor})^{-1}\circ\Psi_C$.
The latter condition is equivalent to
that $f$ is isometry with respect to $C$.
Thus, the first condition is equivalent to the second.
Because 
$\Psi_{C}^{-1}\circ \Psi_{h'}=\kappa\circ f$
and
$\Psi_{h'}^{-1}\circ\Psi_C
=f^{-1}\circ\kappa$,
the first condition is equivalent to the third.
It is easy to see that 
the fourth condition is equivalent to
both the second and third.
\hfill\qed

\vspace{.1in}

We note that
$s(h,h')^{1/2}$ induces an isometry $(V,h')\simeq (V,h)$,
i.e.,
$h(s(h,h')^{1/2}u,s(h,h')^{1/2}v)=h'(u,v)$ for any $u,v\in V$.

\begin{lem}
\label{lem;22.9.7.2}
Suppose that $h'$ is compatible with $C$.
We set $\kappa'=\Psi_{C}^{-1}\circ\Psi_{h'}$,
which is a real structure of $V$.
Then, the following holds.
\begin{itemize}
 \item
$s(h,h')^{1/2}$ is an isometry with respect to $C$.
 \item
      $s(h,h')^{1/2}$ induces an isomorphism
      $(V,\kappa')\simeq (V,\kappa)$,
      i.e.,
      $s(h,h')^{1/2}\circ \kappa'=\kappa\circ s(h,h')^{1/2}$.
 \item
      There exists an $\real$-isomorphism
      $s(h,h')_{\real}^{1/2}:V_{\real,\kappa'}\simeq V_{\real,\kappa}$
      whose complexification equals $s(h,h')^{1/2}$.
      Moreover, $s(h,h')^{1/2}_{\real}$ is an isometry
      $(V_{\real,\kappa'},C_{\real,\kappa'})
      \simeq
      (V_{\real,\kappa},C_{\real,\kappa})$.
\end{itemize} 
\end{lem}
\pf
We obtain the first claim from 
Lemma \ref{lem;22.9.3.1}.
By Lemma \ref{lem;22.8.29.20},
we obtain $\kappa\circ s(h,h')^{1/2}\circ\kappa=s(h,h')^{-1/2}$.
Because $\kappa'=s(h,h')^{-1}\circ \kappa$,
we obtain
$\kappa\circ s(h,h')^{1/2}=s(h,h')^{-1/2}\circ\kappa
=s(h,h')^{1/2}\circ\kappa'$.
The third claim follows from the first and second.
\hfill\qed

\vspace{.1in}

Let $h'\in\Herm(V,C)$.
There exists the decomposition
$V=\bigoplus_{a>0}V_a$
such that
(i) it is orthogonal with respect to both $h'$ and $h$,
(ii) $h'_{|V_a}=a\cdot h_{|V_a}$.
For $a>1$, we set $\Vtilde_a=V_a\oplus V_{a^{-1}}$.
We obtain the decomposition
\begin{equation}
\label{eq;22.9.3.102}
 V=V_1\oplus \bigoplus_{a>1}\Vtilde_a.
\end{equation}
We obtain the following lemma from
Lemma \ref{lem;22.9.3.110}.

\begin{lem}
The decomposition {\rm(\ref{eq;22.9.3.102})}
is orthogonal with respect to $h$, $h'$ and $C$.
It is preserved by $\kappa$.
Moreover, $\kappa$ exchanges $V_a$ and $V_{a^{-1}}$.
\hfill\qed
\end{lem}

\subsubsection{Regular semisimple automorphisms}

Let $C$ be a non-degenerate symmetric bilinear form of $V$.
Let $F$ be an endomorphism of $V$
satisfying the following conditions.
\begin{itemize}
 \item  $F$ is symmetric with respect to $C$, i.e.,
	$C(F\otimes \id_V)=C(\id_V\otimes F)$.
 \item 	$F$ is regular semisimple,
	i.e.,
	the multiplicity of each eigenvalue of $F$ is $1$.
\end{itemize}
We have the eigen decomposition of $F$:
\begin{equation}
\label{eq;22.8.31.1}
 V=\bigoplus_{\alpha\in \Sp(F)} V_{\alpha}.
\end{equation}
Here, $\Sp(F)$ denotes the set of eigenvalues of $F$.
The following lemma is obvious
but useful in this study.
\begin{lem}\mbox{{}}
\label{lem;22.9.3.20}
\begin{itemize}
\item 
The decomposition {\rm(\ref{eq;22.8.31.1})} is orthogonal
with respect to $C$,
and the restriction of $C$ to each $V_{\alpha}$
is non-degenerate. 
\item
There uniquely exists a Hermitian metric $h^{\can}$ of $V$
such that
(i) $h^{\can}$ is compatible with $C$,
(ii) the decomposition {\rm(\ref{eq;22.8.31.1})}
is orthogonal.
The metric $h^{\can}$ is called 
the canonical compatible metric of $(V,C,F)$.
\hfill\qed
 \end{itemize}
\end{lem}

We set
$c_0(F):=\max\bigl\{|\alpha|\,\big|\,\alpha\in\Sp(F)\bigr\}$
and $c_1(F):=
\min\bigl\{
|\alpha-\beta|\,\big|\,
\alpha,\beta\in\Sp(F),\,\,\alpha\neq\beta
\bigr\}$.

\begin{prop}
\label{prop;22.8.31.2}
There exists a constant $B>0$
depending only on $c_i(F)$ $(i=0,1)$ and $\dim V$
such that
the following holds 
for any $h\in\Herm(V,C)$:
\[
 |s(h^{\can},h)|_{h^{\can}}
+
 |s(h^{\can},h)^{-1}|_{h^{\can}}
 \leq B
 \cdot (1+|F|_h)^{\dim V}.
\] 
 \end{prop}
\pf
We set $n:=\dim V$.
We take an ordering $\{\alpha_1,\ldots,\alpha_n\}$
of $\Sp(F)$.
Let $e_i$ be a base of $V_{\alpha_i}$
such that $h^{\can}(e_i,e_i)=1$ and $C(e_i,e_i)=1$.
By the construction,
$\vece=(e_1,\ldots,e_n)$ is an orthonormal base of $V$
with respect to both $h^{\can}$ and $C$.
Let $\vecv=(v_1,\ldots,v_n)$ be a base of $V$
which is orthonormal with respect to
both $h$ and $C$.
(See Corollary \ref{cor;22.9.22.2}.)

Let $K$ be a matrix determined by
$\vecv=\vece\,K$.
Because both $\vecv$ and $\vece$ are orthonormal
with respect to $C$,
we have $\lefttop{t}K=K^{-1}$.

For any $\ell\in\seisuu_{\geq 0}$,
let $A(F^{\ell},\vecv)$ be the matrix representing $F^{\ell}$
with respect to $\vecv$.
Let $\Gamma$ be the diagonal $(n\times n)$-matrix
whose $(i,i)$-th entries are $\alpha_i$.
We have
\[
A(F^{\ell},\vecv)=K^{-1}\Gamma^{\ell} K=\lefttop{t}K\Gamma^{\ell} K.
\]
Let $(A(F^{\ell},\vecv))_{i,j}$
denote the $(i,j)$-entry of $A(F^{\ell},\vecv)$.
We have
\[
 A(F^{\ell},\vecv)_{i,j}
 =\sum_{k=1}^{n}
 \alpha_k^{\ell} K_{k,i} K_{k,j}.
\]
Because $\vecv$ is orthonormal with respect to $h$,
there exists $B_1$ depending only on $n$
such that
\[
  |A(F^{\ell},\vecv)_{i,j}|
  \leq B_1(1+\bigl|F\bigr|_h)^{\ell}.
\]
Let $W(\alpha_1,\ldots,\alpha_n)$ be
the $(n\times n)$-matrix
whose $(i,j)$-th entry is
$\alpha_j^{i-1}$.
It is invertible,
and we obtain
\[
 K_{k,i}^2=
 \sum_{\ell=0}^{n-1}
 \Bigl(
 W(\alpha_1,\ldots,\alpha_n)^{-1}
 \Bigr)_{k,\ell}
\cdot A(F^{\ell},\vecv)_{i,i}.
\]
There exists $B_2>0$ depending only on
$n$ and $c_i(F)$ $(i=0,1)$ such that
\[
 |K_{k,i}|^2
 \leq
 B_2(1+|F|_h)^{n}.
\]

We set $H(h,\vece)_{i,j}=h(e_i,e_j)$.
We obtain the Hermitian matrix $H(h,\vece)=(H(h,\vece)_{i,j})$.
We obtain $s(h^{\can},h)\vece=\vece\cdot \lefttop{t}H(h,\vece)$.
Because $\vecv$ is orthonormal with respect to $h$,
we obtain
$H(h,\vece)
=\lefttop{t}K^{-1}\cdot\overline{K^{-1}}
=K\cdot \lefttop{t}\Kbar$.
Therefore, there exists $B_3>0$
depending only on
$n$ and $c_i(F)$ $(i=0,1)$
such that
\[
 |H(h,\vece)_{i,j}|
 \leq
  B_3(1+|F|_h)^{n}.
\]
Thus,
we obtain the claim for $|s(h^{\can},h)|_{h^{\can}}$.
Because $\lefttop{t}H(h,\vece)^{-1}=H(h,\vece)$,
we also obtain the claim for $|s(h^{\can},h)^{-1}|_{h^{\can}}$.
Thus, we obtain Proposition \ref{prop;22.8.31.2}.
\hfill\qed

\subsubsection{Vector bundles}

Let $M$ be a paracompact $C^{\infty}$-manifold.
Let $V$ be a complex vector bundle on $M$
with a non-degenerate symmetric pairing $C$.
A Hermitian metric $h$ of $V$
is called compatible with $C$
if $h_{|P}$ is compatible with $C_{|P}$
for any $P\in M$.

\begin{lem}
\label{lem;22.9.3.40}  
There exists a Hermitian metric $h$
of $V$ compatible with $C$. 
 \end{lem}
\pf
For any $P\in M$,
there exist a neighbourhood $M_P$ of $P$ around $M$
and a frame $\vecv_P$ of $V_{|M_P}$
which is orthonormal with respect to $C_{|M_P}$.
There exists a Hermitian metric of $V_{|M_P}$
which is compatible with $C_{|M_P}$.

Let $U_i\subset M$ $(i=1,2)$ be open subsets.
Suppose that there exist
Hermitian metrics $h_i$ of $V_{|U_i}$
which are compatible with $C_{|U_i}$.
By using Lemma \ref{lem;22.9.3.30}
and a partition of unity on $U_1\cup U_2$
subordinated to $\{U_1,U_2\}$,
we can construct a Hermitian metric $h$
of $V_{|U_1\cup U_2}$
which is compatible with $C_{|U_1\cup U_2}$.

There exists a locally finite open covering
$\{U_i\,|\,i=1,2,\ldots,\}$ of $X$
such that $V_{|U_i}$ has a Hermitian metric
which is compatible with $C_{|U_i}$.
Then, we can inductively prove
the existence of a Hermitian metric
of $V_{|\bigcup_{i=1}^mU_i}$
which is a compatible with $C_{|\bigcup_{i=1}^mU_i}$.
We can obtain a compatible Hermitian metric of $V$
as a limit.
\hfill\qed

\subsubsection{Appendix:
Compatibility of Hermitian metric and skew-symmetric pairing}
\label{subsection;22.9.4.2}

Let $\omega$ be
a symplectic form of a finite dimensional complex vector space $V$.
We obtain the isomorphism
$\Psi_{\omega}:V\lrarr V^{\lor}$
by $\langle \Psi_{\omega}(u),v\rangle=\omega(u,v)$.
We obtain the induced symplectic form $\omega^{\lor}$
of $V^{\lor}$
by
\[
 \omega^{\lor}(u^{\lor},v^{\lor})
 =\omega\bigl(
 \Psi_{\omega}^{-1}(u^{\lor}),
 \Psi_{\omega}^{-1}(v^{\lor})
 \bigr).
\]
We have
$\Psi_{\omega^{\lor}}\circ\Psi_{\omega}=-\id_V$.
We also have $(\omega^{\lor})^{\lor}=\omega$.
We say that a Hermitian metric $h$ is compatible with $\omega$
if $\Psi_{\omega}$ is isometric with respect to $h$ and $h^{\lor}$.
The following lemma is similar to Lemma \ref{lem;22.8.29.10}.
\begin{lem}
The following conditions are equivalent.
\begin{itemize}
 \item $h$ is compatible with $\omega$.
 \item $\Psi_{h^{\lor}}\circ\Psi_{\omega}
       =\Psi_{\omega^{\lor}}\circ\Psi_h$.
 \item $\omega(u,v)=\overline{\omega^{\lor}(\Psi_h(u),\Psi_h(v))}$
       holds for any $u,v\in V$.
\hfill\qed
\end{itemize}
\end{lem}

From a symplectic form $\omega$ and a Hermitian metric $h$,
we obtain a sesqui-linear isomorphism
$\kappa=\Psi_{h^{\lor}}\circ\Psi_{\omega}$.
\begin{lem}
$h$ is compatible with $\omega$
if and only if $\kappa\circ\kappa=-\id_V$.
If $h$ is compatible with $\omega$,
the following holds.
\begin{itemize}
 \item $h(\kappa(u),\kappa(v))=\overline{h(u,v)}$
       and
       $\omega(u,\kappa(v))=h(u,v)$
       for any $u,v\in V$.
 \hfill\qed
\end{itemize}
\end{lem}

Suppose that $h$ is compatible with $\omega$.
We have $\kappa(\sqrt{-1}u)=-\sqrt{-1}\kappa(u)$
and
$\kappa(\sqrt{-1}\kappa(\sqrt{-1}u))=-u$
for any $u\in V$.
Hence, $V$ is naturally a left quaternionic vector space.
For any $v\in V$ such that $h(v,v)=1$,
we obtain the $2$-dimensional $\cnum$-vector space
$\hyperh\cdot v\subset V$
generated by $v$ and $\kappa(v)$.
We have $\omega(v,\kappa(v))=h(v,v)=1$.
In particular,
the restriction of $\omega$ to $\hyperh\cdot v$
is a symplectic form.
We also have
$h(v,\kappa(v))=-\omega(v,\kappa(v))=0$.
For any $u\in V$,
we have
$h(u,v)=\omega(u,\kappa(v))$
and $h(u,\kappa(v))=-\omega(u,v)$.
The orthogonal complement of $\hyperh\cdot\omega$
with respect to $\omega$
is equal to the orthogonal complement of $\hyperh\cdot\omega$
with respect to $h$.
There exists a decomposition
$V=\bigoplus_{i=1}^m \hyperh\cdot v_i$
such that
(i) the decomposition is orthogonal with respect to
both $\omega$ and $h$,
(ii) each $\hyperh\cdot v_i$ has a base $v_i,\kappa(v_i)$,
(iii) $h(v_i,v_i)=\omega(v_i,\kappa(v_i))=1$
and $h(v_i,\kappa(v_i))=0$.
\begin{lem}
$h$ is compatible with $\omega$
if and only if there exists a base $\vecv$ of $V$
such that 
(i) $\vecv$ is orthonormal with respect to $h$,
(ii) $\vecv$ is symplectic with respect to $\omega$.
\hfill\qed
\end{lem}

Suppose that $\omega$ and $h$ are compatible.
Let $f\in\nbigh_+(V,h)$
such that
$f$ preserves $\omega$,
i.e., $\omega(fu,fv)=\omega(u,v)$ for any $u,v\in V$.
We obtain the eigen decomposition
$V=\bigoplus_{a>0} V(f,a)$ of $f$.
By setting $V(f,a,a^{-1})=V(f,a)\oplus V(f,a^{-1})$,
we obtain a decomposition (\ref{eq;22.9.3.100}).
The following lemma is similar to Lemma \ref{lem;22.9.3.110}.
\begin{lem}
\label{lem;22.9.4.1}
We have $\kappa\circ f=f^{-1}\circ\kappa$
and 
$\kappa\bigl(V(f,a)\bigr)=V(f,a^{-1})$.
The decomposition {\rm(\ref{eq;22.9.3.100})}
is orthogonal with respect to both
$\omega$ and $h$.
Moreover, 
$V(f,a)$ and $V(f,a^{-1})$ are
Lagrangian with respect to $\omega$,
and they are orthogonal with respect to $h$. 
\hfill\qed 
\end{lem}
We also have the converse.
Let $V=\bigoplus_{a\geq 1}U_a$
be a decomposition
which is orthogonal with respect to both
$\omega$ and $h$.
For $a>1$, let $U_{a,1}\subset U_a$ be a Lagrangian subspace
with respect to $\omega_{|U_a}$,
and we set $U_{a,2}=\kappa(U_{a,1})$,
which is the orthogonal complement of $U_{a,1}$ in $U_a$
with respect to $h$.
Let $g$ be the automorphism of $V$
defined by
\[
g=\id_{U_1}\oplus
\bigoplus_{a>1}
\bigl(
 a\id_{U_{a,1}}\oplus a^{-1}\id_{U_{a,2}}
\bigr).
\]
Then, $g$ is Hermitian with respect to $h$,
and $g$ preserves $\omega$.

Let $h'$ be another Hermitian metric of $V$
compatible with $\omega$.
Note that the automorphism $s(h,h')$ is Hermitian with respect to $h$
and that $s(h,h')$ preserves $\omega$.
We obtain the decomposition
$V=\bigoplus_{a>0} V_a$
such that
(i) the decomposition is orthogonal with respect to
both $h$ and $h'$,
(ii) $h'=ah$ on $V_a$.
We set $\Vtilde_a=V_a\oplus V_{a^{-1}}$ for $a>1$.
\begin{lem}
The decomposition
$V=V_1\oplus \bigoplus_{a>1}\Vtilde_a$
is orthogonal with respect to
$h$, $h'$ and $\omega$.
It is preserved by $\kappa$.
The subspaces $V_a$ and $V_{a^{-1}}$ of $\Vtilde_a$
are Lagrangian with respect to $\omega$,
and $\kappa$ exchanges $V_a$ and $V_{a^{-1}}$.
\hfill\qed
\end{lem}

\subsection{Harmonic bundles with real structure}

\subsubsection{$\GL(n,\real)$-harmonic bundles}

Let us recall the notion of $\GL(n,\real)$-harmonic bundle.
Let $Y$ be a Riemann surface.
Let $\pi:\Ytilde\to Y$ be a universal covering.
We fix $Q\in Y$ and $\Qtilde\in\pi^{-1}(Q)$.

Let $V_{\real}$ be an $\real$-vector bundle
of rank $n$ on $Y$ equipped with a flat connection $\nabla_{\real}$.
Let $C_{\real}$ be
a positive definite symmetric bilinear form of $V_{\real}$.
There exists a $\nabla_{\real}$-flat trivialization
$\pi^{-1}(V_{\real})\simeq \Ytilde\times\real^n$.
From $\pi^{-1}C_{\real}$,
we obtain a $\pi_1(Y,Q)$-equivariant map
$F_{C_{\real}}:\Ytilde\to \GL(n,\real)/O(n)$,
where we naturally identify
$\GL(n,\real)/O(n)$
with the space of positive definite symmetric bilinear forms.
If $F_{C_{\real}}$ is harmonic with respect to
a K\"ahler metric $g_{\Ytilde}$ of $\Ytilde$
and the natural Riemannian metric of
$\GL(n,\real)/O(n)$,
$C_{\real}$ is called a harmonic metric of
$(V_{\real},\nabla_{\real})$.
The condition is independent of
the choice of $g_{\Ytilde}$.
Such a tuple $(V_{\real},\nabla_{\real},C_{\real})$
is called a $\GL(n,\real)$-harmonic bundle.

\begin{lem}
\label{lem;22.8.29.31}
$(V_{\real},\nabla_{\real},C_{\real})$
is a $\GL(n,\real)$-harmonic bundle on $Y$
if and only if  
 $(V_{\real}\otimes\cnum,\nabla,h)$ is a harmonic bundle,
i.e., the induced map $\Ytilde\to \GL(n,\cnum)/U(n)$ is harmonic.
Here, $h$ denotes the Hermitian metric 
induced by $C_{\real}$.
\end{lem}
\pf
Because $\GL(n,\real)/O(n)\to \GL(n,\cnum)/U(n)$
is totally geodesic,
we obtain the claim of the lemma.
\hfill\qed

\subsubsection{Real structures of harmonic bundles}
\label{subsection;22.9.23.11}

Let $(E,\delbar_E,\theta,h)$ be a harmonic bundle
on the Riemann surface $Y$.
\begin{df}
\label{df;22.9.3.10}
A real structure of the harmonic bundle $(E,\delbar_E,\theta,h)$ 
is a holomorphic symmetric non-degenerate pairing $C$ of $E$
such that the following conditions are satisfied.
\begin{itemize}
 \item $h_{|Q}$ is compatible with
       $C_{|Q}$ for any $Q\in Y$.
 \item $\theta$ is self-adjoint with respect to $C$,
       i.e.,
       $C(\theta u,v)=C(u,\theta v)$
       for any local sections
       $u$ and $v$ of $E$.
\hfill\qed
\end{itemize}
\end{df}

Let $C$ be a real structure of $(E,\delbar_E,\theta,h)$.
We obtain
the holomorphic isomorphism 
$\Psi_C:E\simeq E^{\lor}$
defined by $\Psi_C(u)(v)=C(u,v)$.
We have
the sesqui-linear isomorphism
$\Psi_h:E\simeq E^{\lor}$
defined by $\Psi_h(u)(v)=h(v,u)$.
Let $\kappa:E\simeq E$
be the sesqui-linear isomorphism
defined by
$\kappa=\Psi_h^{-1}\circ\Psi_C$.
It is a real structure of a complex vector bundle $E$,
i.e., $\kappa\circ\kappa=\id_E$.
Let $E_{\real,\kappa}$ be the $\kappa$-invariant part of $E$.
There exists the natural isomorphism
$\cnum\otimes_{\real}E_{\real,\kappa}\simeq E$.
There exists a positive definite symmetric bilinear form
$C_{\real,\kappa}$ of $E_{\real}$ which induces
both $h$ and $C$.

Let $\nabla_h=\delbar_E+\del_{E,h}$
denote the Chern connection of
$(E,\delbar_E)$ with $h$.
Let $\theta_h^{\dagger}$ denote the adjoint of $\theta$
with respect to $h$.
We obtain the flat connection
$\DD^1_h=\nabla_h+\theta+\theta^{\dagger}_h$ of $E$.
\begin{lem}
\label{lem;22.8.29.30}
$\kappa$ is $\DD^1_h$-flat.
As a result,
there exists a flat connection $\nabla_{\real,\kappa}$
of $E_{\real,\kappa}$
which induces $\DD^1_h$
under the isomorphism
$E\simeq \cnum\otimes_{\real} E_{\real,\kappa}$.
\end{lem}
\pf
Because $C$ is holomorphic,
$\Psi_C\circ\delbar_E=\delbar_{E^{\lor}}\circ\Psi_C$ holds.
Because $\theta$ is self-adjoint with respect to $C$,
we have $\Psi_C\circ\theta=\theta^{\lor}\circ\Psi_C$.
Because $\Psi_C$ is an isometry with respect to $h$ and $h^{\lor}$,
we obtain
$\Psi_C\circ\del_{E,h}=\del_{E^{\lor},h^{\lor}}\circ\Psi_C$
and
$\Psi_C\circ\theta^{\dagger}_h=(\theta^{\lor})^{\dagger}_{h^{\lor}}\circ
\Psi_C$.
We have
$\Psi_h\circ\theta=(\theta_h^{\dagger})^{\lor}\circ\Psi_h$
and
$\Psi_h\circ\theta^{\dagger}_h
 =\theta^{\lor}\circ\Psi_h$.
We have
\[
 \langle\Psi_h(\delbar_{E,h}u),v\rangle
 =\del h(v,u) -h(\del_{E,h}v,u)
 =\del\langle \Psi_h(u),v\rangle
 -\langle\Psi_h(u),\del_{E,h}(v)\rangle
 =\langle \del_{E^{\lor},h^{\lor}}\Psi_h(u),v\rangle.
\]
Hence, we obtain
$\Psi_h\circ\delbar_E=\del_{E^{\lor},h^{\lor}}\circ\Psi_h$.
Similarly, we obtain
$\Psi_h\circ\del_{E,h}=\delbar_{E^{\lor}}\circ\Psi_{h}$.
Because $\kappa=\Psi_C^{-1}\circ\Psi_h$,
we obtain the claim of Lemma \ref{lem;22.8.29.30}.
\hfill\qed

\begin{lem}
$(E_{\real,\kappa},\nabla_{\real,\kappa},C_{\real,\kappa})$
 is a $\GL(n,\real)$-harmonic bundle.
\end{lem}
\pf
It follows from Lemma \ref{lem;22.8.29.31}.
\hfill\qed

\vspace{.1in}
Let $(V_{\real},\nabla_{\real},C_{\real})$
be a $\GL(n,\real)$-harmonic bundle.
We set $V_{\cnum}=\cnum\otimes_{\real}V_{\real}$.
Let $h$ be the Hermitian metric of $V_{\cnum}$
induced by $C_{\real}$,
and let $\nabla_{\cnum}$ denote the induced flat connection.
Then, $(V_{\cnum},\nabla_{\cnum},h)$ is a harmonic bundle.
Let $(V_{\cnum},\delbar_{V_{\cnum}},\theta)$
denote the Higgs bundle
underlying $(V_{\cnum},\nabla,h)$.
Let $C$ denote the holomorphic perfect symmetric pairing
of $V_{\cnum}$ induced by $C_{\real}$.
Then, $C$ is a real structure of
the harmonic bundle  $(V_{\cnum},\delbar_{V_{\cnum}},\theta,h)$.
It is easy to observe the following.
\begin{prop}
By the constructions,
 $\GL(n,\real)$-harmonic bundles
are equivalent to harmonic bundles
equipped with a real structure.
\hfill\qed
\end{prop}

\begin{rem}
Let $U\subset\cnum$ be an open neighbourhood of $0$.
Let $(E,\delbar_E,\theta,h)$ be a harmonic bundle
on $U^{\ast}=U\setminus\{0\}$
equipped with a real structure $C$.
Let $(V_{\real},\nabla_{\real},C_{\real})$
be the corresponding $\GL(n,\real)$-harmonic bundle on $U^{\ast}$.
We obtain the $\real$-local system $L_{\real}$ on $U^{\ast}$
obtained as the sheaf of flat sections of $(V_{\real},\nabla_{\real})$. 
If $(E,\delbar_E,\theta,h)$ is wild at $0$,
the $\cnum$-local system 
$L_{\cnum}=L_{\real}\otimes\cnum$
is equipped with the three level of filtrations
along each ray $\{te^{\sqrt{-1}\theta}\,|\,0<t<\epsilon\}$,
the Stokes filtrations,
the parabolic filtrations,
and the weight filtrations.
They are important for our understanding
of the associated meromorphic flat bundle.
Because the filtrations are described by the growth orders of flat sections
with respect to the harmonic metric, 
we can easily observe that
they are induced by the filtrations of $L_{\real}$ along the ray.
\hfill\qed
\end{rem}

\subsubsection{$\SL(n,\real)$-harmonic bundles}

Let $\underline{\real}_{Y}$
denote the product line bundle
$Y\times\real$.
It has a naturally defined flat connection $\nabla_{\real}$
and a positive definite symmetric pairing $C_{\real}$.
The tuple
$(\underline{\real}_Y,\nabla_{\real},C_{\real})$
is a $\GL(1,\real)$-harmonic bundle.
A $\GL(n,\real)$-harmonic bundle
$(V_{\real},\nabla_{\real},C_{\real})$
is called an $\SL(n,\real)$-harmonic bundle
when an isomorphism
$(\det(V_{\real}),\nabla_{\real},C_{\real})
\simeq
 (\underline{\real}_Y,\nabla_{\real},C_{\real})$
is equipped.

Let $h_{0,Y}$ denote the Hermitian metric
of $\nbigo_Y$ defined by $h_{0,Y}(1,1)=1$.
Let $C_{0,Y}$ be the holomorphic symmetric pairing of
$\nbigo_Y$ defined by $C_{0,Y}(1,1)=1$.
Then,
$(\nbigo_Y,0,h_{0,Y})$ with $C_{0,Y}$
is a harmonic bundle with a real structure on $Y$.

\begin{prop}
$\SL(n,\real)$-harmonic bundles
are equivalent to
harmonic bundles 
$(E,\delbar_E,\theta,h)$
with real structure $C$ 
equipped with 
an isomorphism
$(\det(E),\delbar_{\det(E)},\tr\theta,\det(h),\det(C))
 \simeq 
 (\nbigo_Y,0,h_{0,Y},C_{0,Y})$.
\hfill\qed
\end{prop}

\subsubsection{Compatible harmonic metrics}

Let $(E,\delbar_E,\theta)$ be a Higgs bundle on $Y$.
A holomorphic symmetric pairing $C$ of $E$
is called a symmetric pairing of $(E,\delbar_E,\theta)$
if $\theta$ is self-adjoint with respect to $C$.
(See Definition {\rm\ref{df;22.9.3.10}}.) 
When $(E,\delbar_E,\theta)$ is equipped with
a symmetric pairing $C$,
we say that a harmonic metric $h$ of $(E,\delbar_E,\theta)$
is compatible with $C$
if $C$ is a real structure of
the harmonic bundle $(E,\delbar_E,\theta,h)$,
i.e., $h_{|P}$ is compatible with $C_{|P}$
for any $P\in Y$.
Let $\Harm(E,\delbar_E,\theta;C)$
denote the set of harmonic metrics $h$ of
$(E,\delbar_E,\theta)$ compatible with $C$.

From $h_i\in\Harm(E,\delbar_E,\theta;C)$ $(i=1,2)$,
we obtain the real structures $\kappa_i$ of $E$,
and $\GL(n,\real)$-harmonic bundles
$(E_{\real,\kappa_i},\nabla_{\real,\kappa_i},C_{\real,\kappa_i})$.
\begin{prop}
\label{prop;22.9.7.11}
If $[\theta,s(h_1,h_2)]=0$ and $\delbar_E(s(h_1,h_2))=0$,
then $s(h_1,h_2)^{1/2}$ induces 
an isomorphism of $\GL(n,\real)$-harmonic bundles
$(E_{\real,\kappa_2},\nabla_{\real,\kappa_2},C_{\real,\kappa_2})
 \simeq
 (E_{\real,\kappa_1},\nabla_{\real,\kappa_1},C_{\real,\kappa_1})$.
\end{prop}
\pf
According to Lemma \ref{lem;22.9.7.2},
$s(h_1,h_2)^{1/2}$ induces an isomorphism
$(E_{\real,\kappa_2},C_{\real,\kappa_2})
 \simeq
 (E_{\real,\kappa_1},C_{\real,\kappa_1})$.
Because $s(h_1,h_2)$ is self-adjoint with respect to $h_1$,
we obtain $[\theta^{\dagger}_{h_1},s(h_1,h_2)]=0$,
$\del_{E,h_1}(s(h_1,h_2))=0$,
and $\nabla_{h_1}(s(h_1,h_2))=0$.
The eigenvalues of $s(h_1,h_2)$ are constant.
We obtain the eigen decomposition
$E=\bigoplus_{a>0}E_a$
of $s(h_1,h_2)$,
which is orthogonal with respect to both $h_i$,
and $h_2=ah_1$ on $E_a$.
The decomposition is compatible with
$\delbar_E$ and $\theta$.
We obtain that
$\theta^{\dagger}_{h_2}=\theta^{\dagger}_{h_1}$
and
$\del_{E,h_2}=\del_{E,h_1}$.
We obtain $\DD^1_{h_1}=\DD^1_{h_2}$,
and $\DD^1_{h_1}(s(h_1,h_2)^{1/2})=0$.
Hence, $s(h_1,h_2)^{1/2}$
induces an isomorphism of $\GL(n,\real)$-harmonic bundles.
\hfill\qed

\subsubsection{Canonical harmonic metric
in the regular semisimple case}

\begin{df}
\label{df;22.9.3.9}
The Higgs bundle $(E,\delbar_E,\theta)$
is called regular semisimple if
the following holds for any $P\in Y$.
\begin{itemize}
 \item Let $(Y_P,z)$ be a holomorphic coordinate neighbourhood
       around $P$.
       Let $f_P$ be the endomorphism of $E_{|Y_P}$
       determined by $\theta=f_P\,dz$.
       Then, $f_{P|Q}$ $(Q\in Y_P)$
       are regular semisimple,
       i.e.,
       the multiplicity of each eigenvalue of $f_{P|Q}$
       is  $1$.
\hfill\qed
\end{itemize}
\end{df}

\begin{prop}
\label{prop;22.9.3.21}
Suppose that $(E,\delbar_E,\theta)$ is regular semisimple
and that it is equipped with a non-degenerate
symmetric pairing.
Then, there exists a unique harmonic metric $h^{\can}$
of $(E,\delbar_E,\theta)$
satisfying the following conditions.
\begin{itemize}
 \item $h^{\can}$ is compatible with $C$,
       i.e.,
       $h^{\can}\in\Harm(E,\delbar_E,\theta;C)$.
 \item Let $P$ be any point of $Y$.
       Let $(Y_P,z)$ and $f_P$ be as in
       Definition {\rm\ref{df;22.9.3.9}}.
       Then, the eigen decomposition of $f_{P|Q}$
       is orthogonal with respect to $h_{0|Q}$.
\end{itemize}
The metric $h^{\can}$ is called the canonical metric of
$(E,\delbar_E,\theta,C)$.
\end{prop}
\pf
By Lemma \ref{lem;22.9.3.20},
there exists a Hermitian metric $h^{\can}$ of $E$
satisfying the conditions.
We can easily check that it is a harmonic metric
of the Higgs bundle.
\hfill\qed

\begin{rem}
In general, there are many other harmonic metric
of $(E,\delbar_E,\theta)$ compatible with $C$.
\hfill\qed 
\end{rem}

\subsection{An existence theorem of
compatible harmonic metrics}

\subsubsection{Statement}

Let $X$ be any Riemann surface.
Let $(E,\delbar_E,\theta)$ be a Higgs bundle on $X$.

\begin{df}
The Higgs bundle $(E,\delbar_E,\theta)$
is called generically regular semisimple
if there exists a discrete subset $Z\subset X$
such that $(E,\delbar_E,\theta)_{|X\setminus Z}$
is regular semisimple.
\hfill\qed
\end{df}

We shall prove the following theorem
in \S\ref{subsection;22.9.3.11}--\ref{subsection;22.9.3.12}.
\begin{thm}
\label{thm;22.8.31.10}
We assume that
$(E,\delbar_E,\theta,h)$
is generically regular semisimple,
and that it is equipped with
a non-degenerate symmetric pairing $C$.
Then,
there exists a harmonic metric $h$ of
the Higgs bundle $(E,\delbar_E,\theta)$ 
 which is compatible with $C$,
i.e., $\Harm(E,\delbar_E,\theta;C)\neq\emptyset$. 
\end{thm}

\begin{rem}
If $X$ is compact,
the existence of a non-degenerate pairing
implies $\deg(E)=0$,
and we can obtain the polystability
of $(E,\delbar_E,\theta)$
by using the argument in the proof of Theorem {\rm\ref{thm;22.8.26.12}}.
Hence, we can prove Theorem {\rm\ref{thm;22.8.31.10}}
by using the classical theorem of Hitchin {\rm\cite{Hitchin-self-duality}}
and Simpson {\rm \cite{s1}}
(see Theorem {\rm\ref{thm;22.9.6.30}}).
\hfill\qed
\end{rem}

\subsubsection{Compact case}
\label{subsection;22.9.22.1}

First, let us study the case where $X$ is compact.
Indeed, in this case,
we can also obtain the uniqueness of a harmonic metric
compatible with $C$.
\begin{prop}
\label{prop;23.7.21.1}
Let $(E,\delbar_E,\theta)$ be a generically regular semisimple
Higgs bundle equipped with a non-degenerate symmetric pairing $C$
on a compact Riemann surface $X$.
Then, $(E,\delbar_E,\theta)$ 
is polystable of degree $0$ and
has a unique harmonic metric
compatible with $C$.
\end{prop}
\pf
Because $E$ is equipped with a non-degenerate symmetric pairing $C$,
we obtain $\det(E)\otimes\det(E)\simeq\nbigo_X$,
which implies $\deg(E)=0$.
Let $E'\subset E$ be a subbundle
such that $\theta(E')\subset E'\otimes\Omega^1_X$.
Let $P\in X\setminus Z$.
Let $(X_P,z)$ be a holomorphic coordinate neighbourhood
around $P$ in $X\setminus Z$.
We obtain the endomorphism $f$ of $E_{|X_P}$
such that $\theta=f\,dz$.
There exist holomorphic functions $\alpha_1,\ldots,\alpha_{\rank E}$
on $X_P$
and the eigen decomposition
$(E_{|X_P},f)=\bigoplus_{i=1}^{\rank E}(E_i,\alpha_i\id_{E_i})$.
It is orthogonal with respect to $C_{|X_P}$,
and hence the restriction of $C$ to each $E_i$ is non-degenerate.
Because $E'_{|X_P}$ is the direct sum of some of $E_i$,
the restriction of $C$ to $E'_{|X_P}$ is also non-degenerate.
We obtain that
the restriction $C'$ of $C$ to $E'$
is non-degenerate on $X\setminus Z$.
Hence,
we obtain a monomorphism
$\det(C'):\det(E')\otimes\det(E')\to\nbigo_X$.
It implies that $\deg(E')\leq 0$.
Moreover, if $\deg(E')=0$,
we obtain that $C'$ is non-degenerate on $X$.
We obtain the orthogonal decomposition
$E=E'\oplus E^{\prime\bot}$ with respect to $C$,
which is compatible with $\theta$.
Hence, we obtain that
$(E,\delbar_E,\theta)$ is polystable.

There exists a decomposition
\begin{equation}
\label{eq;22.9.9.1}
(E,\delbar_E,\theta)
=\bigoplus (E_j,\delbar_{E_j},\theta_j)
\end{equation}
into stable Higgs bundles of degree $0$,
which is orthogonal with respect to $C$.
Because of the generic regular semisimplicity,
we obtain
$(E_j,\delbar_{E_j},\theta_j)
\not\simeq
(E_k,\delbar_{E_k},\theta_k)$ $(j\neq k)$.
By the theorem of Hitchin and Simpson,
there exists a harmonic metric $h_1$ of $(E,\delbar_E,\theta)$.
It induces a harmonic metric $h_1^{\lor}$ of
$(E^{\lor},\delbar_{E^{\lor}},\theta^{\lor})$.
The decomposition (\ref{eq;22.9.9.1}) is orthogonal
with respect to
both $h_1$ and $\Psi_C^{\ast}(h_1^{\lor})$.
There exists $c_j>0$ such that
$\Psi_{C}^{\ast}(h_1^{\lor})=c_j^2h_1$ on $E_j$.
We set $h_2=\bigoplus c_j\cdot h_{1|E_j}$.
Then, $h_2$ is compatible with $C$.
The uniqueness is also clear.
\hfill\qed

\subsubsection{Local estimate in the regular semisimple case}
\label{subsection;22.9.3.11}

For $R>0$,
we set $U(R):=\{z\in\cnum\,|\,|z|<R\}$.
Let $(E,\delbar_E,\theta)$ be a Higgs bundle
of rank $r$ on $U(R)$.
We obtain the endomorphism $f$ of $E$
by $\theta=f\,dz$.
Assume the following.
\begin{itemize}
 \item There exist holomorphic functions $\alpha_1,\ldots,\alpha_r$
       and a decomposition
\begin{equation}
\label{eq;22.8.30.1}
 (E,f)=\bigoplus (E_i,\alpha_i\id_{E_i}).
\end{equation}
 \item There exist $0<A_1,A_2$ such that
       $|\alpha_i|<A_1$
       and $|\alpha_i-\alpha_j|>A_2$ $(i\neq j)$.
\end{itemize}

Let $C$ be a holomorphic non-degenerate symmetric bilinear form
of the Higgs bundle $(E,\delbar_E,\theta)$.
The decomposition (\ref{eq;22.8.30.1})
is orthogonal with respect to $C$.
We have the canonical harmonic metric $h^{\can}$
of $(E,\delbar_E,\theta,C)$
as in Proposition \ref{prop;22.9.3.21}.

\begin{prop}
\label{prop;22.9.1.2}
Let $0<R_1<R$.
There exist positive constants
$C_i$ $(i=1,2)$
depending only on $A_j$ $(j=1,2)$ and $R_1,R$
such that the following holds
on $U(R_1)$
for any $h\in\Harm(E,\delbar_E,\theta;C)$:
\[
 |f|_h\leq C_1,
 \quad
 |s(h^{\can},h)|_{h^{\can}}
+|s(h^{\can},h)^{-1}|_{h^{\can}}
 \leq C_2.
\]
\end{prop}
\pf
The estimate for $|f|_h$
is given in \cite[Lemma 2.7]{s5}.
We obtain the estimate for $s(h^{\can},h)$
from Proposition \ref{prop;22.8.31.2}.
\hfill\qed

\subsubsection{Local estimate in the generically regular semisimple case}

Let $X$, $(E,\delbar_E,\theta)$ and $C$ 
be as in Theorem \ref{thm;22.8.31.10}.
Assume that $X$ is non-compact.
We choose a K\"ahler metric $g_{X}$ of $X$.
There exists a discrete subset $Z\subset X$
such that $(E,\delbar_E,\theta)_{|X\setminus Z}$
is regular semisimple.
Let $K_1\subset X$ be any relatively compact open subset of $X$.
For simplicity, we assume that $Z\cap \del K_1=\emptyset$.
Let $K_2$ be a neighbourhood of the closure $\Kbar_1$ in $X$
such that $K_2\cap Z=K_1\cap Z$.
Let $h_{1}$ be any Hermitian metric of $E$
which is compatible with $C$.
(See Lemma \ref{lem;22.9.3.40}.)

\begin{prop}
\label{prop;22.9.22.3}
There exists $B(K_1,K_2)>0$ such that the following inequality holds
on $K_1$
for any $h\in\Harm((E,\delbar_E,\theta;C)_{|K_2})$:
\[
       \bigl|
       s(h_{1},h)
       \bigr|_{h_{1}}
       +
      \bigl|
       s(h_{1},h)^{-1}
       \bigr|_{h_{1}}
       \leq B(K_1,K_2).
\]       
\end{prop}
\pf
Let $P\in \del K_1$.
There exists a relatively compact neighbourhood $X_P$ of $P$
in $K_2$
such that
$(E,\delbar_E,\theta,h)_{|X_P}$ is regular semisimple.
Let $h^{\can}_{X\setminus Z}$ denote the canonical harmonic metric
of $(E,\delbar_E,\theta)_{|X\setminus Z}$
compatible with $C_{|X\setminus Z}$.
By Proposition \ref{prop;22.9.1.2},
there exists $B(X_P)>0$
such that the following holds on $X_P$
for any $h\in\Harm((E,\delbar_E,\theta;C)_{|K_2})$:
\[
       \bigl|
       s(h^{\can}_{X\setminus Z},h)
       \bigr|_{h^{\can}_{X\setminus Z}}
       +
      \bigl|
       s(h^{\can}_{X\setminus Z},h)^{-1}
       \bigr|_{h^{\can}_{X\setminus Z}}
       \leq B(X_P).
\]
Let $N$ be a relatively compact neighbourhood of
$\del K_1$ in $\bigcup_{P\in\del K_1}X_P$.
There exists $B(N)>0$ depending only on $N$
such that the following holds on $N$
for any $h\in\Harm((E,\delbar_E,\theta;C)_{|K_2})$:
\[
      \bigl|
       s(h^{\can}_{X\setminus Z},h)
       \bigr|_{h^{\can}_{X\setminus Z}}
       +
      \bigl|
       s(h^{\can}_{X\setminus Z},h)^{-1}
       \bigr|_{h^{\can}_{X\setminus Z}}
       \leq B(N).
\]
There exists $B'(N)>0$ depending only on $N$
such that the following holds on $N$
for any 
$h\in\Harm((E,\delbar_E,\theta;C)_{|K_2})$:
\[
      \bigl|
       s(h_{1},h)
       \bigr|_{h_{1}}
       +
      \bigl|
       s(h_{1},h)^{-1}
       \bigr|_{h_{1}}
       \leq B'(N).
\]

We recall the following inequality \cite[Lemma 3.1]{s1}
on $K_2$
for any $h\in\Harm((E,\delbar_E,\theta;C)_{|K_2})$:
\[
 \sqrt{-1}\Lambda\delbar\del
 \log\Tr\bigl(
 s(h_1,h)
 \bigr)
 \leq
 \bigl|
 \Lambda F(h_1)
 \bigr|_{h_1}.
\]
There exists a function $\beta$ on $K_2$
such that
$\sqrt{-1}\Lambda\delbar\del\beta=
\bigl| \Lambda F(h_1)
\bigr|_{h_1}$.
We obtain
\[
 \sqrt{-1}\Lambda\delbar\del
 \Bigl(
 \log\Tr\bigl(
 s(h_1,h)
 \bigr)
 -\beta
 \Bigr)
 \leq 0
\]
on $K_2$.
By the maximum principle,
we obtain 
\[
 \max_{K_1}
 \Bigl(
  \log\Tr\bigl(
 s(h_1,h)
 \bigr)
 -\beta
 \Bigr)
 \leq
 \max_{N}
 \Bigl(
  \log\Tr\bigl(
 s(h_1,h)
 \bigr)
 -\beta
 \Bigr).
\]
We also note that $\det(s(h_1,h))=1$.
Hence, we obtain the claim of the lemma.
\hfill\qed

\subsubsection{Proof of Theorem \ref{thm;22.8.31.10}}
\label{subsection;22.9.3.12}

We have the isomorphism of the Higgs bundles
$\Psi_C:(E,\delbar_E,\theta)
\simeq
 (E^{\lor},\delbar_{E^{\lor}},\theta^{\lor})$.

Let $\{X_i\}$ be a smooth exhausting family of $X$
as in \cite[Definition 2.5]{Li-Mochizuki2}.
Let $h_1$ be a Hermitian metric of $E$
which is compatible with $C$.
We have the induced metric $h_1$ on $E^{\lor}$.
The morphism $\Psi_C:E\lrarr E^{\lor}$
is an isometry with respect to $h_1$ and $h_1^{\lor}$.

By \cite[Theorem 2]{Donaldson-boundary-value}
and \cite[Proposition 2.1]{Li-Mochizuki2},
there exists a unique harmonic metric $h_{X_i}$
of $(E,\delbar_E,\theta)_{|X_i}$
such that $h_{X_i|\del X_i}=h_{1|\del X_i}$.
We have the induced harmonic metric
$h_{X_i}^{\lor}$ of 
$(E^{\lor},\delbar_{E^{\lor}},\theta^{\lor})_{|X_i}$.
Because $\Psi_{C|\del X_i}$ is isometric
with respect to
$h_{X_i|\del X_i}$
and
$h_{X_i|\del X_i}^{\lor}$,
we obtain that
$\Psi_{C|X_i}$  is isometric
with respect to
$h_{X_i}$ and $h_{X_i}^{\lor}$.
By \cite[Proposition 2.6]{Li-Mochizuki2}
and Proposition \ref{prop;22.9.22.3},
$\{h_{X_i}\}$ has a convergent subsequence
whose limit is denoted by $h$.
Then, $h$ is a harmonic metric
compatible with $C$.
\hfill\qed

\subsection{Some sufficient conditions
  for generically regular semisimplicity}
\label{subsection;22.10.10.1}

\subsubsection{Spectral curves}

Let $(E,\delbar_E,\theta)$ be a Higgs bundle
of rank $r$ on a Riemann surface $X$.
We may regard $E$ as a module over
the sheaf of algebras $\Sym\Theta_X$,
where $\Theta_X$ denote the tangent sheaf of $X$.
There exists the coherent $\nbigo_{T^{\ast}X}$-module
$\nbigf(E,\theta)$
with an isomorphism of $\Sym\Theta_X$-modules
$\pi_{\ast}\nbigf(E,\theta)\simeq E$
where $\pi:T^{\ast}X\to X$ denote the projection.
The support $\Sigma_{E,\theta}$ of
$\nbigf(E,\theta)$ is called the spectral curve of $(E,\theta)$.
(See \cite{Beauville-Narasimhan-Ramanan, Hitchin-stable-bundles}.)
Let us recall some related notions in a way convenient to us.

For any $P\in X$,
let $\iota_P:T^{\ast}_PX\to T^{\ast}X$ denote the inclusion.
We obtain the $\nbigo_{T^{\ast}_PX}$-module
$\iota_P^{\ast}\nbigf(E,\theta)$.
We have the decomposition by the supports
\[
 \iota_P^{\ast}\nbigf(E,\theta)
 =\bigoplus_{Q\in T^{\ast}_PX\cap \Sigma_{E,\theta}}
 \iota_P^{\ast}\nbigf(E,\theta)_Q,
\]
where the support of $\iota_P^{\ast}\nbigf(E,\theta)_Q$
is $\{Q\}$.
Each $\iota_P^{\ast}\nbigf(E,\theta)_Q$
is naturally a finite dimensional $\cnum$-vector space.
We set $\gminim(Q):=\dim_{\cnum}\iota_P^{\ast}\nbigf(E,\theta)_Q$.
We obtain a map
$\gminim:\Sigma_{E,\theta}\to \seisuu_{>0}$.
We have $\sum_{\pi(Q)=P}\gminim(Q)=r$ for any $P\in X$.
In particular,
we obtain
$|\Sigma_{E,\theta}\cap T_P^{\ast}X|\leq r$.
We set
$n(E,\theta)=\max_{P\in X}|\Sigma(E,\theta)\cap T_P^{\ast}X|\leq r$.
The generically regular semisimplicity condition
for $(E,\delbar_E,\theta)$ is equivalent to $n(E,\theta)=r$.

We set
$D(E,\theta):=\bigl\{P\in X\,\big|\,
 |\Sigma_{E,\theta}\cap T_P^{\ast}X|<n(E,\theta)\bigr\}$.
It is easy to see that $D(E,\theta)$ 
is discrete in $X$.
We set
$\Sigma_{E,\theta}^{\circ}:=
\Sigma_{E,\theta}\setminus \pi^{-1}(D(E,\theta))$.
It is a closed complex submanifold of
$T^{\ast}(X\setminus D(E,\theta))$,
and the projection
$\Sigma_{E,\theta}^{\circ}\to X\setminus D(E,\theta)$
is a local homeomorphism.
The function $\gminim$ is locally constant
on $\Sigma_{E,\theta}^{\circ}\to\seisuu_{>0}$.
We set
$(\Sigma_{E,\theta}^{\circ})_k:=
\Sigma_{E,\theta}^{\circ}
\cap
\gminim^{-1}(k)$.
Let $(\Sigma_{E,\theta})_k$
denote the closure of $(\Sigma_{E,\theta}^{\circ})_k$
in $\Sigma_{E,\theta}$.
We have
$\Sigma_{E,\theta}=\bigcup_{k}
(\Sigma_{E,\theta})_k$.

A union of
some irreducible components of $\Sigma_{E,\theta}$
is called a spectral subcurve of $\Sigma_{E,\theta}$.
For example, $(\Sigma_{E,\theta})_k$
are spectral subcurves.
A spectral subcurve of $\Sigma_{E,\theta}$
is a purely $1$-dimensional complex analytic closed subset
of $\Sigma_{E,\theta}$.
For a spectral subcurve $\Sigma_1$,
we set 
$\Sigma_1^{\circ}:=
\Sigma_1\setminus \pi^{-1}(D(E,\theta))$.
It is easy to see that
(i) $\Sigma_1^{\circ}$ is a smooth complex submanifold of
$T^{\ast}(X\setminus D(E,\theta))$,
(ii) $\Sigma_1^{\circ}
\to X\setminus D(E,\theta)$ is a local homeomorphism,
(iii) $\Sigma_1$ is the closure of $\Sigma_1^{\circ}$ in $T^{\ast}X$.
For any $P\in X\setminus D(E,\theta)$,
we set $n(\Sigma_1):=\bigl|\Sigma_1\cap T_{P}^{\ast}X\bigr|$
and $r(\Sigma_1):=
\sum_{Q\in T_P^{\ast}X\cap \Sigma_1} \gminim(Q)$,
which are independent of the choice of $P$.
In particular,
we set
$n(k):=n((\Sigma_{E,\theta})_k)$
and 
$r(k):=r((\Sigma_{E,\theta})_k)$. 
Note that $r(k)=k\cdot n(k)$.

\subsubsection{Characteristic polynomials and
some Higgs subsheaves in the local case}

Let $(U,z)$ be a connected holomorphic chart of $X$,
i.e., $U$ is a connected open subset of $X$,
and $z$ is a holomorphic coordinate on $U$.
Let $f_U$ be the endomorphism of $E_{|U}$
defined by $\theta=f_U\,dz$.
We obtain the characteristic polynomial
$\det(T\id_{E_{|U}}-f_U)=\sum_{j=0}^r a_j(z)T^j$
which is a monic.
Then, the polynomial
\[
P_{E,\theta,U}(T)=\sum_{j=0}^r a_j(z)(dz)^{r-j}T^j
\in \bigoplus_{j=0}^r H^0(U,K_U^{r-j})T^{j}
\]
is independent of the choice of a coordinate $z$ on $U$.
It is called the characteristic polynomial of
$(E,\delbar_E,\theta)_{|U}$.
Under the isomorphism
$U\times \cnum\simeq T^{\ast}U$ induced by
$(z,T)\longmapsto (z,T\,dz)$,
the spectral curve 
$\Sigma_{E,\theta}\cap T^{\ast}U$
of $(E,\delbar_E,\theta)_{|U}$
is equal to 
$\{(z,T)\in U\times\cnum\,|\,\sum a_j(z)T^j=0\}$.
For $(z_0,\alpha\,dz)\in \Sigma_{E,\theta}\cap T^{\ast}U$,
$\gminim(z_0,\alpha\,dz)$ is equal to
the multiplicity of
the root $\alpha$ of the polynomial $\sum a_j(z_0)T^j$.

\vspace{.1in}

Let $\Sigma_1\subset \Sigma_{E,\theta}\cap T^{\ast}U$
be a spectral subcurve.
For $z_0\in U\setminus D(E,\theta)$, we obtain the complex numbers
$c_j^{\Sigma_1}(z_0)$ $(0\leq j\leq n(\Sigma_1))$
and
$\ctilde_j^{\Sigma_1}(z_0)$ $(0\leq j\leq r(\Sigma_1))$
by setting
\[
 \sum_{j=0}^{n(\Sigma_1)} c_j^{\Sigma_1}(z_0)T^j
 =\prod_{(z_0,\alpha\,dz)\in \Sigma_1}
 (T-\alpha),
\quad\quad
 \sum_{j=0}^{r(\Sigma_1)} \ctilde_j^{\Sigma_1}(z_0)T^j
 =\prod_{(z_0,\alpha\,dz)\in \Sigma_1}
 (T-\alpha)^{\gminim(z_0,\alpha\,dz)}.
\]
We obtain holomorphic functions
$c_j^{\Sigma_1}$ and $\ctilde_j^{\Sigma_1}$
on $U\setminus D(E,\theta)$.
Because the eigenvalues of $f_U$ are locally bounded,
they extend to holomorphic functions on $U$,
which are also denoted by 
$c_j^{\Sigma_1}$ and $\ctilde_j^{\Sigma_1}$.
We obtain the following endomorphisms
of the $\nbigo_U$-module $E_{|U}$:
\[
 F^{\Sigma_1}=\sum_{j=0}^{n(\Sigma_1)}c^{\Sigma_1}_jf_U^j,
 \quad
 \quad
 \Ftilde^{\Sigma_1}=
 \sum_{j=0}^{r(\Sigma_1)}\ctilde^{\Sigma_1}_jf_U^j.
\]
We obtain the $\nbigo_U$-submodules
$\nbigf(\Sigma_1,0)=\Ker F^{\Sigma_1}$
and
$\nbigf(\Sigma_1,1)=\Ker\Ftilde^{\Sigma_1}$.
Note that $\nbigf(\Sigma_1,i)\neq 0$ $(i=0,1)$.

\begin{lem}
\label{lem;22.10.10.2}
$\nbigf(\Sigma_1,i)$ $(i=0,1)$
are holomorphic subbundles of $E_{|U}$,
i.e., 
$E_{|U}\big/\nbigf(\Sigma_1,i)$ $(i=0,1)$
are torsion-free $\nbigo_U$-modules. 
We also have
$f_U(\nbigf(\Sigma_1,i))\subset\nbigf(\Sigma_1,i)$.
\end{lem}
\pf
We obtain the first claim from 
$E_{|U}\big/\nbigf(\Sigma_1,0)\simeq
\Image F^{\Sigma_1}\subset E_{|U}$
and
$E_{|U}\big/\nbigf(\Sigma_1,1)\simeq
\Image \Ftilde^{\Sigma_1}\subset E_{|U}$.
We obtain the second claim from the commutativity
$[f,F^{\Sigma_1}]=[f,\Ftilde^{\Sigma_1}]=0$.
\hfill\qed

\vspace{.1in}
In particular,
$\nbigf(
\Sigma_1,i)$ are Higgs subbundles
of $(E,\theta)_{|U}$.
Let $\theta(\Sigma_1,i)$ denote the Higgs field
of $\nbigf(\Sigma_1,i)$ induced by $\theta_{|U}$.

\begin{lem}
\label{lem;22.10.10.3}
 We have
$\det(T\id_{\nbigf(\Sigma_1,1)}-f_{U|\nbigf(\Sigma_1,1)})
=\sum_{j=0}^{r(\Sigma_1)}\ctilde^{\Sigma_1}_jT^j$.
As a result, we obtain
$P_{\nbigf(\Sigma_1,1),\theta(\Sigma_1,1)}(T)=
 \sum_{j=0}^{r(\Sigma_1)}\ctilde^{\Sigma_1}_j
 (dz)^{r(\Sigma_1)-j}T^j$.
\end{lem}
\pf
It is easy to see
$\det(T\id_{\nbigf(\Sigma_1,1)}-f_{U|\nbigf(\Sigma_1,1)})
_{|U\setminus D(E,\theta)}
=\Bigl(
\sum_{j=0}^{r(\Sigma_1)}\ctilde^{\Sigma_1}_jT^j
\Bigr)_{|U\setminus D(E,\theta)}$.
Then, we obtain the claim of the lemma.
\hfill\qed

\begin{rem}
The rank of $\nbigf(\Sigma_1,0)$ is not necessarily
equal to $n(\Sigma_1)$.
The characteristic polynomial of
$f_{U|\nbigf(\Sigma_1,0)}$ is not necessarily equal to
$\sum_{j=0}^{n(\Sigma_1)}c_j^{\Sigma_1}T^j$.
\hfill\qed
\end{rem}

We set $\Sigma_{E,\theta,k,U}:=(\Sigma_{E,\theta})_k\cap T^{\ast}U$
and 
\[
 P_{E,\theta,k,U}(T)=
 \sum_{j=0}^{n(k)}
 c^{\Sigma_{E,\theta,k,U}}_{j}(z)(dz)^{n(k)-j}T^j.
\]
Note that by definition, 
\[P_{\nbigf(\Sigma_{E,\theta,k,U},1),\theta(\Sigma_{E,\theta,k,U},1),U}(T)= \sum_{j=0}^{r(k)}\ctilde^{\Sigma_{E,\theta,k,U}}_{j}(z)(dz)^{r(k)-j}T^j=P_{E,\theta,k,U}(T)^k.
\]
We then have the factorization
$P_{E,\theta,U}(T)=
\prod_{k\geq 1}
 P_{E,\theta,k,U}(T)^k$.
 
\subsubsection{Some sufficient conditions for
generically regular semisimplicity}

There uniquely exists
\[
 P_{E,\theta}(T)
 \in \bigoplus_{j=0}^r H^0(X,K_X^{r-j})T^{j}
 \quad
 \left(
 \mbox{\textrm resp. }
 P_{E,\theta,k}(T)
 \in 
 \bigoplus_{j=0}^{n(k)}H^0(X,K_X^{n(k)-j})T^j
 \right)
\]
such that the restriction of $P_{E,\theta}(T)$
(resp. $P_{E,\theta,k}(T)$) to $U$
is $P_{E,\theta,U}(T)$
(resp. $P_{E,\theta,k,U}(T)$)
for any holomorphic chart $(U,z)$.
We have the factorization
$P_{E,\theta}(T)=\prod_{k\geq 1}P_{E,\theta,k}(T)^k$.
The polynomial $P_{E,\theta}(T)$ is called
the characteristic polynomial of
the Higgs bundle $(E,\delbar_E,\theta)$.

\begin{prop}
\label{prop;22.10.11.1}
The Higgs bundle $(E,\delbar_E,\theta)$ is generically regular semisimple
if the characteristic polynomial is irreducible.
\end{prop}
\pf
The irreducibility of
$P_{E,\theta}(T)$ implies that
$P_{E,\theta}(T)=P_{E,\theta,1}(T)$.
Hence, $(E,\delbar_E,\theta)$ is 
generically regular semisimple.
\hfill\qed 

\begin{prop}
\label{prop;22.10.11.2}
Suppose that
$(E,\delbar_E,\theta)$ has no proper Higgs subbundle.
Then, the characteristic polynomial
$P_{E,\theta}$ is irreducible.
In particular, $(E,\delbar_E,\theta)$ is 
generically regular semisimple.
\end{prop}
\pf
Let 
$\Sigma_1\subset \Sigma_{E,\theta}$
be a non-empty spectral subcurve.
We obtain non-zero Higgs subbundles
$(\nbigf(\Sigma_1,i),\theta(\Sigma_1,i))$
of $(E,\theta)$
such that
$\nbigf(\Sigma_1,i)_{|U}
=\nbigf(\Sigma_1\cap \pi^{-1}(U),i)$
for any holomorphic chart $(U,z)$,
where $\nbigf(\Sigma_1\cap \pi^{-1}(U),i)$
are as in Lemma \ref{lem;22.10.10.2}.
Because $E$ has no proper Higgs subbundle,
we obtain $\nbigf(\Sigma_1,0)=\nbigf(\Sigma_1,1)=E$.
In particular, we obtain $\Sigma_1=\Sigma_{E,\theta}$,
i.e., $\Sigma_{E,\theta}$ is irreducible.
There exists $k\geq 1$ such that
$\Sigma_{E,\theta}=(\Sigma_{E,\theta})_k$.
Because $\Sigma_{E,\theta}$ is irreducible,
the polynomial $P_{E,\theta,k}(T)$ is irreducible.
It remains to prove $k=1$.

Let $\varphi:\Xtilde\to \Sigma_{E,\theta}$ be a normalization.
Because $\Sigma_{E,\theta}$ is irreducible,
$\Xtilde$ is connected.
We set $\varphitilde:=\pi\circ\varphi:\Xtilde\to X$.
We obtain the Higgs bundle
$(\Etilde,\delbar_{\Etilde},\thetatilde):=
\varphitilde^{\ast}(E,\delbar_E,\theta)$
on $\Xtilde$.
The spectral curve
$\Sigma_{\Etilde,\thetatilde}$
is equal to the image of
$\Sigma_{E,\theta}\times_X\Xtilde
\subset \varphitilde^{\ast}(T^{\ast}X)$
by the naturally induced morphism
$\varphitilde^{\ast}(T^{\ast}X)\to T^{\ast}\Xtilde$.
We obtain the section $s_{\varphi}:
\Xtilde\to \Sigma_{\Etilde,\thetatilde}$
induced by $\varphi$.
Let $\Sigmatilde_1$ denote the image of
$s_{\varphi}$.
We obtain the Higgs subbundle
$\nbigftilde:=\nbigf(\Sigmatilde_1,0)\subset\Etilde$
with the induced Higgs field
$\theta_{\nbigftilde}$.

\begin{lem}
There exists a holomorphic subbundle
$\nbigl\subset\nbigftilde$ such that $\rank\nbigl=1$. 
\end{lem}
\pf
First, we consider the case when $\Xtilde$ is non-compact.
Because any holomorphic vector bundle on $\Xtilde$
has a holomorphic trivialization
(see \cite[Theorem 30.1]{Forster-book}),
we obtain the claim of the lemma in this case.
If $\Xtilde$ is compact,
let $\nbigo_{\Xtilde}(1)$ be an ample line bundle.
There exists $N>0$ such that
$\nbigftilde(N):=
\nbigftilde\otimes\nbigo_{\Xtilde}(N)$
has a non-trivial section $u$.
We obtain the $\nbigo_{\Xtilde}$-submodule
$\nbigl'\subset\nbigftilde(N)$
generated by $u$.
Let $\nbigl''\subset\nbigftilde(N)$
denote the $\nbigo_{\Xtilde}$-submodule
obtained as the pull back of the torsion-part of
$\nbigftilde(N)\bigl/\nbigl'$ by the projection.
Then,
$\nbigl=\nbigl''\otimes\nbigo_{\Xtilde}(-N)$
satisfies the desired condition.
\hfill\qed

\vspace{.1in}

We may naturally regard $s_{\varphi}$
as a holomorphic $1$-form on $\Xtilde$.
Because $\thetatilde-s_{\varphi}\id_{\Etilde}$
vanishes on $\nbigftilde$,
$\nbigl$ is a Higgs subbundle of $\Etilde$.

We set
$D:=D(E,\theta)$
and $\Dtilde:=\varphitilde^{-1}(D)$.
Let $\nbigo_X(\ast D)$ denote the sheaf of
meromorphic functions on $X$ whose poles are contained in $D$.
Similarly, 
let $\nbigo_{\Xtilde}(\ast \Dtilde)$ denote the sheaf of
meromorphic functions on $\Xtilde$
whose poles are contained in $\Dtilde$.
The generalized eigen decomposition induces
the projection
$\varphitilde^{\ast}(E\otimes\nbigo_X(\ast D),\theta)
\to
 (\nbigftilde\otimes\nbigo_{\Xtilde}(\ast\Dtilde),\theta_{\nbigftilde})$.
By the adjunction,
we obtain
$\psi:
E\otimes\nbigo_X(\ast D)
\lrarr
 \varphitilde_{\ast}(\nbigftilde\otimes\nbigo_{\Xtilde}(\ast \Dtilde))
 =\varphitilde_{\ast}(\nbigftilde)\otimes\nbigo_X(\ast D)$
compatible with the Higgs fields.
Because the restriction
$\psi_{|X\setminus D}$ is an isomorphism,
we obtain that $\psi$ is an isomorphism
of locally free $\nbigo_X(\ast D)$-modules.
We obtain the coherent $\nbigo_X(\ast D)$-submodule
$\varphitilde_{\ast}(\nbigl)\otimes\nbigo_X(\ast D)\subset
E\otimes\nbigo_X(\ast D)$.
We obtain the holomorphic subbundle
$\nbigg:=
\bigl(
\varphitilde_{\ast}(\nbigl)\otimes\nbigo_X(\ast D)
\bigr)
\cap E\subset E$
which is preserved by $\theta$.
Because $E$ has no proper Higgs subbundle,
we obtain $\nbigg=E$.
By considering the restriction of the characteristic polynomial
to $X\setminus D$,
we obtain $n(k)=r$ and thus $k=1$,
i.e.,
$(E,\delbar_E,\theta)$ is generically regular semisimple.
\hfill\qed

\section{Good filtered Higgs bundles with symmetric pairing}

\subsection{Pairings of filtered bundles}

\subsubsection{Pairings of locally free $\nbigo_X(\ast D)$-modules}

Let $X$ be a Riemann surface.
Let $D\subset X$ be a discrete subset.
Let $\nbigo_X(\ast D)$ denote the sheaf of
meromorphic functions on $X$ whose poles are contained in $D$.

Let $\nbigv$ be a locally free $\nbigo_X(\ast D)$-module
of finite rank.
Let $\nbigv^{\lor}$ denote the dual of $\nbigv$,
i.e.,
\[
 \nbigv^{\lor}=\nhom_{\nbigo_{X}(\ast D)}(\nbigv,\nbigo_X(\ast D)).
\]
The determinant bundle of $\nbigv$ is denoted by
$\det(\nbigv)$,
i.e.,
$\det(\nbigv)=\bigwedge^{\rank\nbigv}\nbigv$.
There exists a natural isomorphism
$\det(\nbigv^{\lor})\simeq \det(\nbigv)^{\lor}$.
For a morphism $f:\nbigv_1\to\nbigv_2$
of locally free $\nbigo_X(\ast D)$-modules,
$f^{\lor}:\nbigv_2^{\lor}\to\nbigv_1^{\lor}$
denotes the dual of $f$.
If $\rank\nbigv_1=\rank\nbigv_2$,
$\det(f):\det(\nbigv_1)\to\det(\nbigv_2)$
denotes the induced morphism.

For locally free $\nbigo_X(\ast D)$-modules $\nbigv_i$ $(i=1,2)$,
a pairing $C$ of $\nbigv_1$ and $\nbigv_2$
is a morphism
$C:\nbigv_1\otimes_{\nbigo_X(\ast D)}\nbigv_2\to\nbigo_X(\ast D)$.
It induces a morphism
$\Psi_C:\nbigv_1\to\nbigv_2^{\lor}$
by $\Psi_C(u)(v)=C(u\otimes v)$.
Let
$\exchange:\nbigv_1\otimes\nbigv_2\simeq\nbigv_2\otimes\nbigv_1$
denote the natural morphism defined by
$\exchange(v_1\otimes v_2)=v_2\otimes v_1$.
We obtain the pairing
$C\circ\exchange:\nbigv_2\otimes\nbigv_1\to\nbigo_X(\ast D)$.
We have
$\Psi_{C\circ\exchange}=\Psi_C^{\lor}$.
If $\rank\nbigv_1=\rank\nbigv_2$,
we obtain the induced pairing
$\det(C):\det(\nbigv_1)\otimes\det(\nbigv_2)\to\nbigo_X(\ast D)$.
We have $\det(\Psi_C)=\Psi_{\det(C)}$.

A pairing $C$ is called non-degenerate if $\Psi_C$ is an isomorphism.
It is equivalent to the condition that
$C\circ\exchange$ is non-degenerate.
It is also equivalent to the condition that
$\det(C)$ is non-degenerate.
If $C$ is non-degenerate,
we obtain the pairing
$C^{\lor}$ of $\nbigv_2^{\lor}$ and $\nbigv_1^{\lor}$
defined by
$C^{\lor}=C\circ
 \bigl(
 \Psi_{C}^{-1}
 \otimes
 \Psi_{C\circ\exchange}^{-1}
  \bigr)$,
  i.e., the composition of the following morphism:
\[
\begin{CD}
 \nbigv_2^{\lor}\otimes\nbigv_1^{\lor}
 @>{\Psi_C^{-1}\otimes\Psi_{C\circ\exchange}^{-1}}>>
 \nbigv_1\otimes\nbigv_2
 @>{C}>>
 \nbigo_X(\ast D).
\end{CD}   
\]

A pairing $C$ of a locally free $\nbigo_X(\ast D)$-module $\nbigv$
is a morphism $C:\nbigv\otimes\nbigv\to\nbigo_X(\ast D)$,
i.e., a pairing of $\nbigv$ and $\nbigv$.
It is called symmetric if $C\circ\exchange=C$,
i.e.,
$C(v_1\otimes v_2)=C(v_2\otimes v_1)$
for any local sections $v_i$ of $\nbigv$.
We obtain the induced pairing
$\det(C)$ of $\det(\nbigv)$.
If $C$ is non-degenerate,
we obtain the induced pairing
$C^{\lor}$ of $\nbigv^{\lor}$.

\subsubsection{Filtered bundles}

Let us recall the notion of filtered bundles
by following \cite{s1,s2}.
For any sheaf $\nbigf$ on $X$ and $Q\in X$,
let $\nbigf_Q$ denote the stalk of $\nbigf$ at $Q$.
In this subsection,
$\nbigv$ and $\nbigv_i$ $(i=1,2)$
denote locally free $\nbigo_X(\ast D)$-modules of finite rank.

For $P\in D$,
a filtered bundle
$\nbigp_{\ast}(\nbigv_P)$ over $\nbigv_P$
is an increasing sequence of
free $\nbigo_{X,P}$-submodules $\nbigp_a(\nbigv_P)$ $(a\in\real)$
of $\nbigv_P$ such that
(i) $\nbigp_a(\nbigv_P)(\ast P)=\nbigv_P$,
(ii) $\nbigp_a(\nbigv_P)=\bigcap_{b>a}\nbigp_a(\nbigv_P)$ for any $a\in\real$,
(iii) for any $a\in\real$ and $n\in\seisuu$,
we have $\nbigp_{a+n}(\nbigv_P)
=\nbigp_{a}(\nbigv_P)\otimes\nbigo_{X,P}(n\cdot P)$.
When a tuple of filtered bundles
$\nbigp_{\ast}(\nbigv_P)$ $(P\in D)$
is provided,
we obtain locally free $\nbigo_X$-submodules
$\nbigp_{\veca}(\nbigv)\subset \nbigv$ $(\veca\in\real^{D})$
characterized by
the conditions
(i) $\nbigp_{\veca}(\nbigv)\otimes_{\nbigo_X}\nbigo_X(\ast D)
=\nbigv$,
(ii) $\nbigp_{\veca}(\nbigv)_P=\nbigp_{a(P)}(\nbigv_P)$,
where $a(P)$ denotes the $P$-component of $\veca\in\real^D$.
Such a tuple
$\nbigp_{\ast}(\nbigv)=(\nbigp_{\veca}\nbigv\,|\,\veca\in\real^D)$
is called a filtered bundle over $\nbigv$.
We also say that
$\nbigp_{\ast}(\nbigv)$ is a filtered bundle on $(X,D)$.

Let $\nbigv'\subset\nbigv$ be a locally free
$\nbigo_X(\ast D)$-submodule of $\nbigv$.
When a filtered bundle $\nbigp_{\ast}(\nbigv)$ over $\nbigv$ is provided,
by setting $\nbigp_a(\nbigv'_P):=\nbigv'_P\cap\nbigp_a(\nbigv)$,
we obtain the induced filtered bundle
$\nbigp_{\ast}(\nbigv')$ over $\nbigv'$.

When filtered bundles $\nbigp_{\ast}(\nbigv_i)$ $(i=1,2)$
over $\nbigv_i$ are provided,
a morphism $f:\nbigp_{\ast}(\nbigv_1)\to\nbigp_{\ast}(\nbigv_2)$
is an $\nbigo_X(\ast D)$-homomorphism
$f:\nbigv_1\to\nbigv_2$
such that
$f(\nbigp_a(\nbigv_{1,P}))\subset\nbigp_a(\nbigv_{2,P})$
for any $P\in D$ and $a\in\real$.

When filtered bundles $\nbigp_{\ast}(\nbigv_i)$ over $\nbigv_i$
are provided,
for any $a\in\real$ and $P\in D$
we define
\[
 \nbigp_{a}\bigl(
 (\nbigv_{1}\oplus\nbigv_{2})_P
 \bigr)
 =\nbigp_a\nbigv_{1,P}\oplus\nbigp_a\nbigv_{2,P},
\quad
\nbigp_a\bigl(
(\nbigv_{1}\otimes\nbigv_{2})_P
\bigr)
=\sum_{c_1+c_2\leq a}
\nbigp_{c_1}(\nbigv_{1,P})
\otimes\nbigo_{c_2}(\nbigv_{2,P}).
\]
We also define
\[
 \nbigp_a\bigl(
 \nhom_{\nbigo_X(\ast D)}(\nbigv_1,\nbigv_2)_P
 \bigr)
 =\bigl\{
 f\in \nhom_{\nbigo_{X}(\ast D)_P}
 (\nbigv_{1,P},\nbigv_{2,P})\,\big|\,
 f(\nbigp_b(\nbigv_{1,P}))
 \subset
 \nbigp_{b+a}(\nbigv_{2,P})\,\,
 \forall b\in\real
 \bigr\}.
\]
In this way,
we obtain filtered bundles
$\nbigp_{\ast}(\nbigv_1\oplus\nbigv_2)$
over $\nbigv_1\oplus\nbigv_2$,
$\nbigp_{\ast}(\nbigv_1\otimes\nbigv_2)$
over $\nbigv_1\otimes\nbigv_2$,
and 
$\nbigp_{\ast}\nhom(\nbigv_1,\nbigv_2)$
over $\nhom(\nbigv_1,\nbigv_2)$.
The filtered bundles
are denoted by
$\nbigp_{\ast}(\nbigv_1)\oplus\nbigp_{\ast}(\nbigv_2)$,
$\nbigp_{\ast}(\nbigv_1)\otimes\nbigp_{\ast}(\nbigv_2)$,
and $\nhom(\nbigp_{\ast}\nbigv_1,\nbigp_{\ast}\nbigv_2)$,
respectively.

\begin{rem}
Even if $\nbigv=\nbigv'\oplus\nbigv''$ holds
as locally free $\nbigo_X(\ast D)$-modules,
$\nbigp_{\ast}(\nbigv)
=\nbigp_{\ast}(\nbigv')\oplus\nbigp_{\ast}(\nbigv'')$
does not necessarily hold 
for the induced filtered bundles
$\nbigp_{\ast}(\nbigv')$  and $\nbigp_{\ast}(\nbigv'')$.
If $\nbigp_{\ast}(\nbigv)
=\nbigp_{\ast}(\nbigv')\oplus\nbigp_{\ast}(\nbigv'')$ holds,
we say that
the filtered bundle $\nbigp_{\ast}(\nbigv)$ is compatible
with the decomposition $\nbigv=\nbigv'\oplus\nbigv''$.
\hfill\qed
\end{rem}

We set $r=\rank\nbigv$.
From a filtered bundle $\nbigp_{\ast}\nbigv$ over $\nbigv$,
we obtain
a filtered bundle
$\nbigp_{\ast}(\nbigv)^{\otimes r}$
over $\nbigv^{\otimes r}$.
There exists the decomposition
$\nbigv^{\otimes r}=
\det(\nbigv)\oplus(\nbigv^{\otimes r})^{\bot}$
in a way compatible with the natural action of
the $r$-th symmetric group $\gbigs_{r}$.
We obtain
the filtered bundles
$\nbigp_{\ast}(\det(\nbigv))$
and
$\nbigp_{\ast}((\nbigv^{\otimes\,r})^{\bot})$
over $\det(\nbigv)$
and $(\nbigv^{\otimes\,r})^{\bot}$,
respectively,
for which
$\nbigp_{\ast}(\nbigv)^{\otimes r}
=\nbigp_{\ast}(\det(\nbigv))
\oplus
\nbigp_{\ast}((\nbigv^{\otimes\,r})^{\bot})$
holds.

Let $\nbigp^{(0)}_{\ast}(\nbigo_X(\ast D))$
denote the canonical filtered bundle over $\nbigo_X(\ast D)$,
i.e.,
\[
 \nbigp^{(0)}_{\veca}(\nbigo_X(\ast D))
=\nbigo_X\Bigl(\sum_{P\in D}[a(P)]P\Bigr),
\]
where we set $[c]:=\max\{n\in\seisuu\,|\,n\leq c\}$
for any $c\in\real$.
When a filtered bundle $\nbigp_{\ast}(\nbigv)$ over $\nbigv$
is provided,
let $\nbigp_{\ast}(\nbigv^{\lor})$
denote
$\nhom(\nbigp_{\ast}\nbigv,
\nbigp^{(0)}_{\ast}(\nbigo_X(\ast D)))$.

\subsubsection{Pairings of filtered bundles}

Let $\nbigp_{\ast}(\nbigv_i)$ $(i=1,2)$
be filtered bundles on $(X,D)$.
A pairing $C$ of
$\nbigp_{\ast}(\nbigv_1)$
and $\nbigp_{\ast}(\nbigv_2)$
is a morphism
\[
C:\nbigp_{\ast}(\nbigv_1)
\otimes\nbigp_{\ast}(\nbigv_2)
\lrarr
 \nbigp_{\ast}^{(0)}(\nbigo_X(\ast D)).
\]
We obtain the pairing
$C\circ\exchange$
of $\nbigp_{\ast}(\nbigv_2)$
and $\nbigp_{\ast}(\nbigv_1)$.
If $\rank\nbigv_1=\rank\nbigv_2$,
we obtain
\begin{equation}
\label{eq;22.9.3.50}
 \det(C):
 \det(\nbigp_{\ast}\nbigv_1)
 \otimes
 \det(\nbigp_{\ast}\nbigv_2)
\lrarr
 \nbigp^{(0)}_{\ast}(\nbigo_X(\ast D)).
\end{equation}
From a pairing $C$ of $\nbigp_{\ast}(\nbigv_1)$
and $\nbigp_{\ast}(\nbigv_2)$,
we obtain the following morphism of filtered bundles:
\begin{equation}
\label{eq;22.8.29.1}
\Psi_C:\nbigp_{\ast}(\nbigv_1)
\lrarr\nbigp_{\ast}(\nbigv_2^{\lor}).
\end{equation}
\begin{df}
$C$ is called perfect
if the morphism {\rm(\ref{eq;22.8.29.1})}
is an isomorphism of filtered bundles. 
\hfill\qed
\end{df}
Note that if $\rank\nbigv=1$,
$C$ is perfect if and only if (\ref{eq;22.9.3.50})
is an isomorphism.

\begin{lem}
\label{lem;22.8.26.11}
 $C$ is perfect if and only if
the following induced morphism is an isomorphism:
\begin{equation}
\label{eq;22.8.29.2}
 \det(C):
  \det(\nbigp_{\ast}\nbigv_1)
  \otimes
  \det(\nbigp_{\ast}\nbigv_2)
  \lrarr
  \nbigp^{(0)}_{\ast}(\nbigo_X(\ast D)).
\end{equation}
\end{lem}
\pf
The morphism (\ref{eq;22.8.29.1}) is an isomorphism
if and only if
the induced morphism
\[
 \det(\Psi_C):
 \det(\nbigp_{\ast}\nbigv_1)
 \lrarr
 \det(\nbigp_{\ast}\nbigv_2^{\lor})
\]
is an isomorphism,
which is equivalent to the condition that
(\ref{eq;22.8.29.2}) is an isomorphism.
\hfill\qed

\vspace{.1in}

Let $\nbigv_i'\subset\nbigv_i$
be locally free $\nbigo_X(\ast D)$-submodules.
For simplicity, we also assume that
$\nbigv_i'$ are saturated,
i.e.,
$\nbigv_i/\nbigv_i'$ are also locally free.
From a pairing $C$ of $\nbigp_{\ast}\nbigv_1$
and $\nbigp_{\ast}(\nbigv_2)$,
we obtain the induced pairing $C'$
of $\nbigp_{\ast}\nbigv_1'$ and $\nbigp_{\ast}\nbigv_2'$.
There exist the following natural morphisms:
\begin{equation}
 \label{eq;22.8.26.2}
\begin{CD}
 \nbigv_1' @>{i_1}>>
 \nbigv_1 @>{\Psi_C}>>
 \nbigv_2^{\lor}
 @>{i_2^{\lor}}>>
 (\nbigv_2')^{\lor}.
\end{CD}
\end{equation}
Here, $i_1$ denotes the natural inclusion,
and $i_2^{\lor}$ denotes the dual of
the natural inclusion $i_2:\nbigv_2'\to\nbigv_2$.
Note that $\Psi_{C'}=i_2^{\lor}\circ\Psi_C\circ i_1$.
Let $\nbigu_1\subset\nbigv_1$ denote the kernel of
$i_2^{\lor}\circ\Psi_C$.
We have the induced filtered bundle
$\nbigp_{\ast}\nbigu_1$ over $\nbigu_1$.
The following lemma is obvious.

\begin{lem}
\label{lem;22.8.26.10}
If $C$ and $C'$ are perfect,
we have the decomposition of the filtered bundles
$\nbigp_{\ast}\nbigv_1=
 \nbigp_{\ast}\nbigv_1'
 \oplus
\nbigp_{\ast}\nbigu_1$.
\hfill\qed
\end{lem}

\subsubsection{Symmetric pairings of filtered bundles}

Let $C$ be a symmetric pairing of
a filtered bundle $\nbigp_{\ast}\nbigv$ on $(X,D)$.
We have the induced pairing
\begin{equation}
\label{eq;22.8.27.11}
 \det(C):
  \det(\nbigp_{\ast}\nbigv)
  \otimes
  \det(\nbigp_{\ast}\nbigv)
  \to
  \nbigp^{(0)}_{\ast}(\nbigo_X(\ast D)).
\end{equation}
We obtain the following lemma as a special case of
Lemma \ref{lem;22.8.26.11}.
\begin{lem}
\label{lem;22.8.27.31}
$C$ is perfect if and only if
{\rm(\ref{eq;22.8.27.11})} is an isomorphism.
 \hfill\qed
\end{lem}

Let $\nbigv'\subset\nbigv$ be a saturated locally free
$\nbigo_X(\ast D)$-submodule.
Let $(\nbigv')^{\bot\,C}$
denote the kernel of the composition of the following morphisms:
\[
 \begin{CD}
  \nbigv@>{\Psi_C}>>
  \nbigv^{\lor}
  @>{i^{\lor}}>> (\nbigv')^{\lor},
 \end{CD}
\]
where $i^{\lor}$ denote the dual of
the inclusion $i:\nbigv'\to\nbigv$.
We have the filtered bundle
$\nbigp_{\ast}(\nbigv')^{\bot\,C}$
over $(\nbigv')^{\bot\,C}$.
Let $C'$ denote the induced symmetric pairing
of $\nbigp_{\ast}(\nbigv')$.
The following lemma is a special case of
Lemma \ref{lem;22.8.26.10}.

\begin{lem}
\label{lem;22.8.27.32}
If $C$ and $C'$ are perfect,
we obtain the decomposition
$\nbigp_{\ast}\nbigv
 =\nbigp_{\ast}\nbigv'
 \oplus
 \nbigp_{\ast}
 (\nbigv')^{\bot\,C}$.
\hfill\qed
\end{lem}

\begin{cor}
Suppose that
{\rm(\ref{eq;22.8.27.11})} is an isomorphism.
We also assume that
\begin{equation}
 \label{eq;22.8.27.21}
  \det(C'):
  \det(\nbigp_{\ast}\nbigv')
  \otimes
   \det(\nbigp_{\ast}\nbigv')
\to\nbigp^{(0)}_{\ast}(\nbigo_X(\ast D))
\end{equation}
is an isomorphism. 
Then, we have the decomposition of the filtered bundles
$\nbigp_{\ast}\nbigv=
 \nbigp_{\ast}\nbigv'
 \oplus
 \nbigp_{\ast}(\nbigv')^{\bot\,C}$.
\hfill\qed
\end{cor}

\subsubsection{Compact case}

We assume that $X$ is compact.
Recall that
for any filtered bundle $\nbigp_{\ast}\nbigv$ on $(X,D)$,
we obtain $\deg(\nbigp_{\ast}\nbigv)\in\real$
as follows
(see \cite{s1,s2} and also \cite{Mochizuki-KH-Higgs}):
\[
 \deg(\nbigp_{\ast}\nbigv)
 =\deg(\nbigp_{\veca}\nbigv)
 -\sum_{P\in D}
 \sum_{a(P)-1<a\leq a(P)}
 a\cdot
 \dim_{\cnum}\Gr^{\nbigp}_a(\nbigv_P).
\]
Here, we set
$\Gr^{\nbigp}_a(\nbigv_P):=
\nbigp_{a}(\nbigv_P)\Big/\sum_{b<a}\nbigp_b(\nbigv_P)$.
It is independent of the choice of $\veca\in\real^D$.
Note that
$\deg(\nbigp_{\ast}\nbigv)
=\deg(\det(\nbigp_{\ast}\nbigv))$.
\begin{lem}
\label{lem;22.8.28.4}
Let $\nbigp_{\ast}\nbigv_i$ $(i=1,2)$
be filtered bundles of rank $1$ on $(X,D)$.
If there exists a non-zero morphism
$f:\nbigp_{\ast}\nbigv_1\to\nbigp_{\ast}\nbigv_2$,
then $\deg(\nbigp_{\ast}\nbigv_1)\leq\deg(\nbigp_{\ast}\nbigv_2)$.
If moreover $\deg(\nbigp_{\ast}\nbigv_1)=\deg(\nbigp_{\ast}\nbigv_2)$,
then $f$ is an isomorphism.
\end{lem}
\pf
Though this is well known,
we include a sketch of the proof for the convenience of the readers.
By the morphism $f$,
we can regard $\nbigv_1$ as an $\nbigo_X(\ast D)$-submodule of
$\nbigv_2$.
There exists a finite subset $Z\subset X\setminus D$
and function $m:Z\to \seisuu_{>0}$
such that
$\nbigv_1=\nbigv_2(-\sum_{Q\in Z}m(Q)Q )$.
For each $P\in D$,
we take a non-zero element $v_P\in\nbigv_{1,P}=\nbigv_{2,P}$.
We obtain the numbers
$b(P,i):=\min\{b\in\real\,|\,v_P\in\nbigp_{b}(\nbigv_{i,P})\}$.
Let $\nbigv'_i\subset\nbigv_i$
be the locally free $\nbigo_X$-submodules
determined by the conditions
(i) $\nbigv'_i\otimes\nbigo_X(\ast D)=\nbigv_i$,
(ii) $\nbigv'_{i,P}=\nbigo_{X,P}v_P$.
By the definition, we obtain
\[
 \deg(\nbigp_{\ast}\nbigv_i)
=\deg(\nbigv_i')-\sum_{P\in D} b(P,i).
\]
Because $\nbigv_1'=\nbigv_2'(-\sum m(Q)Q)$
and $b(P,1)\geq b(P,2)$ $(P\in D)$,
we obtain
\[
\deg(\nbigp_{\ast}\nbigv_1)
=\deg(\nbigv_1')-\sum_{P\in D}b(P,1)
=\deg(\nbigv_2')-\sum_{Q\in Z}m(Q)
-\sum_{P\in D}b(P,1)
\leq
\deg(\nbigv_2')-\sum_{P\in D}b(P,2)=\deg(\nbigp_{\ast}\nbigv_2).
\]
If $\deg(\nbigp_{\ast}\nbigv_1)=\deg(\nbigp_{\ast}\nbigv_2)$,
we obtain $Z=\emptyset$ and $b(P,1)=b(P,2)$ $(P\in D)$,
and hence $\nbigp_{\ast}\nbigv_1=\nbigp_{\ast}\nbigv_2$.
\hfill\qed

\begin{prop}
\label{prop;22.8.28.3}
Let $C$ be a symmetric pairing of
a filtered bundle $\nbigp_{\ast}\nbigv$ on $(X,D)$.
Then, we have either
(i) $\det(C)=0$,
or (ii) $\det(\nbigp_{\ast}\nbigv)\leq 0$.
If $\det(C)\neq 0$ and $\deg(\nbigp_{\ast}\nbigv)=0$,
then $C$ is perfect.
\end{prop}
\pf
We obtain
\begin{equation}
\label{eq;22.8.28.2}
 \det(C):\det(\nbigp_{\ast}\nbigv)
 \otimes\det(\nbigp_{\ast}\nbigv)
 \to \nbigp^{(0)}_{\ast}(\nbigo_X(\ast D)).
\end{equation}
By Lemma \ref{lem;22.8.28.4},
if $\det(C)\neq 0$,
we obtain 
\[
2\deg(\nbigp_{\ast}\nbigv)
=\deg\bigl(
\det(\nbigp_{\ast}\nbigv)
\otimes
\det(\nbigp_{\ast}\nbigv)
\bigr)
\leq \deg(\nbigp^{(0)}_{\ast}(\nbigo_X(\ast D)))=0.
\]
If moreover $\deg(\nbigp_{\ast}\nbigv)=0$ holds,
then 
(\ref{eq;22.8.28.2}) is an isomorphism.
Hence, $C$ is perfect
by Lemma \ref{lem;22.8.27.31}.
\hfill\qed

\vspace{.1in}
As a complement,
we remark the following.
\begin{lem}
\label{lem;22.9.7.12}
If a filtered bundle $\nbigp_{\ast}\nbigv$ on $(X,D)$
has a perfect symmetric pairing,
then we obtain $\deg(\nbigp_{\ast}\nbigv)=0$.
\end{lem}
\pf
Because
$\det(\nbigp_{\ast}\nbigv)\otimes\det(\nbigp_{\ast}\nbigv)
\simeq \nbigp^{(0)}(\nbigo_X(\ast D))$,
we obtain
$\deg(\nbigp_{\ast}\nbigv)
=\deg(\det(\nbigp_{\ast}\nbigv))
=0$.
\hfill\qed

\subsubsection{Harmonic metrics in the rank one case}

We continue to assume that $X$ is compact.
Let $h^{(0)}$ denote the harmonic metric of
$\nbigo_{X\setminus D}$
defined by $h^{(0)}(1,1)=1$.
It is adapted to $\nbigp_{\ast}^{(0)}(\nbigo_X(\ast D))$,
i.e.,
$\nbigp^{h^{(0)}}_{\ast}(\nbigo_{X\setminus D})
=\nbigp^{(0)}_{\ast}(\nbigo_X(\ast D))$.

Let $\nbigp_{\ast}\nbigv$ be a filtered bundle on $(X,D)$
of rank $1$ such that $\deg(\nbigp_{\ast}\nbigv)=0$
equipped with a non-zero pairing:
\begin{equation}
\label{eq;22.8.29.52}
 C:\nbigp_{\ast}(\nbigv)
  \otimes
  \nbigp_{\ast}(\nbigv)
  \lrarr \nbigp^{(0)}_{\ast}(\nbigo_X(\ast D)).
\end{equation}
Note that (\ref{eq;22.8.29.52}) is an isomorphism
by Proposition \ref{prop;22.8.28.3}.

\begin{lem}
\label{lem;22.8.29.51}
There uniquely exists a Hermitian metric $h$ of
$\nbigv_{|X\setminus D}$
such that
(i) $h$ is harmonic, i.e., the Chern connection of $h$ is flat,
(ii) $h$ is adapted to $\nbigp_{\ast}\nbigv$,
(iii) $h\otimes h=h^{(0)}$ under
the isomorphism
$C:\nbigp_{\ast}\nbigv\otimes\nbigp_{\ast}\nbigv\simeq
\nbigp_{\ast}^{(0)}(\nbigo_X(\ast D))$. 
Note that the condition (iii) means that
$h$ is compatible with $C_{|X\setminus D}$.
\end{lem}
\pf
It is well known that 
there exists a harmonic metric $h_1$
of $\nbigo_{X\setminus D}$
satisfying the conditions (i) and (ii).
We obtain the metric $h_1\otimes h_1$
of $\nbigo_{X\setminus D}$ by the isomorphism
(\ref{eq;22.8.29.52}),
which is adapted to
$\nbigp^{(0)}_{\ast}(\nbigo_X(\ast D))$.
By the uniqueness of harmonic metrics,
there exists $c>0$
such that
$h_1\otimes h_1=ch^{(0)}$.
Then, $h=c^{-1/2}h_1$ satisfies the desired conditions.

\hfill\qed

\subsection{Symmetric pairings of good filtered Higgs bundles}

\subsubsection{Preliminary}

Let $X$ be a Riemann surface.
For each $P\in D$,
let $(X_P,z_P)$ be a holomorphic coordinate neighbourhood
around $P$ such that $z_P(P)=0$.
By the coordinate, we may regard $X_P$ as
a neighbourhood of $0$ in $\cnum$.
For a positive integer $e$,
let $\varphi_e:\cnum\to\cnum$ be defined by
$\varphi_e(\zeta)=\zeta^e$.
Let $X^{(e)}_P=\varphi_e^{-1}(X_P)$.
Let $\varphi_{P,e}:X^{(e)}_P\to X_P$
denote the induced morphism.
Let $P^{(e)}\in X^{(e)}_P$ denote the inverse image of $P$.
On $X^{(e)}_P$,
we choose a holomorphic function $z_{P,e}$
such that 
$z_{P,e}^e=\varphi_{P,e}^{\ast}(z_P)$.
Let $\Gal^{(e)}_{P}$ be the Galois group
of the ramified covering $X^{(e)}_P\to X_P$.

\subsubsection{Meromorphic Higgs bundles}

Let $D\subset X$ be a discrete subset.
Let $\nbigv$ be a locally free $\nbigo_X(\ast D)$-module
of rank $r$.
A Higgs field of $\nbigv$
is a morphism $\theta:\nbigv\to\nbigv\otimes\Omega^1_X$.
Such a tuple $(\nbigv,\theta)$ is called a meromorphic Higgs
bundle on $(X,D)$.
We obtain the Higgs field
$\theta^{\lor}$ of $\nbigv^{\lor}$
as the dual of $\theta$.
We also obtain the Higgs field
$\tr(\theta)$ of $\det(\nbigv)$.
A morphism $g:(\nbigv_1,\theta_1)\to(\nbigv_2,\theta_2)$
of a meromorphic Higgs bundles
is an $\nbigo_X(\ast D)$-homomorphism $g:\nbigv_1\to\nbigv_2$
such that
$\theta_2\circ g=g\circ\theta_1$.

Let $(\nbigv,\theta)$ be a meromorphic Higgs bundle on $(X,D)$.
Let $P\in D$.
We obtain the endomorphism $f_P$
by $\theta=f_P\,dz_P/z_P$ around $P$.
We say that $(\nbigv,\theta)$ is regular at $P$
if there exists
a free $\nbigo_{X,P}$-submodule
$\nbigl\subset\nbigv_{P}$
such that
(i) $\nbigl(\ast P)=\nbigv_{P}$,
(ii) $\theta$ is logarithmic with respect to
$\nbigl$,
i.e.,
$\theta(\nbigl)\subset z_P^{-1}\nbigl\otimes\Omega^1_{X,P}$.
The second condition is equivalent to
that $f_P(\nbigl)\subset\nbigl$.
Note that we have the characteristic polynomial
$\det(T\id_{\nbigv}-f_P)=T^{r}+\sum_{j=0}^{r-1} a_j(z_P)T^j$.
\begin{lem}
$(\nbigv,\theta)$ is regular at $P$
 if and only if $a_j$ are holomorphic at $P$.
\end{lem}
\pf
The ``only if'' part is obvious.
Suppose that $a_j$ are holomorphic at $P$.
We choose a non-zero $v\in\nbigv_P$.
We consider the lattice $\nbigl$ generated by
$f^j(v)$ $(j=0,\ldots,r-1)$.
Because
$f^{r}+\sum_{j=0}^{r-1}a_j(z_P)f^j=0$,
we obtain $f^{r}(v)\in\nbigl$.
Hence, $f(\nbigl)\subset\nbigl$.
\hfill\qed

\vspace{.1in}
The following lemma is well known and easy to prove.

\begin{lem}
If $(\nbigv,\theta)$ is regular at $P$,
there exists a decomposition
\begin{equation}
\label{eq;22.8.29.40}
 (\nbigv_P,\theta)
 =\bigoplus_{\alpha\in\cnum}
 (\nbigv_{P,\alpha},\theta_{P,\alpha}) 
\end{equation}
such that the following holds.
\begin{itemize}
 \item Let $f_{P,\alpha}$ be the endomorphism satisfying
       $\theta_{P,\alpha}-\alpha\,\id_{\nbigv_{P,\alpha}}dz_P/z_P
       =f_{P,\alpha}dz_P/z_P$.
       Let $\det(T\id-f_{P,\alpha})=\sum a_{\alpha,j}(z_P)T^j$
       denote the characteristic polynomial of $f_{P,\alpha}$.
       Then, $a_{\alpha,j}$ are holomorphic at $P$,
       and $a_{\alpha,j}(P)=0$ unless $j=\rank\nbigv_{P,\alpha}$.
\hfill\qed
\end{itemize}
\end{lem}

We recall the following lemma,
which is an analogue of the Hukuhara-Levelt-Turrittin theorem
for meromorphic flat bundles,
but easier to prove by using only standard arguments
in linear algebra.

\begin{lem}
There exist $e\in\seisuu_{>0}$
and a decomposition
\begin{equation}
\label{eq;22.8.29.41}
 (\varphi_{P,e}^{\ast}(\nbigv)_{P^{(e)}},
 \varphi_{P,e}^{\ast}(\theta))
 =\bigoplus_{\gminia\in z_{P,e}^{-1}\cnum[z_{P,e}^{-1}]}
 (\nbigv^{(e)}_{P,\gminia},\theta^{(e)}_{P,\gminia})
 \end{equation}
such that
$(\nbigv^{(e)}_{P,\gminia},\theta^{(e)}_{P,\gminia}-d\gminia\id)$
are regular. 
Indeed, $e$ divides $r!$.
\hfill\qed
\end{lem}

\subsubsection{Good filtered Higgs bundles}

Let $(\nbigv,\theta)$ be a meromorphic Higgs bundle on $(X,D)$.
Let $\nbigp_{\ast}\nbigv$ be a filtered bundle over $\nbigv$.
A filtered Higgs bundle $(\nbigp_{\ast}\nbigv,\theta)$
is called regular at $P\in D$
if $\theta$ is logarithmic with respect to 
each $\nbigp_a(\nbigv_P)$ $(a\in\real)$.
If $(\nbigp_{\ast}\nbigv,\theta)$ is regular at $P$,
the decomposition (\ref{eq;22.8.29.40})
is compatible with the filtration,
i.e.,
we obtain the decomposition of filtered bundles at $P$
for each $a\in\real$:
\[
 \nbigp_a(\nbigv_P)
 =\bigoplus_{\alpha\in\cnum}
 \nbigp_a(\nbigv_{P,\alpha}).
\] 

Recall that
for any $e\in\seisuu_{>0}$,
we obtain the filtered bundle
$\nbigp_{\ast}\bigl(
\varphi_{P,e}^{\ast}(\nbigv)_{P^{(e)}}\bigr)$
over $\varphi_{P,e}^{\ast}(\nbigv)_{P^{(e)}}$
defined as follows:
\[
 \nbigp_a\Bigl(
 \varphi_{P,e}^{\ast}(\nbigv)_{P^{(e)}}
 \Bigr)
 =\sum_{n+eb\leq a}
 z_{P,e}^{-n}
 \varphi_{P,e}^{\ast}\bigl(
 \nbigp_b(\nbigv_P)
 \bigr)
\]
A filtered Higgs bundle $(\nbigp_{\ast}\nbigv,\theta)$
is called good at $P\in D$
if the decomposition 
(\ref{eq;22.8.29.41}) is compatible with
the filtration,
i.e.,
\[
 \nbigp_{a}\Bigl(
 \varphi_{P,e}^{\ast}(\nbigv)_{P^{(e)}}
 \Bigr)
 =\bigoplus
 \nbigp_{a}(\nbigv^{(e)}_{P,\gminia}),
\]
and moreover 
$\theta_{P,\gminia}-d\gminia\id$ are logarithmic
with respect to
$\nbigp_a(\nbigv^{(e)}_{P,\gminia})$
for any $a\in\real$.

\subsubsection{Symmetric pairings of good filtered Higgs bundles}

\begin{df}
 A symmetric pairing $C$ of
a good filtered Higgs bundle $(\nbigp_{\ast}\nbigv,\theta)$
(resp. a meromorphic Higgs bundle $(\nbigv,\theta)$)
is a symmetric pairing $C$ of $\nbigp_{\ast}\nbigv$
(resp. $\nbigv$)
such that $C(\theta\otimes\id)=C(\id\otimes\theta)$. 
 \hfill\qed
\end{df}

In case,
we obtain the induced morphism
$\Psi_C:(\nbigp_{\ast}\nbigv,\theta)
\to (\nbigp_{\ast}\nbigv^{\lor},\theta^{\lor})$
(resp.
$\Psi_C:(\nbigv,\theta)\to(\nbigv^{\lor},\theta^{\lor})$).
We also obtain a symmetric pairing
$\det(C)$ of $(\nbigp_{\ast}\nbigv,\tr\theta)$.

\subsubsection{Wild harmonic bundle with a real structure}
\label{subsection;22.9.3.60}

Let $(E,\delbar_E,\theta,h)$ be a wild harmonic bundle on $(X,D)$.
Around $P\in D$, let $z_P$ denote a holomorphic local coordinate
such that $z_P(P)=0$.
For an open neighbourhood $U$ of $P$,
let $\nbigp^h_a(E)(U)$ denote the space of
holomorphic sections $s$ of $E_{|U\setminus\{P\}}$
satisfying $|s|_h=O(|z_P|^{-a-\epsilon})$ for any $\epsilon>0$.
In this way,
we obtain the associated good filtered Higgs bundle
$(\nbigp^h_{\ast}(E),\theta)$ on $(X,D)$
\cite{s2,Mochizuki-wild}.
Let $C$ be a real structure of the harmonic bundle
$(E,\delbar_E,\theta,h)$.
(See \S\ref{subsection;22.9.23.11}.)

\begin{lem}
\label{lem;22.9.5.5}
$C$ induces a perfect symmetric pairing of
$(\nbigp^h_{\ast}E,\theta)$.
\end{lem}
\pf
Let $h^{\lor}$ be the induced metric of $E^{\lor}$.
The dual of
$\nbigp^{h}_{\ast}(E)$
is naturally isomorphic to
$\nbigp^{h^{\lor}}_{\ast}(E^{\lor})$.
Because $\Psi_C$ extends to an isomorphism
$\nbigp^h_{\ast}(E)\simeq\nbigp^{h^{\lor}}_{\ast}(E^{\lor})$,
we obtain the claim of the lemma.
\hfill\qed

\vspace{.1in}
We recall the following proposition \cite{Mochizuki-wild}.
\begin{prop}
\label{prop;22.9.7.10}
Suppose that $X$ is a compact Riemann surface.
\begin{itemize}
 \item 
 $(\nbigp_{\ast}^h(E),\theta)$
is polystable of degree $0$.
\item
 If $\nbigp^{h_1}_{\ast}(E)=\nbigp^{h_2}_{\ast}(E)$
for $h_i\in\Harm(E,\delbar_E,\theta;C)$,
there exists a decomposition
$(\nbigp^{h_1}_{\ast}(E),\theta)
=\bigoplus_{i=1}^m
(\nbigp_{\ast}\nbigv_i,\theta_i)$
such that
(i) $E=\bigoplus \nbigv_{i|X\setminus D}$ 
is orthogonal with respect to both $h_1$ and $h_2$,
(ii) $h_2=a_ih_1$ for some $a_i>0$ on $\nbigv_{i|X\setminus D}$.
In particular, 
$h_1$ and $h_2$ are mutually bounded,
and we have $\delbar_Es(h_1,h_2)=0$ and $[\theta,s(h_1,h_2)]=0$.
\hfill\qed
\end{itemize}
\end{prop}

We obtain the following corollary
from Proposition \ref{prop;22.9.7.11}.
\begin{cor}
In Proposition {\rm\ref{prop;22.9.7.10}},
the $\GL(n,\real)$-harmonic bundles
associated with $(E,\delbar_E,\theta,h_i;C)$
are naturally isomorphic.
\hfill\qed
\end{cor}

\subsection{Kobayashi-Hitchin correspondence}

Suppose that $X$ is compact.

\subsubsection{Basic polystable objects (1)}

Let $(\nbigp_{\ast}\nbigv,\theta)$
be a stable good filtered Higgs bundle 
of degree $0$ on $(X,D)$
such that
$(\nbigp_{\ast}\nbigv,\theta)
\simeq
(\nbigp_{\ast}\nbigv^{\lor},\theta^{\lor})$.
Let $P_1$ be a pairing
\[
 \nbigp_{\ast}\nbigv
 \otimes
 \nbigp_{\ast}\nbigv
 \to \nbigp^{(0)}_{\ast}(\nbigo_X(\ast D))
\]
such that $\Psi_{P_1}$ induces an isomorphism
$(\nbigp_{\ast}\nbigv,\theta)
\simeq
(\nbigp_{\ast}\nbigv^{\lor},\theta^{\lor})$.
Note that if $P_1'$ is such another pairing,
then there exists $\alpha\in\cnum$
such that $P_1'=\alpha P_1$
by the stability assumption.

\begin{lem}
Either one of
$P_1\circ \exchange=P_1$
or 
$P_1\circ \exchange=-P_1$ holds.
\end{lem}
\pf
Because of the stability condition,
there exists $\beta\in\cnum$
such that
$\Psi_{P_1}^{\lor}=\beta\Psi_{P_1}$.
Because 
$(\Psi_{P_1}^{\lor})^{\lor}=\Psi_{P_1}$,
we obtain $\beta^2=1$.
Because
$\Psi_{P_1\circ\exchange}=\Psi_{P_1}^{\lor}$,
we obtain the claim of the lemma.
\hfill\qed

\vspace{.1in}

Let $C_{\cnum^{\ell}}$
denote the symmetric bilinear form of $\cnum^{\ell}$
defined by
$C_{\cnum^{\ell}}(\vecx,\vecy)=\sum_{i=1}^{\ell} x_iy_i$
for $\vecx=(x_i)$ and $\vecy=(y_i)$.
Let $\omega_{\cnum^{2k}}$
denote the skew-symmetric bilinear form of $\cnum^{2k}$
defined by
$\omega_{\cnum^{2k}}(\vecx,\vecy)
=\sum_{i=1}^k(x_{2i-1}y_{2i}-x_{2i}y_{2i-1})$.
If $P_1$ is symmetric,
we obtain a perfect symmetric pairing
$P_1\otimes C_{\cnum^{\ell}}$
of $(\nbigp_{\ast}\nbigv,\theta)\otimes\cnum^{\ell}$.
If $P_1$ is skew-symmetric,
we obtain a perfect symmetric pairing
$P_1\otimes \omega_{\cnum^{2k}}$
of $(\nbigp_{\ast}\nbigv,\theta)\otimes\cnum^{2k}$.

\begin{lem}
\label{lem;22.9.4.10}
Suppose that
$(\nbigp_{\ast}\nbigv,\theta)\otimes\cnum^{\ell}$
is equipped with a perfect symmetric pairing $C$.
\begin{itemize}
 \item If $P_1$ is symmetric,
       there exists an automorphism $\rho$ of
       $\cnum^{\ell}$
       such that
       $(\id_{\nbigv}\otimes\rho)^{\ast}C=
       P_1\otimes C_{\cnum^{\ell}}$.
 \item If $P_1$ is skew-symmetric,
       $\ell$ is an even number $2k$,
       and there exists an automorphism
       $\rho$ of $\cnum^{2k}$
       such that
       $(\id_{\nbigv}\otimes\rho)^{\ast}C=P_1\otimes\omega_{\cnum^{2k}}$.

\end{itemize} 
\end{lem}
\pf
There exists a non-degenerate bilinear form
$C_1$ of $\cnum^{\ell}$ such that $C=P_1\otimes C_1$.
If $P_1$ is symmetric,
then $C_1$ is symmetric.
By using an orthonormal frame of $\cnum^{\ell}$
with respect to $C_1$,
we obtain the first claim.
If $P_1$ is skew-symmetric,
$C_1$ is skew-symmetric.
In particular, $\ell=2k$ for a positive integer $k$.
By using a symplectic base of $\cnum^{2k}$,
we obtain the second claim.
\hfill\qed

\begin{lem}
There exists a unique harmonic metric $h_0$
of $(\nbigv,\theta)_{|X\setminus D}$
such that
(i) $\Psi_{P_1}$ is isometric with respect to
$h_0$ and $h_0^{\lor}$,
(ii) $h_0$ is adapted to $\nbigp_{\ast}\nbigv$.
\end{lem}
\pf
Let $h_1$ be a harmonic metric
of $(\nbigv,\theta)_{|X\setminus D}$
which is adapted to $\nbigp_{\ast}\nbigv$.
Let $h_1^{\lor}$ denote the induced harmonic metric of
$(\nbigv^{\lor},\theta^{\lor})_{|X\setminus D}$,
which is adapted to $\nbigp_{\ast}\nbigv^{\lor}$.
Note that
both $(\nbigp_{\ast}\nbigv,\theta)$
and $(\nbigp_{\ast}\nbigv^{\lor},\theta^{\lor})$
are stable of degree $0$,
and that
$\Psi_{P_1}:(\nbigp_{\ast}\nbigv,\theta)
\lrarr (\nbigp_{\ast}\nbigv^{\lor},\theta^{\lor})$
is an isomorphism.
Hence, there exists $c>0$ such that
$\Psi_{P_1}^{\ast}(h_1^{\lor})=c^2h_1$.
We set $h_0=ch_1$.
Because $h_0^{\lor}=c^{-1}h_1^{\lor}$,
$h_0$ has the desired property.
The uniqueness is also clear.
\hfill\qed

\vspace{.1in}
The following lemma follows from 
the uniqueness of harmonic metrics
adapted to parabolic structure.
\begin{lem}
\mbox{{}}\label{lem;22.9.4.21}
\begin{itemize}
 \item 
For any Hermitian metric $h_{\cnum^{\ell}}$ of $\cnum^{\ell}$,
$h_0\otimes h_{\cnum^{\ell}}$
is a harmonic metric of
$(\nbigv,\theta)_{|X\setminus D}\otimes\cnum^{\ell}$
which is adapted to $\nbigp_{\ast}\nbigv\otimes\cnum^{\ell}$.
Conversely,  
for any harmonic metric $h$ of
$(\nbigv,\theta)_{|X\setminus D}\otimes\cnum^{\ell}$
which is adapted to $\nbigp_{\ast}\nbigv$,
there exists a Hermitian metric $h_{\cnum^{\ell}}$ of $\cnum^{\ell}$
such that $h=h_0\otimes h_{\cnum^{\ell}}$.
 \item
If $P_1$ is symmetric (resp. skew-symmetric),
a harmonic metric $h_0\otimes h_{\cnum^{\ell}}$
of $(\nbigv,\theta)_{|X\setminus D}\otimes\cnum^{\ell}$
is compatible with $P_1\otimes C_{\cnum^{\ell}}$
(resp. $P_1\otimes \omega_{\cnum^{\ell}}$)
if and only if  $h_{\cnum^{\ell}}$ is compatible with
$C_{\cnum^{\ell}}$
(resp. $\omega_{\cnum^{\ell}}$).
\hfill\qed
\end{itemize}
\end{lem}

See {\rm\S\ref{subsection;22.9.4.4}}
and {\rm\S\ref{subsection;22.9.4.2}}
for the ambiguity of $h_{\cnum^{\ell}}$
compatible with $C_{\cnum^{\ell}}$ or $\omega_{\cnum^{\ell}}$.

\subsubsection{Basic polystable objects (2)}

Let $(\nbigp_{\ast}\nbigv,\theta)$ be 
a stable good filtered Higgs bundle of degree $0$ on $(X,D)$.
Assume that
$(\nbigp_{\ast}\nbigv,\theta)
\not\simeq
(\nbigp_{\ast}\nbigv^{\lor},\theta^{\lor})$.
We set
$\nbigp_{\ast}\nbigvtilde=
\nbigp_{\ast}\nbigv
\oplus
\nbigp_{\ast}(\nbigv^{\lor})$
which is equipped with the Higgs field
$\thetatilde=\theta\oplus\theta^{\lor}$.
We have the naturally defined perfect pairing
\[
 \nbigp_{\ast}\nbigv
 \otimes
 \nbigp_{\ast}(\nbigv^{\lor})
 \lrarr
 \nbigp^{(0)}_{\ast}(\nbigo_X(\ast D)).
\]
It induces a symmetric product
$\Ctilde_{(\nbigp_{\ast}\nbigv,\theta)}$ of
$(\nbigp_{\ast}\nbigvtilde,\thetatilde)$.
\begin{rem}
If
 $(\nbigp_{\ast}\nbigv,\theta)
 \simeq
 (\nbigp_{\ast}\nbigv^{\lor},\theta^{\lor})$,
$(\nbigp_{\ast}\nbigv,\theta)$
 has a symmetric or skew symmetric pairing $P_1$.
If $P_1$ is symmetric (skew-symmetric),
 $(\nbigp_{\ast}\nbigvtilde,\thetatilde,
 \Ctilde_{(\nbigp_{\ast}\nbigv,\theta)})$
is isomorphic to 
$(\nbigp_{\ast}\nbigv,\theta)\otimes\cnum^2$
with $P_1\otimes C_{\cnum^2}$
(resp. $P_1\otimes \omega_{\cnum^2}$).
\hfill\qed
\end{rem}

\begin{lem}
\label{lem;22.9.4.11}
Suppose that
$\bigl(
 (\nbigp_{\ast}\nbigv,\theta)\otimes\cnum^{\ell_1}
 \bigr)
 \oplus
 \bigl(
 (\nbigp_{\ast}\nbigv^{\lor},\theta^{\lor})\otimes\cnum^{\ell_2}
 \bigr)$ 
is equipped with
a perfect symmetric pairing $C$.
Then,
we have $\ell_1=\ell_2$,
and 
there exists 
an isomorphism
 $(\nbigp_{\ast}\nbigvtilde,\thetatilde)
 \otimes\cnum^{\ell_1}
 \simeq
 \bigl(
 (\nbigp_{\ast}\nbigv,\theta)\otimes\cnum^{\ell_1}
 \bigr)
 \oplus
 \bigl(
 (\nbigp_{\ast}\nbigv^{\lor},\theta^{\lor})\otimes\cnum^{\ell_2}
 \bigr)$
under which
 $\Ctilde_{(\nbigp_{\ast}\nbigv,\theta)}
 \otimes C_{\cnum^{\ell_1}}=C$.
\end{lem}
\pf
There exist $1$-dimensional subspaces
$L_i\subset\cnum^{\ell_i}$
such that
the restriction of $C$ to
$\Bigl(
(\nbigp_{\ast}\nbigv,\theta)\otimes L_1\Bigr)
 \oplus\Bigl(
 (\nbigp_{\ast}\nbigv^{\lor},\theta^{\lor})\otimes L_2
 \Bigr)$
is not $0$.
Because
$(\nbigp_{\ast}\nbigv,\theta)$
and 
$(\nbigp_{\ast}\nbigv^{\lor},\theta^{\lor})$
are stable,
and because 
$(\nbigp_{\ast}\nbigv,\theta)\not\simeq
(\nbigp_{\ast}\nbigv^{\lor},\theta^{\lor})$,
the restriction is equal to
$\alpha\cdot \Ctilde_{(\nbigp_{\ast}\nbigv,\theta)}$
for a non-zero $\alpha$.
In particular, it is perfect.
We obtain the decomposition of filtered bundles
which is orthogonal with respect to $C$:
\[
\Bigl(
 \nbigp_{\ast}(\nbigv)\otimes\cnum^{\ell_1}
 \Bigr)
 \oplus
 \Bigl(
 \nbigp_{\ast}(\nbigv^{\lor})\otimes\cnum^{\ell_2}
 \Bigr)
 =\Bigl(
\bigl(
 \nbigp_{\ast}(\nbigv)\otimes L_1
 \bigr)
  \oplus
  \bigl(
   \nbigp_{\ast}(\nbigv^{\lor})\otimes L_2
   \bigr)
   \Bigr)
   \oplus
\nbigp_{\ast}\nbigv'.
\]
It is preserved by the Higgs field,
and 
$\nbigp_{\ast}\nbigv'$ 
with the induced Higgs field
is isomorphic to
$(\nbigp_{\ast}\nbigv,\theta)\otimes\cnum^{\ell_1-1}
 \oplus
 (\nbigp_{\ast}\nbigv,\theta)\otimes\cnum^{\ell_2-1}$.
Hence, 
we obtain the claim of the lemma by an easy induction.
\hfill\qed

\vspace{.1in}
By using $C_{\cnum^{\ell}}$,
we identify $\cnum^{\ell}$
and the dual space $(\cnum^{\ell})^{\lor}$.
Then, the perfect symmetric bilinear form
$\Ctilde_{(\nbigp_{\ast}\nbigv,\theta)}\otimes C_{\cnum^{\ell}}$
on
\[
 (\nbigp_{\ast}\nbigvtilde,\thetatilde)\otimes\cnum^{\ell}
 =
 \bigl(
 (\nbigp_{\ast}\nbigv,\theta)\otimes\cnum^{\ell}
 \bigr)
 \oplus
 \bigl(
 (\nbigp_{\ast}\nbigv,\theta)^{\lor}\otimes(\cnum^{\ell})^{\lor}
 \bigr)
\]
is identified with the symmetric pairing
induced by the natural pairing
\[
 \bigl(
 (\nbigp_{\ast}\nbigv,\theta)\otimes\cnum^{\ell}
  \bigr)
  \otimes
  \bigl(
  (\nbigp_{\ast}\nbigv,\theta)^{\lor}\otimes(\cnum^{\ell})^{\lor}
  \bigr)
  \lrarr
  \nbigp^{(0)}_{\ast}(\nbigo_X(\ast D)).
\]

Let $h_0$ be any harmonic metric of
$(\nbigv,\theta)_{|X\setminus D}$
which is adapted to $\nbigp_{\ast}(\nbigv)$.
Note that 
for any harmonic metric $h'_0$ of
$(\nbigv,\theta)_{|X\setminus D}$
which is adapted to $\nbigp_{\ast}(\nbigv)$,
there exists $c>0$ such that $h'_0=ch_0$.
We obtain the induced harmonic metric
$h_0^{\lor}$ of
$(\nbigv^{\lor},\theta^{\lor})_{|X\setminus D}$
which is adapted to $\nbigp_{\ast}(\nbigv^{\lor})$.

\begin{lem}
\mbox{{}}\label{lem;22.9.4.22}
\begin{itemize}
 \item 
Let $h_{\cnum^{\ell}}$ be any Hermitian metric of $\cnum^{\ell}$.
Let $h^{\lor}_{\cnum^{\ell}}$
denote the induced Hermitian metric on $(\cnum^{\ell})^{\lor}$.
Then,
 $(h_0\otimes h_{\cnum^{\ell}})
 \oplus
 (h_0^{\lor}\otimes h^{\lor}_{\cnum^{\ell}})$
is a harmonic metric of 
$(\nbigvtilde,\thetatilde)_{|X\setminus D}\otimes\cnum^{\ell}$
such that
(i) it is adapted to $\nbigp_{\ast}\nbigvtilde$,
       (ii) it is compatible with
       $\Ctilde_{(\nbigp_{\ast}\nbigv,\theta)}\otimes C_{\cnum^{\ell}}$.
 \item
      Conversely,
      let $h$ be any harmonic metric
      of $(\nbigvtilde,\thetatilde)_{|X\setminus D}\otimes\cnum^{\ell}$
      satisfying the above conditions (i) and (ii).
      Then, there exists a Hermitian metric $h_{\cnum^{\ell}}$
      of $\cnum^{\ell}$
      such that
      $h=(h_0\otimes h_{\cnum^{\ell}})
      \oplus
      (h_0^{\lor}\otimes h^{\lor}_{\cnum^{\ell}})$.
\end{itemize} 
\end{lem}
\pf
The first claim is clear.
Let $h$ be as in the second claim.
By the uniqueness of harmonic metrics to a parabolic structure
(see Proposition \ref{prop;22.9.7.10}),
$\nbigv_{|X\setminus D}\otimes\cnum^{\ell}$
and
$\nbigv^{\lor}_{|X\setminus D}\otimes\cnum^{\ell}$
are orthogonal with respect to $h$.
Let $h_{\nbigv}$
and $h_{\nbigv^{\lor}}$
denote the restrictions of
$h$ to
$\nbigv_{|X\setminus D}\otimes\cnum^{\ell}$
and
$\nbigv^{\lor}_{|X\setminus D}\otimes\cnum^{\ell}$,
respectively.
By the uniqueness again,
there exists a Hermitian metric $h_{\cnum^{\ell}}$
of $\cnum^{\ell}$
such that
$h_{\nbigv}=h_0\otimes h_{\cnum^{\ell}}$.
There also exists a Hermitian metric
$h'_{(\cnum^{\ell})^{\lor}}$ of $(\cnum^{\ell})^{\lor}$
such that 
$h_{\nbigv^{\lor}}=h_0^{\lor}\otimes h'_{(\cnum^{\ell})^{\lor}}$.
By the compatibility condition,
we obtain that
$h'_{(\cnum^{\ell})^{\lor}}
=h_{\cnum^{\ell}}^{\lor}$.
\hfill\qed

\subsubsection{Polystable objects}

Let $(\nbigp_{\ast}\nbigv,\theta)$
be a polystable good filtered Higgs bundle of degree $0$ on $(X,D)$
with a perfect symmetric pairing $C$.

\begin{prop}
\label{prop;22.9.5.20}
There exist stable good filtered Higgs bundles
$(\nbigp_{\ast}\nbigv^{(0)}_i,\theta^{(0)}_i)$ $(i=1,\ldots,p(0))$,
$(\nbigp_{\ast}\nbigv^{(1)}_i,\theta^{(1)}_i)$ $(i=1,\ldots,p(1))$,
and 
$(\nbigp_{\ast}\nbigv^{(2)}_i,\theta^{(2)}_i)$ $(i=1,\ldots,p(2))$
of degree $0$ on $(X,D)$
such that the following holds.
\begin{itemize}
 \item
      $(\nbigp_{\ast}\nbigv^{(0)}_i,\theta^{(0)}_i)$
      is equipped with
      a symmetric perfect pairing $P^{(0)}_i$.
 \item $(\nbigp_{\ast}\nbigv^{(1)}_i,\theta^{(1)}_i)$
      is equipped with
      a skew-symmetric perfect pairing $P^{(1)}_i$.
 \item $(\nbigp_{\ast}\nbigv^{(2)}_i,\theta^{(2)}_i)
       \not\simeq
       (\nbigp_{\ast}\nbigv^{(2)}_i,\theta^{(2)}_i)^{\lor}$.
 \item There exist positive integers
       $\ell(a,i)$  and an isomorphism
\begin{multline}
 (\nbigp_{\ast}\nbigv,\theta)
 \simeq 
 \bigoplus_{i=1}^{p(0)}
 (\nbigp_{\ast}\nbigv^{(0)}_i,\theta^{(0)}_i)
 \otimes\cnum^{\ell(0,i)}
 \oplus
 \bigoplus_{i=1}^{p(1)}
 (\nbigp_{\ast}\nbigv^{(1)}_i,\theta^{(1)}_i)
 \otimes
 \cnum^{2\ell(1,i)}
 \oplus
 \\
 \bigoplus_{i=1}^{p(2)}
 \Bigl(
 \bigl(
 (\nbigp_{\ast}\nbigv^{(2)}_i,\theta^{(2)}_i)
 \otimes\cnum^{\ell(2,i)}
 \bigr)
 \oplus
 \bigl(
 (\nbigp_{\ast}\nbigv^{(2)}_i,\theta^{(2)}_i)^{\lor}
 \otimes(\cnum^{\ell(2,i)})^{\lor}
 \bigr)
 \Bigr).
\end{multline}
       Under this isomorphism,
       $C$ is identified with
       the direct sum of
       $P^{(0)}_1\otimes C_{\cnum^{\ell(0,i)}}$,
       $P^{(1)}_1\otimes\omega_{\cnum^{2\ell(1,i)}}$
       and
       $\Ctilde_{(\nbigp_{\ast}\nbigv^{(2)}_i,\theta^{(2)}_i)}
       \otimes C_{\cnum^{\ell(2,i)}}$.
 \item
      $(\nbigp_{\ast}\nbigv^{(a)}_i,\theta^{(a)}_i)
      \not\simeq
      (\nbigp_{\ast}\nbigv^{(a)}_j,\theta^{(a)}_j)$
      $(i\neq j)$ for $a=0,1,2$,
      and 
      $(\nbigp_{\ast}\nbigv^{(2)}_i,\theta^{(2)}_i)
      \not\simeq
      (\nbigp_{\ast}\nbigv^{(2)}_j,\theta^{(2)}_j)^{\lor}$
      for any $i,j$.
\end{itemize}
\end{prop}
\pf
It follows from 
Lemma \ref{lem;22.9.4.10}
and Lemma \ref{lem;22.9.4.11}.
\hfill\qed

\vspace{.1in}

Let $h^{(a)}_{i}$ $(a=0,1)$ be unique harmonic metrics of
$(\nbigv^{(a)}_i,\theta^{(a)}_i)_{|X\setminus D}$
such that
(i) $h^{(a)}_i$ are adapted to $\nbigp_{\ast}\nbigv^{(a)}_i$,
(ii)
$\Psi_{P^{(a)}_{i}}$ are isometric
with respect to
$h^{(a)}_i$ and $(h^{(a)}_i)^{\lor}$.
Let $h^{(2)}_i$ be harmonic metrics of
$(\nbigv^{(2)}_i,\theta^{(2)}_i)_{|X\setminus D}$
adapted to $\nbigp_{\ast}\nbigv^{(2)}_i$.

\begin{prop}
\label{prop;22.9.4.30}
There exists a harmonic metric $h$
of $(\nbigv,\theta)_{|X\setminus D}$
such that (i) it is adapted to $\nbigp_{\ast}\nbigv$,
(ii) it is compatible with $C$.
More precisely, the following holds.
\begin{itemize}
\item 
Let $h_{\cnum^{\ell(0,i)}}$
      be Hermitian metrics of $\cnum^{\ell(0,i)}$
      compatible with $C_{\cnum^{\ell(0,i)}}$.
Let     $h_{\cnum^{2\ell(1,i)}}$
      be Hermitian metrics of $\cnum^{2\ell(1,i)}$
      compatible with $\omega_{\cnum^{2\ell(0,i)}}$.
      Let $h_{\cnum^{\ell(2,i)}}$
      be any Hermitian metric of
      $\cnum^{\ell(2,i)}$.
      Then,
\begin{equation}
\label{eq;22.9.4.20}
      \bigoplus_{i=1}^{p(0)}
      \bigl(
      h^{(0)}_i\otimes
      h_{\cnum^{\ell(0,i)}}
      \bigr)
      \oplus
      \bigoplus_{i=1}^{p(1)}
      \bigl(
      h^{(1)}_i\otimes
      h_{\cnum^{2\ell(1,i)}}
      \bigr)
      \oplus
      \bigoplus_{i=1}^{p(2)}
      \Bigl(
      \bigl(
      h^{(2)}_i\otimes
      h_{\cnum^{\ell(2,i)}}
      \bigr)
      \oplus
     \bigl(
      (h^{(2)}_i)^{\lor}\otimes
      (h_{\cnum^{\ell(2,i)}})^{\lor}
      \bigr)
      \Bigr)
\end{equation}
      is a harmonic metric of
      $(\nbigv,\theta)_{|X\setminus D}$
      satisfying the conditions (i) and (ii).
 \item Conversely,
       if $h$ is a harmonic metric of
       $(\nbigv,\theta)_{|X\setminus D}$
       satisfying the conditions (i) and (ii),
       $h$ is of the form {\rm(\ref{eq;22.9.4.20})}.
\end{itemize}
\end{prop}
\pf
The first claim is obvious.
The second claim follows from Lemma \ref{lem;22.9.4.21}
and Lemma \ref{lem;22.9.4.22}.
\hfill\qed

\subsubsection{An equivalence}

Let $(E,\delbar_E,\theta,h)$ be a wild harmonic bundle on $(X,D)$
with a real structure $C$.
As explained in \S\ref{subsection;22.9.3.60},
we obtain 
a good filtered Higgs bundle
$(\nbigp_{\ast}^h(E),\theta)$ on $(X,D)$
equipped with a perfect symmetric pairing $C$.
It is polystable of degree $0$.

\begin{thm}
The above construction induces an equivalence
between the following objects.
\begin{itemize}
 \item Wild harmonic bundles on $(X,D)$
       equipped with a real structure.
 \item Polystable good filtered Higgs bundles
       of degree $0$ on $(X,D)$
       equipped with a perfect symmetric pairing.
\end{itemize} 
\end{thm}
\pf
We obtain the converse construction from
Proposition \ref{prop;22.9.4.30}.
\hfill\qed

\section{Generically regular semisimple case}

\subsection{Prolongation of regular semisimple wild Higgs bundles
on a punctured disc}

\subsubsection{Regular semisimple wild Higgs bundles}

Let $U$ be a neighbourhood of $0$ in $\cnum$.
For each $e\in\seisuu_{>0}$,
let $\varphi_e:\cnum\lrarr\cnum$ be given defined by
$\varphi_e(\zeta)=\zeta^e$.
We set $U^{(e)}=\varphi^{-1}(U)$
on which we set $z_e=\zeta$.
We have $z_e^e=\varphi_e^{\ast}(z)$.

Let $(E,\delbar_E,\theta)$ be
a Higgs bundle on $U\setminus\{0\}$ of rank $r$
equipped with a non-degenerate symmetric pairing $C$.
We obtain the endomorphism $f$
determined by $\theta=f\,dz/z$.
We obtain the characteristic polynomial
$\det(T\id_E-f)=T^{r}+\sum_{j=0}^{r-1} a_j(z)T^j$,
where $a_j(z)$ are holomorphic function on $U\setminus\{0\}$.
We assume the following condition.
\begin{itemize}
 \item $(E,\delbar_E,\theta)$ is regular semisimple,
       i.e.,
       for each $z_0\in U\setminus\{0\}$,
       the polynomial
       $T^{r}+\sum_{j=0}^{r-1} a_j(z_0)T^j$
       has $r$-distinct roots.
 \item $(E,\delbar_E,\theta)$ is wild,
       i.e.,
       $a_j(z)$ are meromorphic at $z=0$.
\end{itemize}

There exist a divisor $e$ of $r!$
and holomorphic functions $\alpha_1,\ldots,\alpha_r$ on
$U^{(e)}\setminus\{0\}$
such that
\[
 T^{r}
 +\sum_{j=0}^{r-1}
 \varphi_e^{\ast}(a_j)T^j
=\prod_{j=1}^{r}
 (T-\alpha_j(z_e)),
\]
and that $\alpha_i-\alpha_j$ $(i\neq j)$
are nowhere vanishing on $U^{(e)}\setminus\{0\}$.
Because $a_j(z)$ are meromorphic at $z=0$,
we obtain that $\alpha_j(z)$ are meromorphic at $z=0$.
We set $E^{(e)}:=\varphi_e^{\ast}(E)$,
which is equipped with the endomorphism
$f^{(e)}:=\varphi_e^{\ast}(f)$.
There exists the eigen decomposition
\[
 (E^{(e)},f^{(e)})
 =\bigoplus_{i=1}^r
 (E^{(e)}_{i},\alpha_i\id_{E^{(e)}_{i}}).
\]
It is orthogonal with respect to
$C^{(e)}:=\varphi_e^{\ast}(C)$.
The restriction of $C^{(e)}$
to $E^{(e)}_{i}$
are denoted by $C^{(e)}_{i}$.

Let $\Gal_e$ denote the Galois group of
the covering $U^{(e)}\setminus\{0\}\to U\setminus\{0\}$,
which is a cyclic group of order $e$.
The pull back
$\varphi^{\ast}(E,f,C)$
is naturally equivariant with respect to $\Gal_e$.
There exists the naturally induced $\Gal_e$-action
on $\{1,\ldots,r\}$
such that $b^{\ast}(\alpha_i)=\alpha_{b(i)}$
and
$b^{\ast}(E^{(e)}_{i})=E^{(e)}_{b(i)}$
for any $b\in\Gal_e$.

\subsubsection{Canonical meromorphic extension}

Let $\iota:U\setminus\{0\}\to U$
denote the inclusion.
There exists a natural inclusion
$\nbigo_U(\ast 0)\subset
\iota_{\ast}(\nbigo_{U\setminus\{0\}})$.
We naturally regard $E$ as an $\nbigo_{U\setminus \{0\}}$-module,
and we obtain the $\iota_{\ast}(\nbigo_{U\setminus\{0\}})$-module
$\iota_{\ast}(E)$.
We obtain the morphisms
$\iota_{\ast}(f):\iota_{\ast}(E)\to\iota_{\ast}(E)$
and 
$\iota_{\ast}(C):
\iota_{\ast}(E)
\otimes
\iota_{\ast}(E)
\lrarr
\iota_{\ast}(\nbigo_{U\setminus\{0\}})$.

For any locally free $\nbigo_{U\setminus\{0\}}$-module $\nbigf$,
a locally free $\nbigo_U(\ast 0)$-submodule
$\nbigf_1\subset \iota_{U\ast}\nbigf$
is called a meromorphic extension of $\nbigf$
if $\iota_{U}^{\ast}(\nbigf_1)=\nbigf$.

\begin{prop}
\label{prop;22.9.5.10}
There exists a unique meromorphic extension
$\Etilde$ of $E$
such that
(i) $\iota_{\ast}(f)(\Etilde)\subset\Etilde$,
 (ii) $\iota_{\ast}(C)(\Etilde\otimes\Etilde)
 \subset\nbigo_{U}(\ast 0)$.
 \end{prop}
\pf
Let $\iota_e:U^{(e)}\setminus\{0\}\to U^{(e)}$
denote the inclusion.
A meromorphic extension of
a locally free $\nbigo_{U^{(e)}\setminus\{0\}}$-module
$\nbigf$ 
is defined to be a locally free
$\nbigo_{U^{(e)}}(\ast 0)$-submodule
$\nbigf_1\subset\iota_{e\ast}(\nbigf)$
such that
$\iota_e^{\ast}(\nbigf_1)=\nbigf$.

\begin{lem}
\label{lem;22.9.5.2}
For each $i$,
there exists a unique meromorphic extension
$\Etilde^{(e)}_{i}$ of $E^{(e)}_{i}$
such that $C^{(e)}_{i}$ induces an isomorphism
\begin{equation}
\label{eq;22.9.5.1}
 \Etilde^{(e)}_{i}
 \otimes
  \Etilde^{(e)}_{i}
  \simeq
  \nbigo_{U^{(e)}}(\ast 0).
\end{equation}
\end{lem}
\pf
There exists a holomorphic frame $u^{(e)}_{i}$
of $E^{(e)}_{i}$ on $U^{(e)}\setminus\{0\}$.
We obtain the holomorphic function
$g=C^{(e)}_{i}(u^{(e)}_{i},u^{(e)}_{i})$.
There exist a holomorphic function $g_1$
such that
$g_1^2=g$ or $g_1^2=z_eg$.
We set $\utilde^{(e)}_{i}:=g_1^{-1}u^{(e)}_{i}$,
and
$\Etilde^{(e)}_{i}:=
\nbigo_{U^{(e)}}(\ast 0)\cdot \utilde^{(e)}_{i}
\subset
\iota_{e\ast}(E^{(e)}_{i})$.
Then, $\Etilde^{(e)}_{i}$
satisfies (\ref{eq;22.9.5.1}).

Let $\Etilde^{\prime(e)}_{i}$
be another meromorphic extension of
$E^{(e)}_{i}$
satisfying (\ref{eq;22.9.5.1}).
There exists a holomorphic frame $v^{(e)}_{i}$
of $\Etilde^{\prime(e)}_{i}$
on a neighbourhood $\nbigu$ of $0$ in $U^{(e)}$.
We obtain the holomorphic function
$g_2$ on $\nbigu\setminus\{0\}$
defined by
$v^{(e)}_{i}=g_2\cdot \utilde^{(e)}_{i}$.
Because both
$C^{(e)}_{i}
(v^{(e)}_{i}\otimes v^{(e)}_{i})$
and
$C^{(e)}_{i}
(\utilde^{(e)}_{i}\otimes \utilde^{(e)}_{i})$
are meromorphic at $z_e=0$,
we obtain that $g_2^2$ is meromorphic at $z_e=0$,
and hence $g_2$ is meromorphic at $z_e=0$.
It implies
$\Etilde^{\prime(e)}_{i}
=\Etilde^{(e)}_{i}$.
\hfill\qed

\vspace{.1in}
We obtain a meromorphic extension
$\Etilde^{(e)}:=\bigoplus \Etilde^{(e)}_{i}
\subset
\iota_{e\ast}(E^{(e)})$
of $E^{(e)}$.
It is naturally $\Gal_e$-equivariant.
By the construction,
$\iota_{e\ast}(f^{(e)})\bigl(
\Etilde^{(e)}
\bigr)\subset\Etilde^{(e)}$
and 
$\iota_{e\ast}(C^{(e)})\bigl(
\Etilde^{(e)}\otimes\Etilde^{(e)}
\bigr)\subset\nbigo_{U^{(e)}}(\ast 0)$.
As the descent
(see \cite[\S2.3.2]{Mochizuki-KH-Higgs}),
we obtain a locally free
$\nbigo_U(\ast 0)$-module $\Etilde$
with a $\Gal_e$-equivariant isomorphism
$\varphi_e^{\ast}(\Etilde)\simeq \Etilde^{(e)}$.
It is easy to see that
$\Etilde$ is a meromorphic extension of $E$
with the desired property.

Let $\Etilde'$ be another meromorphic extension
with the desired property.
We obtain a meromorphic extension
$\Etilde^{\prime(e)}:=
\varphi_e^{\ast}(\Etilde')$ of $E^{(e)}$
such that
$\iota_{e\ast}(f^{(e)})\bigl(
\Etilde^{\prime(e)}
\bigr)
\subset
\Etilde^{\prime(e)}$
and
$\iota_{e\ast}(C^{(e)})\bigl(
\Etilde^{\prime(e)}\otimes\Etilde^{\prime(e)}
\bigr)
\subset
\nbigo_{U^{(e)}}(\ast 0)$.
Let $f^{\prime\,(e)}$
and $C^{\prime(e)}$
denote the induced endomorphism
and the pairing of $\Etilde^{\prime(e)}$.

For a sheaf $\nbigf$ on $U^{(e)}$,
let $\nbigf_0$ denote the stalk of $\nbigf$ at $0$.
Note that
$\gbigk=\nbigo_{U^{(e)}}(\ast 0)_0$
is a field.
We obtain a $\gbigk$-vector space
$\Etilde^{\prime(e)}_0$
with the linear endomorphism
$f^{\prime(e)}_0$
and the symmetric bilinear pairing
$C^{\prime(e)}_0$.
Because the characteristic polynomial of
$f^{\prime(e)}_0$
is
$T^r+\sum_{j=0}^{r-1}\varphi_e^{\ast}(a_j)T^j\in\gbigk[T]$,
the eigenvalues of $f^{\prime(e)}_0$
are $\alpha_1,\ldots,\alpha_r\in\gbigk$,
and there exists the eigen decomposition
$\Etilde^{\prime(e)}_0
=\bigoplus_{i=1}^r
(\Etilde^{\prime(e)}_0)_{i}$,
where $f^{\prime(e)}_0-\alpha_i\id$ are $0$
on $(\Etilde^{\prime(e)}_0)_i$. 
Hence, there exists the decomposition
of the $\nbigo_{U}(\ast 0)$-module
$\Etilde^{\prime(e)}
 =\bigoplus
 \Etilde^{\prime(e)}_{i}$
such that
$\iota_{e\ast}(f^{(e)})-\alpha_i\id_{\Etilde^{\prime(e)}}$
are $0$ on $\Etilde^{\prime(e)}_i$.
Each $\Etilde^{\prime(e)}_{i}$
are meromorphic extension of
$E^{(e)}_{i}$,
and 
$C^{\prime(e)}$
induces an isomorphism
$\Etilde^{\prime(e)}_{i}
\otimes
\Etilde^{\prime(e)}_{i}
\simeq\nbigo_{U^{(e)}}(\ast 0)$.
By the uniqueness in Lemma \ref{lem;22.9.5.2},
we obtain
$\Etilde^{\prime(e)}_{i}
=\Etilde_{i}^{(e)}$,
and hence
$\Etilde^{\prime(e)}
=\Etilde^{(e)}$.
It implies $\Etilde'=\Etilde$.
\hfill\qed

\vspace{.1in}
The induced endomorphism $\Etilde\to\Etilde$
and the pairing
$\Etilde\otimes\Etilde\to \nbigo_U(\ast 0)$
are denoted by $\ftilde$ and $\Ctilde$, respectively.
Note that there exists the eigen decomposition
\begin{equation}
\label{eq;22.9.5.3}
 \varphi_e^{\ast}(\Etilde,\ftilde)
=\bigoplus_{i=1}^r
 (\Etilde_{i},\alpha_i\id_{\Etilde_{i}}). 
\end{equation}
It is orthogonal with respect to $\varphi_e^{\ast}\Ctilde$.

\subsubsection{Canonical filtered extension
and regular semisimplicity at $0$}
\label{subsection;22.9.6.20}

\begin{df}
\label{df;22.9.5.4}
A filtered bundle $\nbigp_{\ast}\Etilde$
is called a good filtered extension of
$(E,\delbar_E,\theta)$ with $C$
if the following condition is satisfied:
\begin{description}
 \item[(i)] $(\nbigp_{\ast}\Etilde,\theta)$ is
	    a good filtered Higgs bundle.
 \item[(ii)] $C$ is a perfect pairing of
 $\nbigp_{\ast}(\Etilde)$.
\hfill\qed
\end{description}
\end{df}

\begin{lem}
\label{lem;22.9.4.200}
There exists a unique filtered bundle
$\nbigp^{\can}_{\ast}(\Etilde)$ over $\Etilde$
satisfying the conditions (i), (ii)
and the following additional condition.
\begin{description}
 \item[(iii)] 
 $\varphi_e^{\ast}(\nbigp^{\can}_{\ast}\Etilde)$
is compatible with the decomposition
{\rm(\ref{eq;22.9.5.3})},
i.e., 
\begin{equation}
\label{eq;22.9.1.10}
 \varphi_e^{\ast}(\nbigp^{\can}_{\ast}\Etilde)
=\bigoplus
 \nbigp_{\ast}\Etilde^{(e)}_{i}.
\end{equation}
\end{description}
The filtered bundle $\nbigp^{\can}_{\ast}(\Etilde)$ is called
the canonical filtered extension of $(E,\delbar_E,\theta)$ with $C$.
\end{lem}
\pf
There exists a unique filtered bundle
$\nbigp_{\ast}\Etilde^{(e)}_{i}$
such that
$\nbigp_{\ast}\Etilde^{(e)}_{i}
\otimes
\nbigp_{\ast}\Etilde^{(e)}_{i}
\simeq
\nbigp_{\ast}^{(0)}(\nbigo_{U^{(e)}}(\ast 0))$.
We obtain the filtered bundle
$\nbigp_{\ast}(\Etilde^{(e)})$
by the right hand side of (\ref{eq;22.9.1.10}).
The uniqueness of such $\nbigp_{\ast}\Etilde^{(e)}_{i}$
implies that $\nbigp_{\ast}(\Etilde^{(e)})$
is $\Gal_e$-equivariant,
and we obtain a filtered bundle
$\nbigp_{\ast}(\Etilde)$ over $\Etilde$
as the descent of $\nbigp_{\ast}(\Etilde^{(e)})$,
which has the desired property.
The uniqueness is clear.
\hfill\qed

\begin{df}
\label{df;22.9.5.30}
We say that $(E,\delbar_E,\theta)$ is
regular semisimple at $0$
if the following condition is satisfied.
\begin{itemize}
 \item $(\alpha_i-\alpha_j)^{-1}$ $(i\neq j)$
       are holomorphic at $0$.
\hfill\qed
\end{itemize}
\end{df}

\begin{prop}
\label{prop;22.9.5.11}
If $(E,\delbar_E,\theta)$ is regular semisimple at $0$,
any good filtered extension of
$(E,\delbar_E,\theta)$ with $C$
is equal to $\nbigp^{\can}_{\ast}(\Etilde)$.
\end{prop}
\pf
If $(E,\delbar_E,\theta)$ is regular semisimple at $0$,
then
$\nbigp_{\ast}(\Etilde^{(e)})$
has to be compatible with the decomposition (\ref{eq;22.9.5.3}).
Hence, 
the proposition is clear.
\hfill\qed

\begin{rem}
If $(E,\delbar_E,\theta)$ is not regular semisimple at $0$,
there may exist many good filtered extensions of
$(E,\delbar_E,\theta)$ with $C$, in general.
(See Proposition {\rm\ref{prop;22.9.6.2}}.)
\hfill\qed
\end{rem}

\subsubsection{Compatible harmonic metrics}

Let $h\in\Harm(E,\delbar_E,\theta;C)$.
We obtain the good filtered Higgs bundle
$(\nbigp^h_{\ast}E,\theta)$ on $(U,0)$
with a perfect symmetric pairing $C$
as in Lemma \ref{lem;22.9.5.5}.
We also obtain
the locally free $\nbigo_U(\ast 0)$-module
$\nbigp^h(E)=\bigcup_{a\in\real}\nbigp^h_a(E)$.

\begin{prop}
\label{prop;22.9.5.12}
The $\nbigo_U(\ast 0)$-module
$\nbigp^h E$
equals to the canonical meromorphic extension
$\Etilde$.
If $(E,\delbar_E,\theta)$
is regular semisimple at $0$,
then we have
 $\nbigp^h_{\ast}(E)
 =\nbigp^{\can}_{\ast}(\Etilde)$.
\end{prop}
\pf
Because $C$ is a perfect pairing of
$\nbigp^h_{\ast}(E)$,
we obtain Proposition \ref{prop;22.9.5.12}
from the uniqueness in
Proposition \ref{prop;22.9.5.10}
and Proposition \ref{prop;22.9.5.11}.
\hfill\qed

\subsection{Classification of harmonic metrics by good filtered extensions}
\label{subsection;22.9.5.50}

\subsubsection{Setting}

Let $X$ be a compact Riemann surface
with a finite subset $D$.
Let $(E,\delbar_E,\theta)$ be a Higgs bundle on $X\setminus D$
with a non-degenerate symmetric pairing $C$.
We assume the following conditions.
\begin{itemize}
 \item The Higgs bundle is wild at each point of $D$.
 \item $(E,\delbar_E,\theta)$ is generically regular semisimple.
\end{itemize}

\begin{lem}
There exists a finite subset $Z_1\subset X\setminus D$
such that 
$(E,\delbar_E,\theta)_{|X\setminus (D\cup Z_1)}$ is regular semisimple.
\end{lem}
\pf
Let $P$ be any point of $D$.
Let $(X_P,z)$ be a holomorphic coordinate neighbourhood around $P$.
We obtain the endomorphism $f$ of $E_{|X_P\setminus\{P\}}$
by $\theta=f\,dz/z$.
We obtain the characteristic polynomial
$P_z(T)=T^r+\sum_{j=0}^{r-1} a_j(z)=\det(T\id_E-f)$.
Let $\Disc_T(P_z(T))$ be
the discriminant of $P_z(T)$,
which is a holomorphic function on $X_P\setminus\{P\}$.
Because $(E,\delbar_E,\theta)$ is generically regular semisimple,
there exists a discrete subset $Z_P\subset X_P\setminus\{P\}$
such that $\Disc_T(P_z(T))\neq 0$ unless $z\in Z_P$.
Because $a_j(z)$ are meromorphic at $z=0$,
$\Disc_T(P_z(T))$ is meromorphic at $z=0$.
Hence, we obtain the finiteness of $Z_P$.
\hfill\qed

\vspace{.1in}
Let $\iota:X\setminus D\to X$ denote the inclusion.
A meromorphic extension of 
a locally free $\nbigo_{X\setminus D}$-module $\nbigf$
is defined to
a locally free $\nbigo_{X}(\ast D)$-submodule
$\nbigf_1\subset\iota_{\ast}(\nbigf)$
such that $\iota^{\ast}(\nbigf_1)=\nbigf$.
We obtain the following proposition from
Proposition \ref{prop;22.9.5.10}.
\begin{prop}
There exists a unique meromorphic extension $\Etilde$
of $E$ such that
$\iota_{\ast}(\theta)(\Etilde)\subset
\Etilde\otimes\Omega^1_X$
and
 $\iota_{\ast}(C)(\Etilde\otimes\Etilde)
 \subset\nbigo_X(\ast D)$.
Such $\Etilde$ is called the canonical 
meromorphic extension of $(E,\delbar_E,\theta)$ with $C$.
\hfill\qed
\end{prop}

Let $\thetatilde:\Etilde\to\Etilde\otimes\Omega^1_X$
and $\Ctilde:\Etilde\otimes\Etilde\to\nbigo_X(\ast D)$
denote the induced morphisms.

\subsubsection{Good filtered extensions}
\label{subsection;22.9.23.10}

\begin{df}
A good filtered extension of
$(E,\delbar_E,\theta)$ with $C$
is a filtered bundle $\nbigp_{\ast}\Etilde$
over $\Etilde$
such that 
(i) $(\nbigp_{\ast}\Etilde,\thetatilde)$ is
a good filtered Higgs bundle of degree $0$,
(ii) $\Ctilde$ is a symmetric pairing of
 $\nbigp_{\ast}\Etilde$.
\hfill\qed
\end{df}
In this definition, we do not assume that
$(\nbigp_{\ast}\Etilde,\thetatilde)$ is polystable
nor that $\Ctilde$ is perfect.
Note that if $\Ctilde$ is perfect
the condition $\deg(\nbigp_{\ast}\Etilde)=0$ is automatically satisfied,
as remarked in Lemma \ref{lem;22.9.7.12}.

\begin{thm}
\label{thm;22.8.26.12}
The following holds.
\begin{itemize}
 \item $(\nbigp_{\ast}\Etilde,\thetatilde)$ is polystable,
       and $\Ctilde$ is perfect.
 \item Let $\Etilde_1\subset\Etilde$ be a saturated
       locally free $\nbigo_X(\ast D)$-submodule
       such that
       (i) $\thetatilde(\Etilde_1)\subset
       \Etilde_1\otimes\Omega^1_X$,
       (ii) $\deg(\nbigp_{\ast}\Etilde_1)=0$.
       Then, the restriction of $\Ctilde$
       to $\nbigp_{\ast}\Etilde_1$
       is also perfect,
       and we obtain the decomposition of
       good filtered Higgs bundles
\begin{equation}
\label{eq;22.8.26.3}
       (\nbigp_{\ast}\Etilde,\thetatilde)
       =(\nbigp_{\ast}\Etilde_1,\thetatilde_1)
       \oplus
       (\nbigp_{\ast}\Etilde_1^{\bot\,C},\thetatilde_1^{\bot\,C}).
\end{equation}
 \end{itemize} 
\end{thm}
\pf
Let $\Etilde_1\subset\Etilde$
be a saturated locally free $\nbigo_X(\ast D)$-submodule
such that $\theta(\Etilde_1)\subset\Etilde_1$.
We obtain the induced pairing
$\Ctilde_1:\nbigp_{\ast}(\Etilde_1)
\otimes\nbigp_{\ast}(\Etilde_1)\to
\nbigp_{\ast}^{(0)}\bigl(
\nbigo_X(\ast D)\bigr)$.

There exists a finite subset $Z_2\subset X\setminus D$
such that $(E,\delbar_E,\theta)_{|X\setminus (D\cup Z_2)}$
is regular semisimple,
and that $C_{|X\setminus (D\cup Z_2)}$ is non-degenerate.

\begin{lem}
\label{lem;22.8.26.1}
$C_{1|X\setminus (D\cup Z_2)}$ is a non-degenerate symmetric pairing
of $\Etilde_{1|X\setminus (D\cup Z_2)}$.
\end{lem}
\pf
Let $Q$ be any point of $X\setminus (D\cup Z_2)$.
Let $(X_Q,z)$ be a holomorphic coordinate neighbourhood around $Q$.
Let $f$ be the endomorphism of $\Etilde_{|X_Q}$
defined by $\theta=f\,dz$.
There exists the eigen decomposition
$E_{|Q}=\bigoplus_{\alpha\in\cnum} E_{Q,\alpha}$
of $f_{|Q}$.
The decomposition is orthogonal with respect to $C_{|Q}$.
Because $C$ is non-degenerate,
the restriction of $C$ to each $E_{Q,\alpha}$
is also non-degenerate.
Because $\Etilde_{1|Q}$ is a direct sum of
some eigen spaces, $\Ctilde_{1|Q}$ is non-degenerate.
\hfill\qed

\vspace{.1in}
We set $r_1:=\rank\Etilde_1$.
We obtain
$\det(\Etilde_1)\subset\bigwedge^{r_1}\Etilde$,
and the filtered bundle
$\nbigp_{\ast}\det(\Etilde_1)$.
We obtain the induced pairing
\begin{equation}
\label{eq;22.8.27.30}
\det(C_1):
 \nbigp_{\ast}\det(\Etilde_1)
 \otimes
 \nbigp_{\ast}\det(\Etilde_1)
 \lrarr
 \nbigp^{(0)}_{\ast}\bigl(
 \nbigo_X(\ast D)\bigr).
\end{equation}
It is non-zero by Lemma \ref{lem;22.8.26.1}.
Hence, by Proposition \ref{prop;22.8.28.3},
we obtain $\deg(\nbigp_{\ast}\Etilde_1)\leq 0$.
If moreover $\deg(\nbigp_{\ast}\Etilde_1)=0$ holds,
then $\Ctilde_1$ is perfect.
In particular, $\Ctilde$ is perfect.

If $\Etilde_1\neq \Etilde$,
we obtain the decomposition (\ref{eq;22.8.26.3})
by Lemma \ref{lem;22.8.27.32}.
By an easy induction,
we obtain the polystability of
$(\nbigp_{\ast}\Etilde,\thetatilde)$,
and Theorem \ref{thm;22.8.26.12} is proved.
\hfill\qed

\vspace{.1in}
There exists a filtered bundle $\nbigp^{\can}_{\ast}(\Etilde)$
over $\Etilde$
such that
for each $P\in D$
the restriction of $\nbigp^{\can}_{\ast}(\Etilde)$
to a neighbourhood of $P$
is equal to the canonical filtered extension in 
Lemma \ref{lem;22.9.4.200}.

\begin{lem}
We have $\deg(\nbigp^{\can}_{\ast}\Etilde)=0$.
Hence, $\nbigp^{\can}_{\ast}(\Etilde)$
is a good filtered extension of $(E,\delbar_E,\theta)$.
\end{lem}
\pf
By the construction, $\Ctilde$ induces
an isomorphism
$\det(\nbigp^{\can}_{\ast}\Etilde)
\otimes
\det(\nbigp^{\can}_{\ast}\Etilde)
\simeq
\nbigp^{(0)}_{\ast}(\nbigo_X(\ast D))$.
Hence, we obtain
$\deg(\nbigp^{\can}_{\ast}\Etilde)=0$.
\hfill\qed

\subsubsection{Classification and uniqueness}

\begin{thm}
\label{thm;22.9.5.40}
There exists the bijection between the following objects.
The correspondence is induced by $h\mapsto \nbigp^h_{\ast}(E)$
as in Proposition {\rm\ref{prop;22.9.5.12}}.
 \begin{itemize}
 \item Harmonic metrics of $(E,\delbar_E,\theta)$ compatible with $C$.
  \item Good filtered extensions 
	of $(E,\delbar_E,\theta)$ with $C$.
 \end{itemize}
\end{thm}
\pf
Let $\nbigp_{\ast}\Etilde$ be a good filtered extension of
$(E,\delbar_E,\theta)$ with $C$.
Note that
$(\nbigp_{\ast}\Etilde,\theta)$ is polystable of degree $0$.
There exists a decomposition
\[
 (\nbigp_{\ast}\Etilde,\theta)
 =\bigoplus_{i=1}^m
 \nbigp_{\ast}(\nbigv_i,\theta_i)
\]
into stable good filtered Higgs bundles of degree $0$.
Because $(E,\delbar_E,\theta)$ is generically regular semisimple,
we have $\nbigp_{\ast}(\nbigv_i,\theta_i)
\not\simeq
\nbigp_{\ast}(\nbigv_j,\theta_j)$
unless $i=j$.
Because the spectral curves of
$(\nbigv_i,\theta_i)$
and $(\nbigv_i^{\lor},\theta_i^{\lor})$
are the same,
we obtain that
$\nbigp_{\ast}(\nbigv^{\lor}_i,\theta^{\lor}_i)
\not\simeq
\nbigp_{\ast}(\nbigv_j,\theta_j)$
unless $i=j$.
By Proposition \ref{prop;22.9.5.20},
each $(\nbigp_{\ast}\nbigv_i,\theta_i)$
has a symmetric perfect pairing $P_i$,
and we may assume that
$\Ctilde$ is the direct sum of $P_i$.
By Proposition \ref{prop;22.9.4.30}
together with Lemma \ref{lem;22.9.4.21},
there exists a unique harmonic metric $h$
of $(E,\delbar_E,\theta)$ compatible with $C$.
\hfill\qed

\begin{thm}
\label{thm;22.9.5.100}
If moreover $(E,\delbar_E,\theta)$ is regular semisimple
at each point of $D$ in the sense of 
Definition {\rm\ref{df;22.9.5.30}},
there exists a unique harmonic metric $h$ of $(E,\delbar_E,\theta)$
compatible with $C$.
In this case,
$\nbigp^h_{\ast}(E)$  is equal to
$\nbigp^{\can}_{\ast}(\Etilde)$.
\end{thm}
\pf
It follows from Proposition \ref{prop;22.9.5.11}
and Theorem \ref{thm;22.9.5.40}.
\hfill\qed

\subsubsection{Complement}

Let $(\nbigp_{\ast}\nbigv,\theta)$
be a good filtered Higgs bundle of degree $0$
on $(X,D)$
with a symmetric pairing $C$.
Let $(E,\delbar_E,\theta)$ be the Higgs bundle
on $X\setminus D$
obtained as the restriction of $(\nbigv,\theta)$.
By applying the argument in the proof of
Theorem \ref{thm;22.8.26.12},
we obtain the following theorem.
\begin{thm}
Assume that there exists
$Q\in X\setminus D$
such that the following holds.
\begin{itemize}
 \item $C_{|Q}$ is non-degenerate,
       and $\theta$ is regular semisimple around $Q$.      
\end{itemize}
Then, the following holds.
\begin{itemize}
 \item $(\nbigp_{\ast}\nbigv,\theta)$ is polystable,
       and
       $C$ is perfect.
       In particular,
       $C_{|X\setminus D}$ is non-degenerate.
 \item $\nbigp_{\ast}\nbigv$
       is a good filtered extension of
       $(E,\delbar_E,\theta)$ with $C$.
\end{itemize}
As a result, there exists a unique harmonic metric $h$
of $(E,\delbar_E,\theta)$  compatible with $C$
such that $\nbigp^h_{\ast}E=\nbigp_{\ast}\nbigv$.
\hfill\qed
\end{thm}

\subsection{Examples}

Let $X$ be any Riemann surface.
We set $K_X=\Omega_X^1$.
For $r\in\seisuu_{>1}$,
we set
$\hyperk_{X,r}:=\bigoplus_{i=1}^r K_X^{(r-2i+1)/2}$.
Let $q_j$ $(j=2,\ldots,r)$ be
holomorphic $j$-differentials on $X$.
The multiplication of $q_j$ induces
\[
 K_X^{(r-2i+1)/2}\to
K_X^{(r-2i+2(j-1)+1)/2}\otimes K_X
 \quad
 (j\leq i\leq r).
\]
We also have the multiplications of $i(r-i)/2$
for $i=1,\ldots,r-1$:
\[
 K_X^{(r-2i+1)/2}
 \to
 K_X^{(r-2(i+1)+1)/2}\otimes K_X.
\]
They define a Higgs field
$\theta(\vecq)$ of $\hyperk_{X,r}$.
The natural pairing
$K_X^{(r-2i+1)/2}
\otimes
K_X^{-(r-2i+1)/2}
\to \nbigo_X$
induces a non-degenerate
symmetric bilinear form
$C_{\hyperk,X,r}$ of
$\hyperk_{X,r}$.
It is a symmetric pairing of
$(\hyperk_{X,r},\theta(\vecq))$.
We have the following corollary of
Theorem \ref{thm;22.8.31.10}.
\begin{cor}
\label{cor;22.9.6.102}
If the Higgs bundle $(\hyperk_{X,r},\theta(\vecq))$
is generically regular semisimple,
then there exists a harmonic metric
of $(\hyperk_{X,r},\theta(\vecq))$
which is compatible with $C_{\hyperk,X,r}$.
It is induced by an $\SL(r,\real)$-harmonic bundle.
\hfill\qed
\end{cor}
In some case, we can obtain
a precise classification by using
Theorem \ref{thm;22.9.5.40} and Theorem \ref{thm;22.9.5.100}.
Let us explain the case $r=3$ and $X=\cnum$.

\subsubsection{Example 1}
\label{subsection;22.9.6.110}

Let $\alpha_i\in\cnum[z]$ be polynomials
such that
$\alpha_1+\alpha_2+\alpha_3=0$
and that 
$\alpha_i-\alpha_j$ $(i\neq j)$ are not constantly $0$.
We set
\[
 q_2=
 -\frac{1}{2}(\alpha_1\alpha_2+\alpha_2\alpha_3+\alpha_3\alpha_1)(dz)^2,
 \quad
 q_3=\alpha_1\alpha_2\alpha_3(dz)^3.
\]
Let $f$ be the endomorphism of $\hyperk_{\cnum,3}$
defined by $\theta(q_2,q_3)=f\,dz$.
Because the eigenvalues of $f$ are $\alpha_i$ $(i=1,2,3)$,
the Higgs bundle is generically semisimple.

\begin{prop}
In this case,
$(\hyperk_{\cnum,3},\theta(q_2,q_3))$
is regular semisimple at $\infty$.
As a result,
it has a unique harmonic metric $h$ compatible with
$C_{\hyperk,\cnum,3}$. 
\end{prop}
\pf
We set $w=z^{-1}$.
Around $w=0$,
let $g$ be the endomorphism defined by
$g\,dw/w=\theta(q_2,q_3)$.
Because $g=-w^{-1}f$,
the eigenvalues of $g$ are
$-w^{-1}\alpha_i(w^{-1})$.
Hence, 
$(\hyperk_{\cnum,3},\theta(q_2,q_3))$
is regular semisimple at $\infty$.
We obtain the uniqueness of compatible harmonic metric
from Theorem \ref{thm;22.9.5.100}.
\hfill\qed

\vspace{.1in}
Because $\nbigp^h_{\ast}\hyperk_{\cnum,3}$
is equal to the canonical filtered extension,
we can study the asymptotic behaviour of $h$
more closely by using a general theory of wild harmonic bundles
\cite{Mochizuki-wild}.

\subsubsection{Example 2}
\label{subsection;22.9.6.112}

Let $\beta$ be a non-constant polynomial.
We set
\[
 q_2:=
 -\frac{1}{2}\Bigl(1-\frac{4}{3}\beta^2\Bigr)(dz)^2,
 \quad\quad
 q_3:=
 -\Bigl(
 \frac{2}{3}\beta-\frac{16}{27}\beta^3
 \Bigr)(dz)^3.
\]
On a neighbourhood of $\infty$,
we set
\[
 \alpha_1=\frac{1}{3}\beta+\beta(1-\beta^{-2})^{1/2},
 \quad
 \alpha_2=\frac{1}{3}\beta-\beta(1-\beta^{-2})^{1/2},
 \quad
 \alpha_3=-\frac{2}{3}\beta.
\]
Then, we have
\[
 T^3+\Bigl(1-\frac{4}{3}\beta^2\Bigr)T
+\frac{2}{3}\beta-\frac{16}{27}\beta^3
=(T-\alpha_1)(T-\alpha_2)(T-\alpha_3).
\]
Note that
\[
 \alpha_2-\alpha_3
 =\beta-\beta(1-\beta^{-2})^{1/2}
=\frac{1}{2}\beta^{-1}+O(\beta^{-2}).
\]
Set $w=z^{-1}$.
Let $g$ be the endomorphism of $\hyperk_{\cnum,3|\cnum\setminus\{0\}}$
defined by $\theta(q_2,q_2)=g(dw/w)$.
The eigenvalues of $g$ are
$-w^{-1}\alpha_i(w^{-1})$.

\begin{lem}
If $\deg(\beta)=1$,
then $(\hyperk_{\cnum,3},\theta(q_2,q_3))$
is regular semisimple at $\infty$.
If $\deg(\beta)\geq 2$,
$(\hyperk_{\cnum,3},\theta(q_2,q_3))$
is not regular semisimple at $\infty$.
\hfill\qed
\end{lem}

If $\deg(\beta)=1$,
$(\hyperk_{\cnum,3},\theta(q_2,q_3))$
has a unique harmonic metric
by Theorem \ref{thm;22.9.5.100}.
In the case $\deg(\beta)=2$,
we can classify good filtered extensions of 
$(\hyperk_{\cnum,3},\theta(q_2,q_3))$
by using Proposition \ref{prop;22.9.6.2} below,
which implies the classification of
harmonic metrics
of $(\hyperk_{\cnum,3},\theta(q_2,q_3))$
by Theorem \ref{thm;22.9.5.40}.
We set
\[
 \qtilde_2:=-\frac{1}{2}\Bigl(
  1-\frac{4}{3}\beta^2
  \Bigr),
  \quad
  v_i=\left(
  \begin{array}{c}
   \alpha_i^2-\qtilde_2 \\ \alpha_i \\ 1
  \end{array}
   \right).
\]
We have $g(v_i)=-w^{-1}\alpha_i\cdot v_i$.
We also have $C(v_2,v_2)=-2+2\beta^2-2\beta^2(1-\beta^{-2})^{1/2}$
and $C(v_3,v_3)=1$.
For any good filtered extension
$\nbigp_{\ast}\widetilde{\hyperk}_{\cnum,3}$,
there exists the decomposition of filtered bundles
\[
\nbigp_{\ast}\widetilde{\hyperk}_{\cnum,3}
=\nbigp_{\ast}(\nbigo_{U_{\infty}}(\ast\infty)v_1)
\oplus
\nbigp_{\ast}\Bigl(
\nbigo_{U_{\infty}}(\ast\infty)v_2
\oplus
\nbigo_{U_{\infty}}(\ast\infty)v_3
\Bigr).
\]
By an appropriate normalization,
$\nbigo_{U_{\infty}}(\ast\infty)v_2
\oplus\nbigo_{U_{\infty}}(\ast\infty)v_2$
with the induced Higgs field and the symmetric product
is isomorphic to
the Higgs bundle
$(\nbigv,\theta)$ with the symmetric product $C_{0,0}$
in \S\ref{subsection;22.9.6.21}.
Hence, there are good filtered extensions
of types (I), (II) and (III-1)
as in Proposition \ref{prop;22.9.6.2}.
In particular, the uniqueness
of harmonic metrics of $(\hyperk_{\cnum,3},\theta(q_2,q_3))$
does not hold in this case.

\subsubsection{Example 3}
\label{subsection;22.9.6.111}

For $a\in\cnum$,
we set $q_{2,a}=a z^2\sin(z)(dz)^2$ and $q_3=(z+1)^4\cos(z)(dz)^3$.
Because $(z^2\sin(z))_{|z=0}=0$
and $\bigl((z+1)^4\cos(z)\bigr)_{|z=0}=1$,
it is easy to check that
$(\hyperk_{\cnum,3},\theta(q_{2,a},q_3))$
is generically regular semisimple.
Hence, there exists a harmonic metric
of $(\hyperk_{\cnum,3},\theta(q_{2,a},q_3))$
compatible with $C_{\hyperk,\cnum,3}$.

\begin{rem}
If $a=0$, we can also prove the uniqueness of
such a harmonic metric by using a result in
{\rm\cite{Li-Mochizuki2}}.
\hfill\qed
\end{rem}

\subsection{Appendix: Classification of
regular filtered extensions in an easy case}
\label{subsection;22.9.6.21}

Let $U$ be a neighbourhood of $0$ in $\cnum$.
We set
$\nbigv=\nbigo_{U}(\ast \{0\})e_1
\oplus\nbigo_{U}(\ast \{0\})e_2$.
We consider the Higgs field
$\theta$ of $\nbigv$
given by
\[
 \theta e_1= e_1\,\,dz,
 \quad
 \theta e_2=e_2\,\,(-1)dz.
\]
For $\vecm=(m_1,m_2)\in\{0,1\}^2$,
let $C_{\vecm}$ be the symmetric pairing of $(\nbigv,\theta)$
determined by
$C_{\vecm}(e_i,e_i)=z^{m_i}$ $(i=1,2)$
and $C_{\vecm}(e_1,e_2)=0$.

\subsubsection{Logarithmic lattices}

For $\vecn=(n_1,n_2)\in\seisuu^2$,
we set
\[
 \nbigv^{(\vecn)}:=
 \nbigo_{U}z^{-n_1}e_1
 \oplus
 \nbigo_{U}z^{-n_2}e_2
 \subset\nbigv.
\]
For $\vecalpha=(\alpha_1,\alpha_2)\in\cnum^2\setminus\{(0,0)\}$
and $\vecn\in\seisuu^2$,
we set
\[
 \nbigv^{(\vecn,\vecalpha)}:=
 z\nbigv^{(\vecn)}
 +\nbigo_{U}\cdot
 \bigl(
 \alpha_1z^{-n_1}e_1+\alpha_2z^{-n_2}e_2
 \bigr)
 \subset\nbigv^{(\vecn)}.
\]
Note that
$\nbigv^{(\vecn,(\alpha_1,0))}=\nbigv^{(n_1,n_2-1)}$
and
$\nbigv^{(\vecn,(0,\alpha_2))}=\nbigv^{(n_1-1,n_2)}$.
We also note that
$\nbigv^{(\vecn,\vecalpha)}=
\nbigv^{(\vecn,\gamma\vecalpha)}$ for $\gamma\neq 0$.

A lattice of $\nbigv$ means
a locally free $\nbigo_U$-submodule $\nbigu\subset\nbigv$
such that $\nbigu(\ast 0)=\nbigv$.
It is called logarithmic if
$\theta(\nbigu)\subset\nbigu\otimes\Omega^1_U(\log 0)$.
If $\nbigu$ is a logarithmic lattice,
we obtain the endomorphism $\Res_{\nbigu}(\theta)$
of $\nbigu_{|0}=\nbigu/z\nbigu$.

Both $\nbigv^{(\vecn)}$ and $\nbigv^{(\vecn,\vecalpha)}$
are logarithmic.
We have
$\Res_{\nbigv^{(\vecn)}}(\theta)=0$.
If $\vecalpha\in(\cnum^{\ast})^2$,
we also have
$\Res_{\nbigv^{(\vecn,\vecalpha)}}(\theta)\neq 0$
and
$\Res_{\nbigv^{(\vecn,\vecalpha)}}(\theta)^2=0$.

\begin{lem}
Let $\nbigu\subset\nbigv$ be a logarithmic lattice.
Then, we have
$\nbigu=\nbigv^{(\vecn)}$ for some
$\vecn\in\seisuu^2$
or
$\nbigu=\nbigv^{(\vecn,\vecalpha)}$
 for some
 $(\vecn,\vecalpha)\in\seisuu^2\times
 \bigl((\cnum^2)\setminus\{(0,0)\}\bigr)$.
\end{lem}
\pf
Let $L_i$ $(i=1,2)$ be the image of
$\nbigu$ by the projection to $\nbigo_{U}(\ast 0)\cdot e_i$.
There exists $\vecn\in\seisuu^2$ such that
$L_i=\nbigo_U\cdot z^{-n_i}e_i$.
There exists $\vecalpha\in(\cnum^{\ast})^2$
and a holomorphic function $\beta$ on $U$
such that
(i) $\beta(0)\neq 0$,
(ii) $\alpha_1z^{-n_1}e_1+\alpha_2z^{-n_2}\beta e_2$
is contained in $\nbigu$.
Because $\theta$ is logarithmic with respect to $\nbigu$,
$\alpha_1z^{-n_1+1}e_1-\alpha_2z^{-n_2+1}\beta e_2$
is also contained in $\nbigu$.
We obtain $z\nbigv^{(\vecn)}\subset\nbigu$,
and hence
$\nbigv^{(\vecn,\vecalpha)}\subset\nbigu
 \subset
 \nbigv^{(\vecn)}$.
It implies
$\nbigu=
\nbigv^{(\vecn,\vecalpha)},
\nbigv^{(\vecn)}$.
\hfill\qed

\begin{lem}
\label{lem;22.9.6.1}
Let $\nbigu$ be a logarithmic lattice.
If $z\nbigv^{(\vecn)}\subsetneq\nbigu\subsetneq\nbigv^{(\vecn)}$
for some $\vecn$,
then there exists $\vecalpha\in\cnum^2\setminus\{(0,0)\}$
such that $\nbigu=\nbigv^{(\vecn,\vecalpha)}$.
If $z\nbigv^{(\vecn,\vecalpha)}
\subsetneq\nbigu
\subsetneq\nbigv^{(\vecn,\vecalpha)}$
for some $(\vecn,\vecalpha)\in\seisuu^2\times(\cnum^{\ast})^2$,
then we have
$\nbigu=z\nbigv^{(\vecn)}$. 
\end{lem}
\pf
The first claim is clear.
Let us study the second claim.
Let $L\subset \nbigv^{(\vecn,\vecalpha)}_{|0}$
be the image of $\nbigu$.
Because $\nbigu$ is a logarithmic lattice,
$L$ is preserved by $\Res_{\nbigv^{(\vecn,\vecalpha)}}(\theta)$.
Then, the second claim follows.
\hfill\qed

\subsubsection{Regular filtered Higgs bundles}

Let $(\vecn,\vecalpha)\in\seisuu^2\times(\cnum^2\setminus\{(0,0)\})$.
For $b\in\real$,
let $\nbigp^{(\vecn;b)}_{\ast}(\nbigv)$
and $\nbigp^{(\vecn,\vecalpha;b)}_{\ast}(\nbigv)$
be the filtered bundles over $\nbigv$
defined as follows $(c\in\real)$:
\[
 \nbigp^{(\vecn;b)}_{c}(\nbigv)
=z^{-[c-b]}\nbigv^{(\vecn)},
\quad\quad
 \nbigp^{(\vecn,\vecalpha;b)}_{c}(\nbigv)
=z^{-[c-b]}\nbigv^{(\vecn,\vecalpha)}.
\]
Here, $[d]:=\max\{n\in\seisuu\,|\,n\leq d\}$.
For $\vecb=(b_1,b_2)\in\real^2$
such that $b_1-1<b_2<b_1$,
let $\nbigp^{(\vecn,\vecalpha;\vecb)}_{\ast}(\nbigv)$
be the filtered bundle defined as follows
($c\in\real$, $m\in\seisuu$):
\[
 \nbigp^{(\vecn,\vecalpha;\vecb)}_{c}(\nbigv)
 =\left\{
 \begin{array}{ll}
  z^{-m}\nbigv^{(\vecn,\vecalpha)}&
 (m+b_2\leq c<m+b_1)\\
  \\
  z^{-m}\nbigv^{(\vecn)}&
 (m+b_1\leq c<m+1+b_2).
 \end{array}
\right.
\]
We set $\vecdelta=(1,1)$.
We can check the following lemma directly from the definitions.
\begin{lem}
\label{lem;22.9.6.10}
The following holds for any $m\in\seisuu$.
\[
\nbigp^{(\vecn;b)}_{\ast}(\nbigv)
=\nbigp^{(\vecn-m\vecdelta;b-m)}_{\ast}(\nbigv),
 \quad
\nbigp^{(\vecn,\vecalpha;b)}_{\ast}(\nbigv)
 =\nbigp^{(\vecn-m\vecdelta,\vecalpha;b-m)}_{\ast}(\nbigv),
\quad 
\nbigp^{(\vecn,\vecalpha;\vecb)}_{\ast}(\nbigv)
=\nbigp^{(\vecn-m\vecdelta,\vecalpha;\vecb-m\vecdelta)}.
\]
The following holds:
\[
\nbigp^{((n_1,n_2),(\alpha_1,0);b)}_{\ast}(\nbigv)
 =\nbigp^{((n_1,n_2-1);b)}_{\ast}(\nbigv),
 \quad
 \nbigp^{((n_1,n_2),(0,\alpha_2);b)}_{\ast}(\nbigv)
 =\nbigp^{((n_1-1,n_2);b)}_{\ast}(\nbigv).
\]
The following holds:
\[
 \nbigp^{((n_1,n_2),(\alpha_1,0);(b_1,b_2))}_{\ast}(\nbigv)
=\nbigp^{((n_1,n_2-1),(0,\alpha_2);(b_2,b_1-1))}_{\ast}(\nbigv).
\]
We also have
$\nbigp^{(\vecn,\vecalpha;b)}_{\ast}(\nbigv)
=\nbigp^{(\vecn,\gamma\vecalpha;b)}_{\ast}(\nbigv)$
and  
$\nbigp^{(\vecn,\vecalpha;\vecb)}_{\ast}(\nbigv)
=\nbigp^{(\vecn,\gamma\vecalpha;\vecb)}_{\ast}(\nbigv)$
for any $\gamma\in\cnum^{\ast}$.
\hfill\qed 
\end{lem}

\begin{prop}
Let $\nbigp_{\ast}(\nbigv)$ be a filtered bundle
such that $(\nbigp_{\ast}\nbigv,\theta)$
is a regular filtered Higgs bundle.
Then, $\nbigp_{\ast}(\nbigv)$
equals one of
$\nbigp^{(\vecn;b)}_{\ast}(\nbigv)$,
$\nbigp^{(\vecn,\vecalpha;b)}_{\ast}(\nbigv)$
or
$\nbigp^{(\vecn,\vecalpha;\vecb)}_{\ast}(\nbigv)$. 
\end{prop}
\pf
If there exists $-1<b\leq 0$
such that
$\Gr^{\nbigp}_c(\nbigv)=0$
unless $c-b\in\seisuu$,
$\nbigp_{\ast}(\nbigv)$
equals either
$\nbigp^{(\vecn;b)}_{\ast}(\nbigv)$
or $\nbigp^{(\vecn,\vecalpha;b)}_{\ast}(\nbigv)$.
If not,
there exists $\vecb\in\real^2$
such that (i) $b_1-1<b_2<b_1$,
(ii)
$\Gr^{\nbigp}_c(\nbigv)=0$
if $c-b_i\not\in\seisuu$ $(i=1,2)$.
By Lemma \ref{lem;22.9.6.1},
one of $\nbigp_{b_i}(\nbigv)$
is of the form $\nbigv^{(\vecn)}$ for some $\vecn\in\seisuu^2$.
By exchanging $(b_1,b_2)$ with $(b_2,b_1-1)$ if necessary,
we may assume that
$\nbigp_{b_1}(\nbigv)=\nbigv^{(\vecn)}$.
We obtain $\nbigp_{b_2}(\vecv)=\nbigv^{(\vecn,\vecalpha)}$
for some $\vecalpha\in \cnum^2\setminus\{(0,0)\}$.
Then, we obtain
$\nbigp_{\ast}(\nbigv)
=\nbigp^{(\vecn,\vecalpha;\vecb)}_{\ast}(\nbigv)$.
\hfill\qed

\subsubsection{Compatibility with the symmetric form}

\begin{prop}
\label{prop;22.9.6.2}
Let $\vecn\in\seisuu^2$ and $b\in\real$.
\begin{description}
 \item[(I)] $\nbigp^{(\vecn;b)}_{\ast}(\nbigv)$
       is compatible with $C_{\vecm}$
       if and only if
       $n_1-\frac{m_1}{2}=n_2-\frac{m_2}{2}=b$.
 \item[(II)] If $\vecalpha\in(\cnum^{\ast})^2$,
	    $\nbigp^{(\vecn,\vecalpha;b)}_{\ast}(\nbigv)$
       is compatible with $C_{\vecm}$
       if and only if
       $n_1-\frac{1}{2}m_1=n_2-\frac{1}{2}m_2=b+\frac{1}{2}$
       and
       $\alpha_1^2+\alpha_2^2=0$.
 \item[(III-1)]
	    If $\vecalpha\in(\cnum^{\ast})^2$,
	    $\nbigp^{(\vecn,\vecalpha;\vecb)}_{\ast}(\nbigv)$
       is compatible with $C_{\vecm}$
       if and only if
       $n_1-\frac{1}{2}m_1=n_2-\frac{1}{2}m_2=\frac{1}{2}(b_1+b_2)$
       and
       $\alpha_1^2+\alpha_2^2=0$.
 \item[(III-2)]
       $\nbigp^{(\vecn,(0,\alpha_2);\vecb)}_{\ast}(\nbigv)$
       is compatible with $C_{\vecm}$
	   if and only if
	   $n_1-\frac{1}{2}m_1=b_1$ and
	   $n_2-\frac{1}{2}m_2=b_2$.
\end{description}
\end{prop}
\pf
As a preliminary, let us study the dual lattices.
Let $e_1^{\lor},e_2^{\lor}$
be the dual frame of $\nbigv^{\lor}$.
For $\vecn\in\seisuu^2$
and $\vecalpha\in\cnum^2\setminus\{(0,0)\}$,
we set
\[
 (\nbigv^{\lor})^{(\vecn)}
=\nbigo_U\cdot z^{-n_1}e_1^{\lor}
\oplus 
 \nbigo_U\cdot z^{-n_2}e_2^{\lor},
 \quad\quad
 (\nbigv^{\lor})^{(\vecn,\vecalpha)}
=z(\nbigv^{\lor})^{(\vecn)}
+\nbigo_{U}\cdot
\bigl(
 \alpha_1z^{-n_1}e_1^{\lor}
+\alpha_2z^{-n_2}e_2^{\lor}
\bigr).
\]
For any lattice $\nbigu$ of $\nbigv$,
we set
$\nbigu^{\lor}=\nhom_{\nbigo_U}(\nbigu,\nbigo_U)$.
We have
$(\nbigv^{(\vecn)})^{\lor}
 =(\nbigv^{\lor})^{(-\vecn)}$.
\begin{lem}
\label{lem;22.9.6.3}
For $\vecalpha\in\cnum^2\setminus\{(0,0)\}$,
let $\vecbeta\in\cnum^2\setminus\{(0,0)\}$
such that $\alpha_1\beta_1+\alpha_2\beta_2=0$.
 Then, we have
\[
(\nbigv^{(\vecn,\vecalpha)})^{\lor}
=(\nbigv^{\lor})^{(-\vecn+\vecdelta,\vecbeta)}.
\] 
\end{lem}
\pf
Because
$\nbigv^{(\vecn-\vecdelta)}
\subset
\nbigv^{(\vecn,\vecalpha)}$,
we obtain
$(\nbigv^{(\vecn,\vecalpha)})^{\lor}
\subset
(\nbigv^{(\vecn-\vecdelta)})^{\lor}
=(\nbigv^{\lor})^{(-\vecn+\vecdelta)}$.
Because
\[
\bigl\langle
 \beta_1z^{n_1-1}e_1^{\lor}+\beta_2z^{n_2-1}e_2^{\lor},
 \alpha_1z^{-n_1}e_1+\alpha_2z^{-n_2}e_2
 \bigr\rangle
=z^{-1}\bigl(\alpha_1\beta_1+\alpha_2\beta_2\bigr),
\]
$\beta_1z^{n_1-1}e_1^{\lor}+\beta_2z^{n_2-1}e_2^{\lor}$
is contained in
$(\nbigv^{(\vecn,\vecalpha)})^{\lor}$
if and only if
$\beta_1\alpha_1+\beta_2\alpha_2=0$.
\hfill\qed

\vspace{.1in}

Let us return to the proof of Proposition \ref{prop;22.9.6.2}.
Because
$\Psi_{C_{\vecm}}(\nbigv^{(\vecn)})
=(\nbigv^{\lor})^{(\vecn-\vecm)}$,
$\nbigp^{(\vecn;b)}_{\ast}(\nbigv)$
is compatible with $C_{\vecm}$
if and only if $2n_1-m_1=2n_2-m_2=2b$,
i.e.,
$n_1-\frac{1}{2}m_1=n_2-\frac{1}{2}m_2=b$.

Let us consider the case of
$\nbigp_{\ast}^{(\vecn,\vecalpha;b)}(\nbigv)$
for $\vecalpha\in(\cnum^{\ast})^2$.
Note that
$\Psi_{C_{\vecm}}(\nbigv^{(\vecn,\vecalpha)})
=(\nbigv^{\lor})^{(\vecn-\vecm,\vecalpha)}$.
Because it cannot be
$(\nbigv^{\lor})^{(\vecp)}$ for any $\vecp\in\seisuu^2$,
it equals to
$(\nbigv^{\lor})^{(-\vecn+\ell\vecdelta,\vecbeta)}$
for some $\ell\in\seisuu$ and
$\vecbeta\in(\cnum^{\ast})^2$ such that
$\alpha_1\beta_1+\alpha_2\beta_2=0$
by Lemma \ref{lem;22.9.6.3}.
Hence,
we obtain $\alpha_1^2+\alpha_2^2=0$
and $n_1-\frac{1}{2}m_1=n_2-\frac{1}{2}m_2$.
Note the following equalities:
\[
 \bigl\langle
  \alpha_1z^{-n_1}e_1+\alpha_2z^{-n_2}e_2,\,\,
  \alpha_1z^{-n_1+m_1}e_1^{\lor}+\alpha_2z^{-n_2+m_2}e_2^{\lor}
 \bigr\rangle
=0,
\]
\[
  \bigl\langle
  \alpha_1z^{-n_1}e_1+\alpha_2z^{-n_2}e_2,\,\,
  \alpha_1z^{-n_1+1+m_1}e_1^{\lor}-\alpha_2z^{-n_2+1+m_2}e_2^{\lor}
  \bigr\rangle
=z^{-2n_1+1+m_1}(\alpha_1^2-\alpha_2^2),
\]
\[
 \bigl\langle
  \alpha_1z^{-n_1+1}e_1-\alpha_2z^{-n_2+1}e_2,\,\,
  \alpha_1z^{-n_1+m_1}e_1^{\lor}+\alpha_2z^{-n_2+m_2}e_2^{\lor}
 \bigr\rangle
  =z^{-2n_1+1+m_1}(\alpha_1^2-\alpha_2^2),
\]
\[
  \bigl\langle
  \alpha_1z^{-n_1+1}e_1-\alpha_2z^{-n_2+1}e_2,\,\,
  \alpha_1z^{-n_1+1+m_1}e_1^{\lor}-\alpha_2z^{-n_2+1+m_2}e_2^{\lor}
  \bigr\rangle
=0.
\]
Hence,
under the condition $n_1-\frac{1}{2}m_1=n_2-\frac{1}{2}m_2$
and $\alpha_1^2+\alpha_2^2=0$,
$\nbigp^{(\vecn,\alpha;b)}_{\ast}(\nbigv)$
is compatible with $C_{\vecm}$
if and only if $2b=2n_1-1-m_1$.
Similarly,
if $\vecalpha\in(\cnum^{\ast})^2$,
and if $\nbigp^{(\vecn,\vecalpha;\vecb)}_{\ast}(\nbigv)$
is compatible with $C_{\vecm}$,
we obtain $n_1-\frac{1}{2}m_1=n_2-\frac{1}{2}m_2$
and $\alpha_1^2+\alpha_2^2=0$,
and under these conditions,
the compatibility condition is
$(b_1-1)+b_2=2n_1-1-m_1$,
i.e.,
$n_1-\frac{1}{2}m_1=\frac{1}{2}(b_1+b_2)$.
If $\alpha_1=0$,
then the filtration
$\nbigp^{(\vecn,\vecalpha;\vecb)}_{\ast}(\nbigv)$
is compatible with the decomposition
$\nbigv=\nbigo_U(\ast 0)e_1\oplus\nbigo_U(\ast 0)e_2$.
Hence,
it is compatible with $C_{\vecm}$
if and only if it equals to the canonical filtered extension
in \S\ref{subsection;22.9.6.20}
by Lemma \ref{lem;22.9.4.200}.
It is equivalent to
$n_i-\frac{1}{2}m_i=b_i$.
\hfill\qed

\vspace{.1in}
Suppose $m_1=m_2=:m$.
Then, the filtered bundle
$\nbigp^{((0,0);-m/2)}_{\ast}(\nbigv)$
of type I is unique
by Lemma {\rm\ref{lem;22.9.6.10}},
which equals the canonical filtered extension.
There are two filtered bundles
$\nbigp_{\ast}^{((0,0),\vecalpha;-(1+m)/2)}(\nbigv)$
of type II corresponding to
$\vecalpha=(1,\pm\sqrt{-1})$.
The filtered bundles
$\nbigp^{((0,0),\vecalpha;(b_1,-b_1-m))}_{\ast}(\nbigv)$
of type (III-1)
are parameterized by
$\vecalpha=(1,\pm\sqrt{-1})$
and
$-\frac{m+1}{2}<b_1<-\frac{m}{2}$,
where we set $\vecn=(0,0)$.
The filtered bundle of type (III-2) does not exist.

Suppose $m_1\neq m_2$.
Then, there does not exist
the filtered bundle of type I, II nor (III-1).
There is a unique filtered bundle of type (III-2),
which is equal to the canonical filtered extension.

\paragraph{Acknowledgements}

The authors are grateful to the reviewer
for his/her careful reading and valuable comments.

T.M. is grateful to
Michael McBreen,
Natsuo Miyatake,
Franz Pedit
Martin Traizet
and Hitoshi Fujioka
for interesting discussions.
T.M. is partially supported by
the Grant-in-Aid for Scientific Research (A) (No. 21H04429),
the Grant-in-Aid for Scientific Research (A) (No. 22H00094),
the Grant-in-Aid for Scientific Research (A) (No. 23H00083),
and the Grant-in-Aid for Scientific Research (C) (No. 20K03609),
Japan Society for the Promotion of Science.
T.M. is also partially supported by the Research Institute for Mathematical
Sciences, an International Joint Usage/Research Center located in Kyoto
University.

Q.L. is partially supported by the National Key R\&D Program of China No. 2022YFA1006600, the Fundamental Research
Funds for the Central Universities and Nankai Zhide foundation.

\end{document}